\documentclass{amsart}
\usepackage{amsfonts}
\usepackage{amssymb, amsmath, amsthm}
\usepackage{amscd}
\usepackage{verbatim}
\usepackage{enumerate}
\usepackage{setspace}

\allowdisplaybreaks
\newtheorem{theorem}{Theorem}[section]

\newtheorem{corollary}[theorem]{Corollary}

\newtheorem{definition}[theorem]{Definition}

\newtheorem{lemma}[theorem]{Lemma}

\newtheorem{proposition}[theorem]{Proposition}

\theoremstyle{remark}
\newtheorem{remark}[theorem]{Remark}
\newtheorem{example}[theorem]{Example}
\numberwithin{equation}{section}

\newcommand{\beq}{\begin{equation}}
\newcommand{\eeq}{\end{equation}}
\newcommand{\beqa}{\begin{eqnarray}}
\newcommand{\eeqa}{\end{eqnarray}}

\input{tcilatex}

\begin{document}
\title[Regularized family of two-phase flows]{On a regularized family of
models for homogeneous incompressible two-phase flows}
\author{Ciprian G. Gal}
\address{Department of Mathematics, Florida International University, Miami,
FL 33199, USA}
\email{cgal@fiu.edu}
\author{T. Tachim Medjo}
\address{Department of Mathematics, Florida International University, Miami,
FL 33199, USA}
\email{tachimt@fiu.edu}

\begin{abstract}
We consider a general family of regularized models for incompressible
two-phase flows based on the Allen-Cahn formulation in $n$-dimensional
compact Riemannian manifolds for $n=2,3$. The system we consider consists of
a regularized family of Navier-Stokes equations (including the Navier-Stokes-%
$\alpha $-like model, the Leray-$\alpha $ model, the Modified Leray-$\alpha $
model, the Simplified Bardina model, the Navier-Stokes-Voight model and the
Navier-Stokes model) for the fluid velocity $u$ suitably coupled with a
convective Allen-Cahn equation for the order (phase) parameter $\phi $. We
give a unified analysis of the entire three-parameter family of two-phase
models using only abstract mapping properties of the principal dissipation
and smoothing operators, and then use assumptions about the specific form of
the parameterizations, leading to specific models, only when necessary to
obtain the sharpest results. We establish existence, stability and
regularity results, and some results for singular perturbations, which as
special cases include the inviscid limit of viscous models and the $\alpha
\rightarrow 0$ limit in $\alpha $-models. Then, we also show the existence
of a global attractor and exponential attractor for our general model, and
then establish precise conditions under which each trajectory $\left( u,\phi
\right) $ converges to a single equilibrium by means of a Lojasiewicz-Simon
inequality. We also derive new results on the existence of global and
exponential attractors for the regularized family of Navier--Stokes
equations and magnetohydrodynamics models which improve and complement the
results of \cite{HLT}. Finally, our analysis is applied to certain
regularized Ericksen-Leslie (RSEL) models for the hydrodynamics of liquid
crystals in $n$-dimensional compact Riemannian manifolds.
\end{abstract}

\maketitle
\tableofcontents

\section{Introduction}

\label{s:intro}

Modelling and simulating the behavior of binary fluid mixtures in various
turbulent regimes can be rather challenging \cite{AM}. A possible approach
is based on the so-called diffuse-interface method (see \cite{AM,Br,Mo} and
their references). This method consists in introducing an order parameter,
accounting for the presence of two species, whose dynamics interacts with
the fluid velocity. For incompressible fluids with matched densities, a
well-known model consists of the classical Navier-Stokes equation suitably
coupled with either a convective Cahn-Hilliard or Allen-Cahn equation (see 
\cite{Bl, DLRW, DLL,FHL, GPV, HH, RN, TLK,YFLS} cf. also \cite{BCB, CV, JV,
LM, LT, O, Si}). Denoting by $u=\left( u_{1},...,u_{n}\right) ,$ $n\geq 2$,
the velocity field and by $\phi $ the order parameter, where we suppose that 
$\phi $ is normalized in such a way that the two pure phases of the fluid
are $-1$ and $+1,$ respectively, the Cahn-Hiliard-Navier-Stokes and the
Allen-Cahn-Navier-Stokes systems can be written in a unified form. Indeed,
if additionally we assume that the viscosity of fluid is constant, and
temperature differences are negligible, we have 
\begin{align}
& \partial _{t}u+u\cdot \nabla u-\nu \Delta u+\nabla p=-\varepsilon \text{div%
}(\nabla \phi \otimes \nabla \phi )+g,  \label{1.1} \\
& \text{div}\left( u\right) =0,  \label{1.2} \\
& \partial _{t}\phi +u\cdot \nabla \phi +A_{K}\mu =0,  \label{1.3} \\
& \mu =-\varepsilon \Delta \phi +\varepsilon ^{-1}f\left( \phi \right) ,
\label{1.4}
\end{align}%
in $\Omega \times \left( 0,+\infty \right) ,$ where $\Omega $ is a bounded
domain in $\mathbb{R}^{n},$ $n=2,3$, with a sufficiently smooth boundary $%
\Gamma ,$ $\varepsilon >0$ is a parameter related to the thickness of the
interface separating the two fluids, and $g=g\left( t\right) $ is an
external body force. Moreover, the operator $A_{K}$ has a two-fold
definition according to the case $K=CH$ (Cahn-Hilliard fluid) or $K=AC$
(Allen-Cahn fluid), namely, 
\begin{equation}
A_{CH}\mu =-m\Delta \mu ,\qquad A_{AC}\mu =\mu ,  \label{CHAC}
\end{equation}%
where $m>0$ is the mobility of the mixture. The so-called chemical potential 
$\mu $ is obtained, under an appropriate choice of boundary conditions, as a
variational derivative of the following free energy functional 
\begin{equation}
\mathcal{F}\left( \phi \right) =\int\limits_{\Omega }\left( \frac{%
\varepsilon }{2}\left\vert \nabla \phi \right\vert ^{2}+\varepsilon
^{-1}F\left( \phi \right) \right) dx,  \label{1.5}
\end{equation}%
where $F(r)=\int_{0}^{r}f(y)dy$, $r\in \mathbb{R}$. Here, the potential $F$
is either a double-well logarithmic-type function%
\begin{equation}
F\left( r\right) =\gamma _{1}\left( \left( 1+r\right) \log \left( 1+r\right)
+\left( 1-r\right) \log \left( 1-r\right) \right) +\gamma _{2}\left(
1-r^{2}\right) ,\quad r\in (-1,1),  \label{sing-pot}
\end{equation}%
for some $\gamma _{1},\gamma _{2}>0$, or a polynomial approximation of the
type%
\begin{equation}
F\left( r\right) =\gamma _{3}\left( r^{2}-1\right) ^{2},  \label{reg-pot}
\end{equation}%
for some $\gamma _{3}>0.$

Both the two systems (\ref{1.1})-(\ref{1.4}) for $K=AC$ and $K=CH$ capture
basic features of binary fluid behavior. There are several key differences
between the two formulations:

\textbf{(I)} With the Allen-Cahn formulation, both singular and regular
potentials can be treated since a maximum principle holds under appropriate
assumptions on $F$ when $K=AC.$ In this case, we recall that the phase-field 
$\phi $ takes values in a given bounded interval (i.e., the domain of $F$ in
the singular potential case (\ref{sing-pot})) or $\phi \in \left[ -1,1\right]
$ in the regular potential case, when $F$ is of the form (\ref{reg-pot}),
see \cite[Section 6]{GG1}. On the other hand, with the Cahn-Hilliard
formulation the latter property is only true when $F$ is a singular
potential, like (\ref{sing-pot}), see, for instance, \cite{Ab, Ab2}. Indeed,
it is well-known that in the case of the regular potential (\ref{reg-pot})
the order parameter does \emph{not} remain in the physically relevant
interval $\left[ -1,1\right] $, see \cite{CMZ}. Numerical simulations show
that the system (\ref{1.1})-(\ref{1.4}) for $K=AC$ captures basic features
of two-phase flow behavior, including vesicle dynamics or drop formation
processes (cf. \cite{DLL, Sy, YFLS} and references therein). Moreover, from
the numerical point of view it is easier to implement the system based on
the Allen-Cahn equation than the system (\ref{1.1})-(\ref{1.4}) for $K=CH,$
see \cite{FHL, TLK}.

\textbf{(II)} When $K=CH,$ the phase-field component $\phi $ enjoys "good"
regularity properties and, therefore, the system (\ref{1.1})-(\ref{1.4}) is
only \emph{weakly} coupled through the Korteweg force in (\ref{1.1}) even in
the inviscid case when $\nu \equiv 0$, see \cite{CaoG}. Theoretical aspects
(i.e., well-posedness, regularity and asymptotic behavior as time goes to
infinity) for the system \eqref{1.1}-\eqref{1.4} in the case $K=CH$ have
been investigated in a sufficiently large number of papers in both two and
three dimensions. Well-posedness results for $K=CH$ when $F$ is a smooth
polynomial potential can be found in \cite{Bo1, Bo2, Bo3, CaoG, St}, and in 
\cite{Ab, Ab2, BF} when $F$ is a singular potential of the form (\ref%
{sing-pot}). Regarding the longtime behavior when $K=CH$, results about the
stability of stationary solutions were given in \cite{Ab2, Bo2, GG3}, while
theorems about the convergence to single equilibria and existence of global
and exponential attractors were proven in \cite{Ab, GG0, GG2, ZWH}. Various
numerical aspects of (\ref{1.1})-(\ref{1.4}) when $K=CH$ were investigated
in \cite{BBCV, BCB, Ja, KSW, KKL,LS, LM}. However, for the system when $K=AC$
one expects lower regularity for the $\phi $-component and, hence, in this
case (\ref{1.1})-(\ref{1.4}) is \emph{strongly} coupled. This feature has
already been present in \cite{GG1, Mej, X, XZL, ZGH}. Finally, some
existence results in the compressible case for (\ref{1.1})-(\ref{1.4}) in
both formulations (\ref{CHAC})\ are contained in \cite{AF, FPRS, DLL2, Ko2}.
A comparison of these models, providing further insight on the behavior of
the full problem (\ref{1.1})-(\ref{1.4}), is given in \cite{LBK}. The
relationship between (\ref{1.1})-(\ref{1.4}) and the standard sharp
interface models (which are obtained by taking $\varepsilon \rightarrow
0^{+} $ in (\ref{1.1})-(\ref{1.4})) is discussed in \cite{AM, CaoG}.

\textbf{(III)} For both these approaches, one expects a less or more
incomplete theory for (\ref{1.1})-(\ref{1.4}) in three dimensions because a
full mathematical theory for the 3D Navier-Stokes equation (NSE) is still
lacking at present. Moreover, as noted in \cite{HLT} direct numerical
simulation of the 3D NSE for many physical applications with high Reynolds
number flows is "\textit{intractable even using state-of-the-art numerical
methods on the most advanced supercomputers available nowadays}". Recently,
many applied mathematicians have developed regularized turbulence models for
the 3D NSE\ as an attempt to overcome this simulation barrier. Their aim is
to capture "\textit{the large, energetic eddies without having to compute
the smallest dynamically relevant eddies, by instead modelling the effects
of small eddies in terms of the large scales in the 3D NSE}". Since 1998,
many such regularized models have been proposed, tested and investigated
from both the numerical and the mathematical point of views. Among these
models, one can find the globally well-posed 3D Navier-Stokes-$\alpha $ (NS-$%
\alpha $) equations (also known as the viscous Camassa-Holm equations and
Lagrangian averaged Navier-Stokes-$\alpha $ model), the 3D Leray-$\alpha $
models, the modified 3D Leray-$\alpha $ models, the simplified 3D Bardina
models,~the 3D Navier-Stokes-Voight (NSV) equations, and their inviscid
counterparts. For instance, it has been observed that computational
simulations of the 3D Navier-Stokes-$\alpha $ (NS-$\alpha $) equations are
statistically indistinguishable from the simulations of the Navier-Stokes
equations. Furthermore, the 3D Navier-Stokes-$\alpha $ model provides
tremendous computational savings as shown in simulations of both forced and
decaying turbulence. Finally, the 3D Navier-Stokes-$\alpha $ model arises
from a variational principle in the same fashion as the Navier-Stokes
equations. We refrain from giving an exhaustive list of references but we
refer the reader to \cite{HLT} for a complete bibliography and detailed
description of the results available for these regularized models.

In this paper, upon taking the point of view described in (III), we consider
first the following prototype of initial value problem for two-phase
incompressible flows on an $n$-dimensional compact Riemannian manifold $%
\Omega $ with or without boundary, when $n=2,3$:%
\begin{equation}
\left\{ 
\begin{array}{l}
\partial _{t}u+A_{0}u+(Mu\cdot \nabla )(Nu)+\chi \nabla (Mu)^{T}\cdot
(Nu)+\nabla p=-\varepsilon \text{div}(\nabla \phi \otimes \nabla \phi )+g,
\\ 
\partial _{t}\phi +Nu\cdot \nabla \phi +\varepsilon A_{1}\phi +\varepsilon
^{-1}f\left( \phi \right) =0, \\ 
\text{div}\left( u\right) =0, \\ 
u\left( 0\right) =u_{0}, \\ 
\phi \left( 0\right) =\phi _{0},%
\end{array}%
\right.  \label{e:pde}
\end{equation}%
where $A_{0}$, $A_{1}$, $M$, and $N$ are linear operators having certain
mapping properties, and where $\chi $ is either $1$ or $0$. All kinds of
boundary conditions (i.e., periodic, no-slip, no-flux, Navier boundary
conditions, etc)\ can be treated and are included in our analysis; they will
be incorporated in the weak formulation for the problem (\ref{e:pde}), see
Section \ref{s:prelim}. We introduce three parameters which control the
degree of smoothing in the operators $A_{0}$, $M$ and $N$, namely $\theta $, 
$\theta _{1}$ and $\theta _{2}$, while $A_{1}$ is a differential operator of 
\emph{second} order. Thus we will only focus on the case when $K=AC$ in (\ref%
{CHAC}), which is actually the harder case (see (II) above). Some examples
of operators $A_{0}$, $A_{1}$, $M$, and $N$ which satisfy the mapping
assumptions we will need in this paper are 
\begin{equation}
A_{0}=\nu (-\Delta )^{\theta },\quad A_{1}=-\Delta \text{, \ \ }M=(I-\alpha
^{2}\Delta )^{-\theta _{1}},\quad N=(I-\alpha ^{2}\Delta )^{-\theta _{2}},
\label{param}
\end{equation}%
for fixed positive real numbers $\alpha ,\nu $ and for specific choices of
the real parameters $\theta $, $\theta _{1}$, and $\theta _{2}$. We note
that the Korteweg force in (\ref{e:pde}) can be equivalently rewritten in
the following form%
\begin{align}
-\varepsilon \text{div}(\nabla \phi \otimes \nabla \phi )& =\varepsilon \mu
\nabla \phi -\nabla (\frac{\varepsilon ^{2}}{2}|\nabla \phi |^{2}+F\left(
\phi \right) )  \label{1.3ter} \\
& =-\varepsilon ^{2}\Delta \phi \nabla \phi -\frac{\varepsilon ^{2}}{2}%
\nabla (|\nabla \phi |^{2}),  \notag
\end{align}%
where $\mu $ is given by (\ref{1.4}). As in \cite{HLT}, we emphasize that
the abstract mapping assumptions we employ are more general, and as a result
do not require any specific form of the parametrizations of $A_{0}$, $A_{1},$
$M$, and $N$. This abstraction allows~(\ref{e:pde}) to recover some of the
existing models that have been previously studied, as well as to represent a
much larger three-parameter family of models that have not been explicitly
studied in detail. For clarity, some of the specific regularization models
recovered by (\ref{e:pde}) for particular choices of the operators $%
A_{0},M,N $ and $\chi $ are listed in Table~\ref{t:spec}.

{\small 
\begin{table}[th]
\caption{Some special cases of the model (\protect\ref{e:pde}) with $\protect%
\alpha >0$, and with $S=(I-\protect\alpha \Delta )^{-1}$ and $S_{\protect%
\theta _{2}}~=~[I~+~(-\protect\alpha \Delta )^{\protect\theta _{2}}]^{-1}$. }
\label{t:spec}
\begin{center}
{\small $%
\begin{tabular}{|l||lllllll|}
\hline
Model & \multicolumn{1}{||l|}{NSE-AC} & \multicolumn{1}{l|}{Leray-AC-$\alpha 
$} & \multicolumn{1}{l|}{ML-AC-$\alpha $} & \multicolumn{1}{l|}{SBM-AC} & 
\multicolumn{1}{l|}{NSV-AC} & \multicolumn{1}{l|}{NS-AC-$\alpha $} & NS-AC-$%
\alpha $-like \\ \hline
$A_{0}$ & \multicolumn{1}{||l|}{$-\nu \Delta $} & \multicolumn{1}{l|}{$-\nu
\Delta $} & \multicolumn{1}{l|}{$-\nu \Delta $} & \multicolumn{1}{l|}{$-\nu
\Delta $} & \multicolumn{1}{l|}{$-\nu \Delta S$} & \multicolumn{1}{l|}{$-\nu
\Delta $} & $\nu \left( -\Delta \right) ^{\theta }$ \\ \hline
$M$ & \multicolumn{1}{||l|}{$I$} & \multicolumn{1}{l|}{$S$} & 
\multicolumn{1}{l|}{$I$} & \multicolumn{1}{l|}{$S$} & \multicolumn{1}{l|}{$S$%
} & \multicolumn{1}{l|}{$S$} & $S_{\theta _{2}}$ \\ \hline
$N$ & \multicolumn{1}{||l|}{$I$} & \multicolumn{1}{l|}{$I$} & 
\multicolumn{1}{l|}{$S$} & \multicolumn{1}{l|}{$S$} & \multicolumn{1}{l|}{$S$%
} & \multicolumn{1}{l|}{$I$} & $I$ \\ \hline
$\chi $ & \multicolumn{1}{||l|}{$0$} & \multicolumn{1}{l|}{$0$} & 
\multicolumn{1}{l|}{$0$} & \multicolumn{1}{l|}{$0$} & \multicolumn{1}{l|}{$0$%
} & \multicolumn{1}{l|}{$1$} & $1$ \\ \hline
\end{tabular}%
$  }
\end{center}
\end{table}
}Recall that $\alpha $-models of turbulence were intended as a basis for
regularizing numerical schemes for simulating turbulence in single-like
fluids \cite{HLT} (see the point (III) above). Thus, it is important to
verify whether the ad~hoc smoothed systems from Table \ref{t:spec} inherit
some of the original properties of the Navier-Stokes-Allen-Cahn (NSE-AC)
system. In particular, one would like to see if the natural energy of the
smoothed systems can be identified with the energy of the original NSE-AC
system under suitable boundary conditions. For the NSE-AC\ system there is
one essential ideal invariant (for instance,\ under rectangular periodic
boundary conditions or in the whole space), namely, the energy%
\begin{equation*}
\mathcal{E}_{0}=\frac{1}{2}\int_{\Omega }\left( |u\left( x\right)
|^{2}+\varepsilon |\nabla \phi (x)|^{2}\right) dx+\varepsilon
^{-1}\int_{\Omega }F\left( \phi \right) dx.
\end{equation*}%
In the case of the $\alpha $-models from Table \ref{t:spec}, the
corresponding ideal invariant is the energy%
\begin{equation*}
\mathcal{E}_{\alpha }=\frac{1}{2}\int_{\Omega }\left( u\left( x\right) \cdot
Nu\left( x\right) +\varepsilon |\nabla \phi (x)|^{2}\right) dx+\varepsilon
^{-1}\int_{\Omega }F\left( \phi \right) dx,
\end{equation*}%
which reduces, as $\alpha \rightarrow 0$, to the dissipated energy $\mathcal{%
E}_{0}$ of the NSE-AC equations.

Our main goal in this paper is to develop well-posedness and long-time
dynamics results for the entire three-parameter family of models, and then
subsequently recover the existing results of this type for the specific
regularization models that have been previously studied. We first aim to
establish a number of results for the entire three-parameter family,
including results on existence, regularity, uniqueness, continuous
dependence with respect to initial data, linear and nonlinear perturbations
(with the inviscid and $\alpha \rightarrow 0$ limits as special cases),
existence and finite dimensionality of global attractors, and existence of
exponential attractors (also known as inertial sets). Elaborating further on
the latter issue, we recall that in the global attractor theory, it is
usually extremely difficult (if not impossible) to estimate and to express
the rate of convergence of trajectories to the global attractor in terms of
the physical parameters of the system considered. This constitutes the main
drawback of the theory. Simple examples show that the rate of convergence
can be arbitrarily slow and non-uniform with respect to the parameters of
the system considered. As a consequence, the global attractor becomes
sensitive to small perturbations and, moreover, it may miss important
transient behaviors because the global attractor consists only of states in
the final stage. Another suitable object which always contains the global
attractor, and thus is more structurally rich in content than the global
attractor is the so-called exponential attractor. The concept of \emph{%
exponential attractor} overcomes the difficulties we mentioned earlier.
Indeed, in contrast to the global attractor theory, the relevant constants
can be explicitly found in terms of the physical parameters, and the
exponential attractor theory can provide a direct way to estimate the
fractal dimension of the global attractor even when the classical machinery
fails. Furthermore, an exponential attractor attracts bounded subsets of the
energy phase-space at an exponential rate, which makes it a more useful
object in numerical simulations than the global attractor. We refer the
reader for more details to the survey article \cite{MZ}. Following \cite{HLT}%
, our main goal is to analyze a generalized model based on abstract mapping
properties of the principal operators $A_{0}$, $M$, and $N$ allowing for a
simple analysis that helps bring out the core common structure of the
various regularized and unregularized Navier-Stokes-Allen-Cahn systems. In 
\cite{GM1}, a direct relationship between the long-term dynamics of the 3D
NSE-AC system and the three dimensional system based on the NS-$\alpha $%
-model coupled with the Allen-Cahn equation (NS-AC-$\alpha $) was
established. Note that the NS-AC-$\alpha $ system corresponds to a subset of
those problems studied here. For the other regularized models considered in
Table \ref{t:spec}, as far as we know mathematical and numerical results
have \emph{not} been previously established in the literature. The global
existence, uniqueness and regularity of solutions for these models have been
mostly known only for the two-dimensional NSE-AC system \cite{GG1, X, ZGH}
and the three-dimensional NS-AC-$\alpha $ model \cite{GM1}. Here, as a
consequence of a more general result, we develop complete well-posedness and
global regularity results for models of turbulence in two-phase flows
described by (\ref{e:pde}). Furthermore, we establish various convergence
results for the global weak solutions of (\ref{e:pde}) as either one of the
parameters $\alpha ,\nu $ goes to zero. In addition, we prove results on the
existence of finite-dimensional global and exponential attractors. Then, by
the Lojasiewicz--Simon technique, we also establish the convergence of any
bounded solution of (\ref{e:pde}) to single steady states, provided that $F$
is a real analytic function, and that the time-dependent body force $g$ is
asymptotically decaying in a precise way, i.e.,%
\begin{equation*}
\int_{t}^{\infty }\left\Vert g\left( s\right) \right\Vert _{H^{-\theta
-\theta _{2}}}^{2}ds\lesssim \left( 1+t\right) ^{-\left( 1+\delta \right) },%
\text{ for all }t\geq 0.
\end{equation*}%
In particular, for any fixed initial datum $\left( u_{0},\phi _{0}\right) $
the corresponding trajectory satisfies the estimate%
\begin{equation*}
\left\Vert u\left( t\right) \right\Vert _{H^{-\theta _{2}}}+||\phi (t)-\phi
_{\ast }||_{H^{1}}\lesssim (1+t)^{-\xi },\text{ for all }t\geq 0,
\end{equation*}%
for some $\xi \in \left( 0,1\right) ,$ depending on $\phi _{\ast }$\ and $%
\delta >0$, where $\phi _{\ast }$ is a steady-state of $A_{1}\phi _{\ast
}+\varepsilon ^{-2}f\left( \phi _{\ast }\right) =0$, $\varepsilon >0$.

We emphasize again that all these results are all new for the models
mentioned in Table \ref{t:spec}, and that the abstract mapping assumptions
we employ for (\ref{e:pde}) are more general, and as a result do not require
any specific form of the parametrizations of $A_{0}$, $M$, and $N,$ as in (%
\ref{param}). As a consequence, the framework we exploit allows~us to derive
new results for a much larger three-parameter family of models that have not
been included in Table \ref{t:spec} and explicitly studied anywhere in
detail. Finally, it is also worth emphasizing that any nondissipative (i.e., 
$\theta =0$) regularized Navier-Stokes equation in (\ref{e:pde}) can be
thought as an inviscid regularization of the usual (unregularized)
Navier-Stokes equation (NSE). Thus, in contrast to the case $\theta >0$ the
regularized system (\ref{e:pde}) for $\theta =0$ is even more strongly
coupled than before. Consequently, this feature will make the analysis even
more delicate, especially in the treatment of the long-term dynamic behavior
as time goes to infinity. Indeed, in this case the Korteweg force (\ref%
{1.3ter}) can be less regular than the convective term $(Mu\cdot \nabla
)(Nu) $ from (\ref{e:pde}), especially in three space dimensions. Besides,
our analysis can be applied \emph{verbatim }to certain regularized
simplified Ericksen-Leslie (RSEL) models for the hydrodynamics of liquid
crystals in $n$-dimensional compact Riemannian manifolds with or without
boundary. Recently, the simplified (unregularized) Ericksen-Leslie system
which consists of the $n$-dimensional NSE coupled with the Allen-Cahn
equation (in this context, also known as the Ginzburg-Landau equation) for
the orientation parameter $\phi \in \mathbb{R}^{n}$ was considered in \cite%
{Bos, Er, Le, LL, Sh, Wu} (and the references therein). Moreover, the same
system where the 3D Navier-Stokes equation is replaced by the Lagrangian
averaged 3D Navier-Stokes-$\alpha $ model was also considered in \cite[%
Section 7]{Sh}. We note that the problems studied in these references
correspond to only a subset of those regularized (RSEL) systems contained
here (see Section \ref{Leslie}).

It is also important to note that the general framework of \cite{HLT}, also
exploited and extended further here, allows for the development of new
results for certain (regularized or un-regularized) Navier-Stokes equations
and magnetohydrodynamics (MHD) models. For instance, our results on the
existence of exponential attractors, and the existence of global attractors
in non-dissipative systems (e.g., when $\theta =0$ in (\ref{e:pde}), and
when there is no coupling) are completely new and complementary to the
results of \cite[Section 5]{HLT}. In fact, in this paper we will show how to
close a gap in the proof of \cite[Section 5, Corollary 5.4]{HLT} whose
assumptions can only be verified in the case $\theta >0$. Indeed, one can
easily observe that when $\theta =0$, the assumptions of \cite[Theorem 5.1,
(b)]{HLT} do not longer provide the existence of a compact absorbing set as
claimed on \cite[pg. 550]{HLT}. For instance, the 3D Navier-Stokes-Voight
model, or any other non-dissipative system when $\theta =0$, is no longer
covered by the result of \cite[Corollary 5.4]{HLT} (see also Section \ref%
{rem-mhd}).

The remainder of the paper is structured as follows. In~Section \ref%
{s:prelim}, we establish our notation and give some basic preliminary
results for the operators appearing in the general regularized model.
In~Section \ref{s:well}, we build some well-posedness results for the
general model; in particular, we establish existence results (Section \ref%
{ss:exist}), regularity results (Section \ref{ss:reg}), and uniqueness and
continuous dependence results (Section \ref{ss:stab}). In~Section \ref%
{s:pert} we establish some results for singular perturbations, which as
special cases include the inviscid limit of viscous models and the $\alpha
\rightarrow 0$ limit in $\alpha $ models; this involves a separate analysis
of the linear (Section \ref{ss:pert-lin}) and nonlinear (Section \ref%
{ss:pert-nonlin}) terms. In~Section \ref{s:longbehav}, we show existence of
a global attractor for the general model by dissipation arguments (Sections %
\ref{ss:attr-diss} and \ref{ss:attr-nondiss}), and then by employing the
approach from~\cite{GMPZ, MZ}, to show the existence of exponential
attractors (Sections \ref{ss:attr-diss} and \ref{ss:attr-nondiss}).
In~Section \ref{s:convss}, we establish asymptotic stability results as time
goes to infinity of solutions to our regularized models, with the help from
a Lojasiewicz--Simon technique. Section \ref{rem-mhd} contains several
important new theorems and remarks for the systems considered by \cite{HLT}.
Section \ref{Leslie} contains some additional remarks on a regularized
system for the simplified Ericksen-Leslie model for the hydrodynamics of
liquid crystals. To make the paper sufficiently self-contained, our final
Section \ref{ss:app} contains supporting material on Sobolev and Gr\"{o}%
nwall-type inequalities, and several other abstract results which are needed
to prove our main results.

\section{Preliminary material}

\label{s:prelim}

We follow the same framework and notation as in \cite{HLT}. To this end, let 
$\Omega $ be an $n$-dimensional smooth compact manifold with or without
boundary and equipped with a volume form, and let $E\rightarrow \Omega $ be
a vector bundle over $\Omega $ equipped with a Riemannian metric $h=\left(
h_{ij}\right) _{n\times n}$. With $C^{\infty }(E)$ denoting the space of
smooth sections of $E$, let $\mathcal{V}\subseteq C^{\infty }(E)$ be a
linear subspace, let $A_{0}:\mathcal{V}\rightarrow \mathcal{V}$ be a linear
operator, and let $B_{0}:\mathcal{V}\times \mathcal{V}\rightarrow \mathcal{V}
$ be a bilinear map. At this point $\mathcal{V}$ is conceived to be an
arbitrary linear subspace of $C^{\infty }(E)$; however, later on we will
impose restrictions on $\mathcal{V}$ implicitly through various conditions
on certain operators such as $A_{0}$. Furthermore, we let $\mathcal{W}%
\subseteq C^{\infty }(\Omega )$ be a linear subspace and let $A_{1}:\mathcal{%
W}\rightarrow \mathcal{W}$ be a linear operator satisfying various
assumptions below. In order to define the variational setting for the
phase-field component we also need to introduce the bilinear operators $%
R_{0}:\mathcal{W}\times \mathcal{W\rightarrow V}$, $B_{1}:\mathcal{V}\times 
\mathcal{W}\rightarrow \mathcal{W}$, as follows:%
\begin{equation}
B_{1}\left( u\left( x\right) ,\phi \left( x\right) \right) :=Nu\left(
x\right) \cdot \nabla \phi \left( x\right) ,\text{ }R_{0}\left( \psi \left(
x\right) ,\phi \left( x\right) \right) :=\psi \left( x\right) \nabla \phi
\left( x\right) .  \label{b00}
\end{equation}%
Recalling (\ref{1.3ter}) and assuming that $\varepsilon =1$ (for the sake of
simplicity), the initial data $u_{0}\in \mathcal{V}$, $\phi _{0}\in \mathcal{%
W}$ and forcing term $g\in C^{\infty }(0,T;\mathcal{V})$ with $T>0$,
consider the following system%
\begin{equation}
\left\{ 
\begin{array}{l}
\partial _{t}u+A_{0}u+B_{0}(u,u)=R_{0}\left( A_{1}\phi ,\phi \right) +g, \\ 
\partial _{t}\phi +B_{1}\left( u,\phi \right) +\mu =0, \\ 
\mu =A_{1}\phi +f\left( \phi \right) , \\ 
u\left( 0\right) =u_{0},\phi \left( 0\right) =\phi _{0},%
\end{array}%
\right.  \label{e:op}
\end{equation}%
on the time interval $[0,T]$. Bearing in mind the model (\ref{e:pde}), we
are mainly interested in bilinear maps of the form 
\begin{equation}
B_{0}(v,w)=\bar{B}_{0}(Mv,Nw),  \label{e:b-def}
\end{equation}%
where $M$ and $N$ are linear operators in $\mathcal{V}$ that are in some
sense regularizing and are relatively flexible, and $\bar{B}_{0}$ is a
bilinear map fixing the underlying nonlinear structure of the fluid
equation. In the following, let $P:C^{\infty }(E)\rightarrow \mathcal{V}$ be
the $L^{2}$-orthogonal projector onto $\mathcal{V}$. Finally, concerning the
derivative $f$\ of the function $F$ in (\ref{1.5}) we will focus mostly on
the regular potential case when $f\in C^{2}\left( \mathbb{R},\mathbb{R}%
\right) $ satisfies $f\left( 1\right) \geq 0,$ $f\left( -1\right) \leq 0$
and obeys the following condition%
\begin{equation}
\underset{\left\vert r\right\vert \rightarrow \infty }{\lim \inf }%
f^{^{\prime }}\left( r\right) >0.  \label{condf}
\end{equation}%
However, when $f$ is a singular potential, see Remark \ref{sing-rem}.

We will study the regularized system (\ref{e:op}) by extending it to
function spaces that have weaker differentiability properties. To this end,
we interpret (\ref{e:op}) in distribution sense, and need to continuously
extend $A_{0},$ $A_{1}$ and $B_{0},B_{1}$ and $R_{0}$ to appropriate
smoothness spaces. Namely, we employ the spaces $V^{s}=\mathrm{clos}_{H^{s}}%
\mathcal{V}$, $W^{s}=\mathrm{clos}_{H^{s}}\mathcal{W}$, which will
informally be called Sobolev spaces in the following. The pair of spaces $%
V^{s}$ and $V^{-s}$ are equipped with the duality pairing $\left\langle
\cdot ,\cdot \right\rangle $, that is, the continuous extension of the $%
L^{2} $-inner product on $V^{0}$. Same applies to the triplet $W^{s}\subset
W^{0}=\left( W^{0}\right) ^{\ast }\subset W^{-s}.$ Moreover, we assume that
there are self-adjoint \emph{positive} operators $\Lambda $ and $A_{1},$
respectively, such that $\Lambda ^{s}:V^{s}\rightarrow V^{0},$ $%
A_{1}^{s/2}:W^{s}\rightarrow W_{0}$ are isometries for any $s\in \mathbb{R}$%
, and $\Lambda ^{-1}$, $\left( A_{1}\right) ^{-1}$ are compact operators.
For arbitrary real $s$, assume that $A_{0}$, $A_{1}$, $M$, and $N$ can be
continuously extended so that 
\begin{equation}
A_{0}:V^{s}\rightarrow V^{s-2\theta },\quad A_{1}:W^{s}\rightarrow
W^{s-2},\quad M:V^{s}\rightarrow V^{s+2\theta _{1}},\quad \text{and}\quad
N:V^{s}\rightarrow V^{s+2\theta _{2}},  \label{e:bdd-amn}
\end{equation}%
are bounded operators. Again, we emphasize that the assumptions we will need
for $A_{0}$, $M$, and $N$ are more general, and do not require this
particular form of the parametrization (see (\ref{e:coercive-a})-(\ref%
{e:coercive-abis}) below). We will assume $\theta \geq 0$ and no \emph{a
priori} sign restrictions on $\theta _{1}$, $\theta _{2}.$ We remark that $s$
in (\ref{e:bdd-amn}) is assumed to be arbitrary for the purpose of the
discussion in this section; of course, it suffices to assume (\ref{e:bdd-amn}%
) for a limited range of $s$ for most of the results in this paper. The
canonical norm in the Hilbert spaces $V^{s}$ and $W^{s},$ respectively, will
be denoted by the same quantity $\left\Vert \cdot \right\Vert _{s}$ whenever
no further confusion arises, while we will use the notation $\left\Vert
\cdot \right\Vert _{L^{p}}$ for the $L^{p}$-norm. Furthermore, we assume
that $A_{0}$ and $N$ are both self-adjoint, and coercive in the sense that
for $\beta \in \mathbb{R}$, 
\begin{equation}
\left\langle A_{0}v,\Lambda ^{2\beta }v\right\rangle \geq c_{A_{0}}\Vert
v\Vert _{\theta +\beta }^{2}-C_{A_{0}}\Vert v\Vert _{\beta }^{2},\qquad v\in
V^{\theta +\beta },  \label{e:coercive-a}
\end{equation}%
with $c_{A_{0}}=c_{A_{0}}(\beta )>0$, and $C_{A_{0}}=C_{A_{0}}(\beta )\geq 0$%
, and that 
\begin{equation}
\left\langle Nv,v\right\rangle \geq c_{N}\Vert v\Vert _{-\theta
_{2}}^{2},\qquad v\in V^{-\theta _{2}},  \label{e:coercive-n}
\end{equation}%
with $c_{N}>0$. We also assume that%
\begin{equation}
\left\langle A_{0}v,Nv\right\rangle \geq c_{A_{0}}\Vert v\Vert _{\theta
-\theta _{2}}^{2},\qquad v\in V^{\theta -\theta _{2}},
\label{e:coercive-abis}
\end{equation}%
Note that if $\theta =0$, (\ref{e:coercive-a}) is strictly speaking not
coercivity and follows from the boundedness of $A$, and note also that (\ref%
{e:coercive-n}) implies the invertibility of $N$.

As examples, one may typically consider the following operators in various
combinations in (\ref{e:op}).

\begin{example}
\label{x:spaces}

(a) When $\Omega $ is a closed Riemannian manifold, and $E=T\Omega $ the
tangent bundle, an example of $\mathcal{V}$ is $\mathcal{V}_{\mathrm{per}%
}\subseteq \{u\in C^{\infty }(T\Omega ):\mathrm{div}\,u=0\}$, a subspace of
the divergence-free functions. The space of periodic functions with
vanishing mean on a torus $\mathbb{T}^{n}$ is a special case of this
example. In this case, one typically has $A_{0}=(-\Delta )^{\theta }$, $%
M=(I-\alpha ^{2}\Delta )^{-\theta _{1}}$, $N=(I-\alpha ^{2}\Delta )^{-\theta
_{2}}$ and $A_{1}=-\Delta ,$ as operators that satisfy (\ref{e:bdd-amn}),
cf. also \cite[Example 2.1, (a)]{HLT}.

(b) When $\Omega $ is a compact Riemannian manifold with boundary, and again 
$E=T\Omega $ the tangent bundle, a typical example of $\mathcal{V}$ is $%
\mathcal{V}_{\mathrm{hom}}=\{u\in C_{0}^{\infty }(T\Omega ):\mathrm{div}%
\,u=0\}$ the space of compactly supported divergence-free functions. In this
case, one may consider the choices $A_{0}=(-P\Delta )^{\theta }$, $%
A_{1}=-\Delta $, $M=(I-\alpha ^{2}P\Delta )^{-\theta _{1}}$, and $%
N=(I-\alpha ^{2}P\Delta )^{-\theta _{2}}$, respectively, as operators
satisfying (\ref{e:bdd-amn}), cf. also \cite[Example 2.1, (b)]{HLT}.
\end{example}

\begin{example}
\label{y:spaces}

(a) In Example~\ref{x:spaces} above, the bilinear map $\bar{B}_{0}$ can be
taken to be 
\begin{equation}
\bar{B}_{0\chi }(v,w)=P[(v\cdot \nabla )w+\chi (\nabla w^{T})v],
\label{b00bis}
\end{equation}%
which correspond to the models with $\chi \in \left\{ 0,1\right\} $ as
discussed in~the Introduction (see Table \ref{t:spec}).

(b) Let $\Omega $ be connected Riemannian $n$-dimensional manifold with
non-empty (sufficiently smooth) boundary $\partial \Omega .$ Define $%
A_{1}=-\Delta $, as the Laplacian of the metric $h$, acting on 
\begin{equation*}
D\left( A_{1}\right) =\left\{ \phi \in W^{2}:\partial _{\zeta }\phi =0\text{
on }\partial \Omega \right\} ,
\end{equation*}%
where $\zeta $ is an outward unit normal vector field of $\partial \Omega .$
Recall that in local coordinates $\left\{ x_{i}\right\} _{i=1}^{n},$ the
Laplacian reads%
\begin{equation*}
\Delta \left( \cdot \right) =\frac{1}{\sqrt{\det \left( h\right) }}%
\sum\nolimits_{i,j=1}^{n}\partial _{x_{j}}\left( h^{ij}\sqrt{\det \left(
h\right) }\partial _{x_{i}}\left( \cdot \right) \right) ,
\end{equation*}%
where the matrix $\left( h^{ij}\right) $ is the inverse matrix of $h.$ We
have that $A_{1}$ is a nonnegative self-adjoint operator on $W^{0}$. Next,
consider $\overline{A}_{1}=A_{1}+\gamma I,$ for some $\gamma >0$ and define $%
f_{1}\left( r\right) =f\left( r\right) -\gamma r.$ In this case, we can
rewrite the second and third equations of (\ref{e:op}) in the form%
\begin{equation*}
\partial _{t}\phi +B_{1}\left( u,\phi \right) +\overline{A}_{1}\phi
+f_{1}\left( \phi \right) =0,
\end{equation*}%
with the function $f_{1}$ still obeying assumption (\ref{condf}). Clearly, $%
\overline{A}_{1}$ is positive and it can now be continuously extended so
that it satisfies the corresponding condition from (\ref{e:bdd-amn}). Hence,
our original restriction that $A_{1}$ is \emph{positive} is indeed not
necessary and, thus, the above framework also allows us to deal with a
nonnegative selfadjoint operator $A_{1}.$
\end{example}

\begin{example}
\label{z:spaces} (\emph{Navier boundary conditions}). Let $\Omega \subset 
\mathbb{R}^{n}$ be a bounded connected domain with sufficiently smooth
boundary $\partial \Omega .$ We take $\mathcal{V}$ as $\mathcal{V}_{\mathrm{%
Nbc}}=\{u\in C^{\infty }(T\Omega ):\mathrm{div}\,u=0,$ $u\cdot \zeta =0$ on $%
\partial \Omega \}$ and recall the classical decomposition%
\begin{equation*}
V^{0}=V^{1}\oplus \left( V^{1}\right) ^{\perp },\text{ }\left( V^{1}\right)
^{\perp }:=\left\{ \nabla u:u\in V^{1}\right\} .
\end{equation*}%
When two vectors $u$ and $v$ are divergent free, $u$ satisfies the Navier
boundary condition%
\begin{equation}
2\left[ \left( Du\right) n\right] \cdot \tau +\sigma u\cdot \tau =0\text{,
on }\partial \Omega ,  \label{bc3}
\end{equation}%
($2Du:=\nabla u+(\nabla u)^{tr}$ is the "usual" deformation tensor and $%
\sigma >0$ is some friction coefficient), and $u,v\in \mathcal{V}_{\mathrm{%
Nbc}}$, the Green formula \cite{SS} yields%
\begin{equation}
\int_{\Omega }\left( -\Delta u\right) \cdot vdx=2\int_{\Omega
}Du:Dvdx+\sigma \int_{\Gamma }\left( u\cdot \tau \right) \left( v\cdot \tau
\right) dS.  \label{Green}
\end{equation}%
On the basis of (\ref{Green}), one can define the bilinear form%
\begin{equation}
\rho _{\sigma }\left( u,v\right) :=2\int_{\Omega }Du:Dvdx+\sigma
\int_{\Gamma }\left( u\cdot \tau \right) \left( v\cdot \tau \right) dS\text{%
, for }u,v\in V^{1}.  \label{form}
\end{equation}%
Note that $\rho _{\sigma }\left( \cdot ,\cdot \right) $ is bounded on $V^{1}$
and $\rho _{\sigma }\left( u,u\right) >0$ for all $u\in V^{1}.$ Then, by
Korn's inequality \cite{Ko} we see that there exists a constant $C=C\left(
\sigma \right) >0$ (independent of $u$) such that%
\begin{equation}
C\left\Vert u\right\Vert _{1}^{2}\leq \rho _{\sigma }\left( u,u\right) \text{%
, for all }u\in V^{1}.  \label{lbform}
\end{equation}%
In view of (\ref{lbform}), we see that the bilinear form $\rho _{\sigma }$
is symmetric and coercive on $V^{1}$; for vector fields that satisfy (\ref%
{bc3}) we set%
\begin{equation*}
V_{\sigma }^{2}:=\left\{ u\in V^{2}:\text{ (\ref{bc3}) holds on }\partial
\Omega \right\} .
\end{equation*}%
The corresponding Stokes operator $A_{0}$ associated with the form $\rho
_{\sigma }$ can be constructed on the basis of first and second
representations, as follows:%
\begin{equation*}
A_{0}u=\Theta _{\sigma }u,\text{ }u\in D\left( A_{\sigma }\right) :=\left\{
u\in V^{1}:\Theta _{\sigma }u\in V^{0}\right\} ,
\end{equation*}%
where $\Theta _{\sigma }\in \mathcal{L(}V^{1},\left( V^{1}\right) ^{\ast })$
is a bounded one-to-one mapping such that%
\begin{equation*}
\rho _{\sigma }\left( u,v\right) =\left\langle u,\Theta _{\sigma
}\right\rangle _{V^{1},\left( V^{1}\right) ^{\ast }}
\end{equation*}%
for all $u,v\in V^{1}$. As a byproduct, we also obtain the following compact
embeddings $D\left( A_{0}\right) \subset V^{1}\subset V^{0}\subset \left(
V_{\sigma }^{1}\right) ^{\ast },$ the operator $A_{0}$ is a positive,
selfadjoint operator on $V^{0}$, and $\left( I+A_{0}\right) ^{-1}$ is
compact. Exploiting the formula (\ref{Green}) once more, we see that $Pv=v$,
for all $v\in V^{1}$ and $\rho _{\sigma }\left( u,v\right) =\left\langle
P\left( -\Delta u\right) ,v\right\rangle ,$ \ for all $u\in V_{\sigma }^{2}$
and $v\in V^{1}$. This implies, on the basis of standard regularity theory 
\cite{SS}\ for the Stokes operator, that $D\left( A_{0}\right) =V_{\sigma
}^{2}$ and $A_{0}u=P\left( -\Delta u\right) ,$ for all $u\in D\left(
A_{0}\right) .$ Finally, fractional powers $A_{0}^{\theta }$ are also
well-defined for all $\theta \geq 0$.
\end{example}

To refer to the above examples, let us further introduce the shorthand
notation:%
\begin{equation}
B_{0\chi }\left( v,w\right) =\overline{B}_{0\chi }\left( Mv,Nw\right) ,\text{
}\chi \in \left\{ 0,1\right\} .  \label{b01}
\end{equation}%
For clarity, we list in Table \ref{t:spec2} the corresponding values of the
parameters and bilinear maps discussed above for special cases listed in
Table \ref{t:spec}.{\small 
\begin{table}[th]
\caption{Values of the parameters $\protect\theta ,$ $\protect\theta _{1}$
and $\protect\theta _{2}$, and the particular form of the bilinear map $%
B_{0} $ for some special cases of the model (\protect\ref{e:op}). (The
bilinear maps $B_{00}$ and $B_{01}$ are as in (\protect\ref{b01})).}
\label{t:spec2}
\begin{center}
{\small $%
\begin{tabular}{|c|c|c|c|c|c|c|c|}
\hline\hline
Model & NSE-AC & Leray-AC-$\alpha $ & ML--AC-$\alpha $ & SBM-AC & NSV-AC & 
NS--AC-$\alpha $ & NS--AC-$\alpha $-like \\ \hline
$\theta $ & 1 & 1 & 1 & 1 & 0 & 1 & $\theta $ \\ 
$\theta _{1}$ & 0 & 1 & 0 & 1 & 1 & 0 & 0 \\ 
$\theta _{2}$ & 0 & 0 & 1 & 1 & 1 & 1 & $\theta _{2}$ \\ 
$B_{0}$ & $B_{00}$ & $B_{00}$ & $B_{00}$ & $B_{00}$ & $B_{00}$ & $B_{01}$ & $%
B_{01}$ \\ \hline\hline
\end{tabular}%
$  }
\end{center}
\end{table}
}

Next, we denote the trilinear forms 
\begin{equation}
b_{0}(u,v,w)=\langle B_{0}(u,v),w\rangle ,\text{ }b_{1}\left( u,\phi ,\psi
\right) =\left\langle B_{1}\left( u,\phi \right) ,\psi \right\rangle ,
\label{b01bis}
\end{equation}%
and similarly the forms $\bar{b}_{0\chi }$ and $b_{0\chi }$, following (\ref%
{b00}), (\ref{b00bis}) and (\ref{b01}). Then our notion of \emph{weak
solution} for problem (\ref{e:op}) can be formulated as follows.

\begin{definition}
\label{weak}\textit{Let }$g\left( t\right) \in L^{2}(0,T;V^{s})$ for some $%
s\in \mathbb{R},$ and%
\begin{equation*}
\left( u_{0},\phi _{0}\right) \in \mathcal{Y}_{\theta _{2}}:=V^{-\theta
_{2}}\times \left( W^{1}\cap \left\{ \phi _{0}\in L^{\infty }\left( \Omega
\right) :\left\vert \phi _{0}\right\vert \leq 1\right\} \right) .
\end{equation*}%
\textit{Find a pair of functions} 
\begin{equation}
\left( u,\phi \right) \in L^{\infty }\left( 0,T;\mathcal{Y}_{\theta
_{2}}\right) \cap L^{2}\left( 0,T;V^{\theta -\theta _{2}}\times W^{2}\right)
\label{1.8}
\end{equation}%
\textit{such that}%
\begin{equation}
\partial _{t}u\in L^{p}\left( 0,T;V^{-\gamma }\right) ,\text{ }\partial
_{t}\phi \in L^{2}\left( 0,T;W^{-2}\right)  \label{1.8bis}
\end{equation}%
for some $p>1$ and $\gamma \geq 0$, \textit{such that }$\left( u,\phi
\right) $\textit{\ fulfills }$u\left( 0\right) =u_{0},$ $\phi \left(
0\right) =\phi _{0}$ \textit{and satisfies}%
\begin{align}
& \int\nolimits_{0}^{T}\left( -\left\langle u\left( t\right) ,w^{^{\prime
}}\left( t\right) \right\rangle +\left\langle A_{0}u\left( t\right) ,w\left(
t\right) \right\rangle +b_{0}\left( u\left( t\right) ,u\left( t\right)
,w\left( t\right) \right) \right) dt  \label{weak1} \\
& =\int\nolimits_{0}^{T}\left( \left\langle g\left( t\right) ,w\left(
t\right) \right\rangle +\left\langle R_{0}\left( A_{1}\phi \left( t\right)
,\phi \left( t\right) \right) ,w\left( t\right) \right\rangle \right) dt, 
\notag
\end{align}%
\begin{equation}
\int\nolimits_{0}^{T}\left( -\left\langle \phi \left( t\right) ,\psi
^{^{\prime }}\left( t\right) \right\rangle +\left\langle \mu \left( t\right)
,\psi \left( t\right) \right\rangle +b_{1}\left( u\left( t\right) ,\phi
\left( t\right) ,\psi \left( t\right) \right) \right) dt=0,  \label{weak2}
\end{equation}%
for any $\left( w,\psi \right) \in C_{0}^{\infty }\left( 0,T;\mathcal{V}%
\times \mathcal{W}\right) $, such that $\mu \left( t\right) =A_{1}\phi
\left( t\right) +f\left( \phi \left( t\right) \right) $ a.e. on $\Omega
\times \left( 0,T\right) .$
\end{definition}

\begin{remark}
As far as the interpretation of the initial conditions $u\left( 0\right)
=u_{0},$ $\phi \left( 0\right) =\phi _{0}$ is concerned, note that
properties (\ref{1.8})-(\ref{1.8bis}) imply that $u\in C(0,T;V^{-\gamma })$
and $\phi \in C(0,T;W^{0})$. Thus, the initial conditions are satisfied in a
weak sense.
\end{remark}

\section{Well-posedness results}

\label{s:well}

Analogous to the theory for the Navier-Stokes-Alle-Cahn system \cite{GG0,
GG1}, we begin to develop a solution theory for the general three-parameter
family of regularized models. We begin by showing energy estimates that will
be used to establish existence and regularity results, and under appropriate
assumptions also uniqueness and stability. At the end of the proof of each
theorem, we give the corresponding conditions for $\left( \theta ,\theta
_{1},\theta _{2}\right) $ which allow us to not only recover old results,
but also establish new results in the literature especially for the cases
listed in Table \ref{t:spec}. Throughout the paper, $C\geq 0$ will denote a 
\emph{generic} constant whose further dependence on certain quantities will
be specified on occurrence. The value of the constant can change even within
the same line. Furthermore, we introduce the notation $a\lesssim b$ to mean
that there exists a constant $C>0$ such that $a\leq Cb.$ This notation will
be used when the constant $C$ is irrelevant and becomes tedious.

\subsection{Existence of weak solutions}

\label{ss:exist}

In this subsection, we establish sufficient conditions for the existence of
weak solutions to the problem (\ref{e:op}) (cf. Definition \ref{weak}). As
noted in the Introduction, in the case $K=AC$ a maximum principle holds for
the phase--field component of any weak solution.

\begin{proposition}
\label{maxp}Suppose that $f$ satisfies (\ref{condf}), $f\left( 1\right) \geq
0,$ $f\left( -1\right) \leq 0$, and $b_{1}\left( v,\psi ,\psi \right) =0$,
for any $v\in V^{\theta -\theta _{2}},$ $\psi \in W^{1}.$ Let $\phi _{0}\in
L^{\infty }\left( \Omega \right) $ such that $\left\vert \phi
_{0}\right\vert \leq 1$ a.e. in $\Omega .$ Then for any weak solution $%
\left( u,\phi \right) $ to problem (\ref{e:op}) in the sense of Definition %
\ref{weak}, we have $\phi \in L^{\infty }\left( 0,T;L^{\infty }\left( \Omega
\right) \right) $ and%
\begin{equation}
\left\vert \phi \left( t\right) \right\vert _{L^{\infty }\left( \Omega
\right) }\leq 1\text{, a.e. on }\left( 0,T\right) .  \label{6.3}
\end{equation}
\end{proposition}

\begin{proof}
For the reader's convenience, a proof of the above statement is contained in 
\cite[Theorem 6.1]{GG1}.
\end{proof}

\begin{theorem}
\label{t:exist} Let the assumptions of Proposition \ref{maxp} and the
following conditions hold.

\begin{itemize}
\item[i)] $\left( u_{0},\phi _{0}\right) \in \mathcal{Y}_{\theta _{2}}$ with
any $\theta _{2}\geq -1$, and $g\in L^{2}(0,T;V^{-\theta -\theta _{2}})$, $%
T>0$.

\item[ii)] $b_{0}(v,v,Nv)=0,$ for any $v\in V^{\theta -\theta _{2}}$;

\item[iii)] $b_{0}:V^{\bar{\sigma}_{1}}\times V^{\bar{\sigma}_{2}}\times V^{%
\bar{\gamma}}\rightarrow \mathbb{R}$ is bounded for some $\bar{\sigma}%
_{i}<\theta -\theta _{2}$, $i=1,2$, and $\bar{\gamma}\geq \gamma $;

\item[iv)] $b_{0}:V^{\sigma _{1}}\times V^{\sigma _{2}}\times V^{\gamma
}\rightarrow \mathbb{R}$ is bounded for some $\sigma _{i}\in \lbrack -\theta
_{2},\theta -\theta _{2}]$, $i=1,2$, and $\gamma \in \lbrack \theta +\theta
_{2},\infty )\cap (\theta _{2},\infty )\cap (\frac{n}{2},\infty )$;

Then, there exists at least one weak solution $\left( u,\phi \right) $
satisfying (\ref{1.8})-(\ref{1.8bis}) such that%
\begin{equation*}
p=\left\{ 
\begin{array}{ll}
\min \{2,\frac{2\theta }{\sigma _{1}+\sigma _{2}+2\theta _{2}}\}, & \text{if 
}\theta >0, \\ 
2, & \text{if }\theta =0.%
\end{array}%
\right.
\end{equation*}
\end{itemize}
\end{theorem}

\begin{proof}
Let $\{V_{m}:m\in \mathbb{N}\}\subset V^{\theta -\theta _{2}},$ $%
\{W_{m}:m\in \mathbb{N}\}\subset D\left( A_{1}\right) \cap L^{\infty }\left(
\Omega \right) $ be sequences of finite dimensional subspaces of $V^{\theta
-\theta _{2}}$ and $D\left( A_{1}\right) ,$ respectively, such that

\begin{enumerate}
\item $V_{m}\subset V_{m+1},$ $W_{m}\subset W_{m+1},$ for all $m\in \mathbb{N%
}$;

\item $\cup _{m\in \mathbb{N}}V_{m}$ is dense in $V^{\theta -\theta _{2}}$,
and $\cup _{m\in \mathbb{N}}W_{m}$ is dense in $D\left( A_{1}\right) ;$

\item For $m\in \mathbb{N}$, with $\widetilde{V}_{m}=NV_{m}\subset V^{\theta
+\theta _{2}}$, the projectors $P_{m}:V^{\theta -\theta _{2}}\rightarrow
V_{m},$ $Q_{m}:D\left( A_{1}\right) \rightarrow W^{m}$, defined by 
\begin{align*}
\left\langle P_{m}v,w_{m}\right\rangle & =\left\langle v,w_{m}\right\rangle
,\qquad w_{m}\in \widetilde{V}_{m},\,v\in V^{\theta -\theta _{2}}, \\
\left\langle Q_{m}\phi ,\psi _{m}\right\rangle & =\left\langle \phi ,\psi
_{m}\right\rangle ,\qquad \psi _{m}\in W_{m},\,\phi \in D\left( A_{1}\right)
,
\end{align*}%
are uniformly bounded as maps from $V^{-\gamma }\rightarrow V^{-\gamma }$
and $W^{-2}\rightarrow W^{-2}$, respectively.
\end{enumerate}

Such sequences can be constructed e.g., by using the eigenfunctions of the
isometries $\Lambda ^{1+\theta }:V^{1+\theta -\theta _{2}}\rightarrow
V^{-\theta _{2}}$, $A_{1}:D\left( A_{1}\right) \rightarrow W^{0}.$ Consider
the problem of finding $\left( u_{m},\phi _{m}\right) \in
C^{1}(0,T;V_{m}\times W_{m})$ such that for all $\left( w_{m},\psi
_{m}\right) \in \widetilde{V}_{m}\times W_{m}$,%
\begin{equation}
\left\{ 
\begin{array}{l}
\langle \partial _{t}u_{m},w_{m}\rangle +\langle A_{0}u_{m},w_{m}\rangle
+b_{0}(u_{m},u_{m},w_{m})=\langle g,w_{m}\rangle +\left\langle R_{0}\left(
\mu _{m},\phi _{m}\right) ,w_{m}\right\rangle , \\ 
\left\langle \partial _{t}\phi _{m},\psi _{m}\right\rangle +\left\langle \mu
_{m},\psi _{m}\right\rangle +b_{1}\left( u_{m},\phi _{m},\psi _{m}\right) =0,
\\ 
\mu _{m}=A_{1}\phi _{m}+Q_{m}f\left( \phi _{m}\right) , \\ 
\langle u_{m}(0),w_{m}\rangle =\langle u_{0},w_{m}\rangle , \\ 
\left\langle \phi _{m}\left( 0\right) ,\psi _{m}\right\rangle =\left\langle
\phi _{0},\psi _{m}\right\rangle .%
\end{array}%
\right.  \label{e:galerkin}
\end{equation}%
Upon choosing a basis for $V_{m}\times W_{m}$, the above becomes an initial
value problem for a system of ODE's, and moreover since $N$ is invertible by
(\ref{e:coercive-n}), the standard ODE theory gives a unique local-in-time
solution. Furthermore, this solution is global if its norm is finite at any
finite time instance. The fourth equality in (\ref{e:galerkin}) gives 
\begin{equation*}
c_{N}\Vert u_{m}(0)\Vert _{-\theta _{2}}^{2}\leq \langle
u_{m}(0),Nu_{m}(0)\rangle =\langle u(0),Nu_{m}(0)\rangle \leq \Vert
u(0)\Vert _{-\theta _{2}}\Vert Nu_{m}(0)\Vert _{\theta _{2}},
\end{equation*}%
so that 
\begin{equation}
\Vert u_{m}(0)\Vert _{-\theta _{2}}\leq \frac{\Vert N\Vert _{-\theta
_{2};\theta _{2}}}{c_{N}}\Vert u(0)\Vert _{-\theta _{2}}.
\end{equation}%
Now in the first and second equalities of (\ref{e:galerkin}), taking $%
w_{m}=Nu_{m}$ and $\psi _{m}=\mu _{m}$, respectively, and using the
condition ii) on $b_{0}$, we get after standard transformations%
\begin{equation}
\begin{split}
& \frac{d}{dt}\left( \langle u_{m},Nu_{m}\rangle +\left\Vert A_{1}^{1/2}\phi
_{m}\right\Vert _{L^{2}}^{2}+2\int_{\Omega }F\left( \phi _{m}\right)
dx\right) +2\langle A_{0}u_{m},Nu_{m}\rangle +\left\Vert \mu _{m}\right\Vert
_{L^{2}}^{2} \\
& =2\langle g_{m},Nu_{m}\rangle \\
& \leq \varepsilon ^{-1}\Vert g\Vert _{-\theta -\theta _{2}}^{2}+\varepsilon
\Vert N\Vert _{-\theta _{2};\theta _{2}}^{2}\Vert u_{m}\Vert _{\theta
-\theta _{2}}^{2},
\end{split}
\label{e:exist-1}
\end{equation}%
for any $\varepsilon >0$. Let us now set%
\begin{equation*}
\mathcal{E}\left( u,\phi \right) :=\left\langle u,Nu\right\rangle
+\left\Vert A_{1}^{1/2}\phi \right\Vert _{L^{2}}^{2}+2\int_{\Omega }F\left(
\phi \right) dx+C_{F},
\end{equation*}%
for some sufficiently large positive constant $C_{F}$ such that $\mathcal{E}%
\geq 0$ (indeed, such a constant exists due to assumption (\ref{condf})).
Choosing $\varepsilon >0$ small enough, we can ensure 
\begin{equation*}
-2\langle A_{0}u_{m},Nu_{m}\rangle +\varepsilon \Vert N\Vert _{-\theta
_{2};\theta _{2}}^{2}\Vert u_{m}\Vert _{\theta -\theta _{2}}^{2}\lesssim
-\Vert u_{m}\Vert _{\theta -\theta _{2}}^{2},
\end{equation*}%
so that by Gr\"{o}nwall's inequality we have%
\begin{equation}
\mathcal{E}\left( u_{m}\left( t\right) ,\phi _{m}\left( t\right) \right)
\lesssim \mathcal{E}\left( u_{m}\left( 0\right) ,\phi _{m}\left( 0\right)
\right) +C\int_{0}^{t}\Vert g\Vert _{-\theta -\theta _{2}}^{2},
\label{e:exist-2}
\end{equation}%
for some $C>0$. For any fixed $T>0$, this gives%
\begin{equation*}
\left( u_{m},\phi _{m}\right) \in L^{\infty }(0,T;V^{-\theta _{2}}\times
W^{1})
\end{equation*}%
with uniformly (in $m$) bounded norm. By Proposition \ref{maxp}, we also
have the uniform bound $\phi _{m}\in L^{\infty }\left( 0,T;L^{\infty }\left(
\Omega \right) \right) $ such that $\left\vert \phi _{m}\left( t\right)
\right\vert \leq 1$ a.e. in $\Omega \times \left( 0,T\right) .$ Moreover,
integrating (\ref{e:exist-1}), and taking into account (\ref{e:exist-2}), we
infer 
\begin{equation}
\int_{0}^{t}\langle A_{0}u_{m}\left( s\right) ,Nu_{m}\left( s\right) \rangle
ds\leq C_{T},\text{ }\int_{0}^{t}\left\Vert \mu _{m}\left( s\right)
\right\Vert _{L^{2}}^{2}ds\leq C_{T},\qquad t\in \left( 0,T\right) ,
\label{e:exist-3}
\end{equation}%
If $\theta >0$, by the coerciveness of $A_{0}$, the above bounds imply $%
u_{m}\in L^{2}(0,T;V^{\theta -\theta _{2}}),$ $\mu _{m}\in L^{2}\left(
0,T;L^{2}\left( \Omega \right) \right) ,$ with uniformly bounded norms. The
latter bound also yields from the third equation in (\ref{e:galerkin}), the
uniform bound 
\begin{equation}
\phi _{m}\in L^{2}\left( 0,T;D\left( A_{1}\right) \right) .
\label{e:exist-3bis}
\end{equation}%
Summarizing, $\left( u_{m},\phi _{m}\right) $ is uniformly bounded in $%
L^{\infty }(0,T;\mathcal{Y}_{\theta _{2}})\cap L^{2}(0,T;V^{\theta -\theta
_{2}}\times D\left( A_{1}\right) )$, and passing to a subsequence, there
exists $\left( u,\phi \right) \in L^{\infty }(0,T;\mathcal{Y}_{\theta
_{2}})\cap L^{2}(0,T;V^{\theta -\theta _{2}}\times D\left( A_{1}\right) )$
such that%
\begin{equation}
\begin{array}{l}
\left( u_{m},\phi _{m}\right) \rightarrow \left( u,\phi _{m}\right) 
\mbox{
weak-star in }L^{\infty }(0,T;\mathcal{Y}_{\theta _{2}}), \\ 
\left( u_{m},\phi _{m}\right) \rightarrow \left( u,\phi \right) 
\mbox{
weakly in }L^{2}(0,T;V^{\theta -\theta _{2}}\times D\left( A_{1}\right) ).%
\end{array}
\label{e:exist-weak-conv}
\end{equation}

As usual, for passing to the limit as $m\rightarrow \infty $ in (\ref%
{e:galerkin}), we will need a strong convergence result, which is obtained
by a standard compactness argument. We proceed by deriving bounds on the
derivatives of $u_{m}$ and $\phi _{m},$ respectively. Note that (\ref%
{e:galerkin}) can also be written as%
\begin{equation}
\left\{ 
\begin{array}{l}
u_{m}^{^{\prime }}+P_{m}A_{0}u_{m}+P_{m}B_{0}(u_{m},u_{m})=P_{m}\left(
g+R_{0}\left( \mu _{m},\phi _{m}\right) \right) , \\ 
\phi _{m}^{^{\prime }}+Q_{m}\mu _{m}+Q_{m}B_{1}\left( u_{m},\phi _{m}\right)
=0, \\ 
u_{m}(0)=P_{m}u(0),\phi _{m}\left( 0\right) =Q_{m}\phi \left( 0\right) .%
\end{array}%
\right.  \label{e:galerkin-abs}
\end{equation}%
Therefore,%
\begin{eqnarray}
\Vert \dot{u}_{m}\Vert _{-\gamma } &\lesssim &\Vert u_{m}\Vert _{\theta
-\theta _{2}}+\Vert P_{m}B_{0}(u_{m},u_{m})\Vert _{-\gamma }+\Vert
P_{m}g\Vert _{-\gamma }+\Vert P_{m}R_{0}\left( \mu _{m},\phi _{m}\right)
\Vert _{-\gamma }  \label{exist-1} \\
&=&:I_{1}+I_{2}+I_{3}+I_{4}.  \notag
\end{eqnarray}%
First, we notice that one has $I_{4}\lesssim \left\Vert \mu _{m}\right\Vert
_{L^{2}}\left\Vert \phi _{m}\right\Vert _{1}$ provided that $\gamma >\frac{n%
}{2}$. Therefore, $R_{0}\left( \mu _{m},\phi _{m}\right) \in L^{2}\left(
0,T;V^{-\gamma }\right) $ with uniform bound on account of (\ref{e:exist-3}%
)-(\ref{e:exist-weak-conv}). By the boundedness of $B_{0}$ (see (iv)), we
have as in \cite[Theorem 3.1]{HLT},%
\begin{equation}
I_{2}\lesssim \Vert u_{m}\Vert _{\sigma _{1}}\Vert u_{m}\Vert _{\sigma _{2}}.
\label{exist0}
\end{equation}%
If $\theta =0$, then the norms in the right hand side are the $V^{-\theta
_{2}}$-norm which is uniformly bounded. If $\theta >0$, since 
\begin{equation}
\Vert u_{m}\Vert _{\sigma _{i}}\lesssim \Vert u_{m}\Vert _{-\theta
_{2}}^{1-\lambda _{i}}\Vert u_{m}\Vert _{\theta -\theta _{2}}^{\lambda
_{i}},\qquad \lambda _{i}=\frac{\sigma _{i}+\theta _{2}}{\theta },\quad
i=1,2,  \label{exist1}
\end{equation}%
by the uniform boundedness of $u_{m}$ in $L^{\infty }(V^{-\theta _{2}})$ we
have 
\begin{equation}
I_{2}\lesssim \Vert u_{m}\Vert _{-\theta _{2}}^{2-\lambda _{1}-\lambda
_{2}}\Vert u_{m}\Vert _{\theta -\theta _{2}}^{\lambda _{1}+\lambda
_{2}}\lesssim \Vert u_{m}\Vert _{\theta -\theta _{2}}^{\lambda _{1}+\lambda
_{2}}.  \label{exist2}
\end{equation}%
Hence, with $\lambda :=\lambda _{1}+\lambda _{2}=\frac{\sigma _{1}+\sigma
_{2}+2\theta _{2}}{\theta }$ if $\theta >0$, and with $\lambda =1$ if $%
\theta =0$, we get 
\begin{align}
\Vert u_{m}^{^{\prime }}\Vert _{L^{p}(V^{-\gamma })}^{p}& \lesssim \Vert
u_{m}\Vert _{L^{p}(V^{\theta -\theta _{2}})}^{p}+\Vert u_{m}\Vert
_{L^{p\lambda }(V^{\theta -\theta _{2}})}^{p}+\Vert g\Vert
_{L^{p}(V^{-\theta -\theta _{2}})}^{p}  \label{exist3} \\
& +\Vert R_{0}\left( \mu _{m},\phi _{m}\right) \Vert _{L^{p}(V^{-\gamma
})}^{p}.  \notag
\end{align}%
The first and last terms on the right-hand side are bounded uniformly when $%
p\leq 2$. The second term is bounded if $p\lambda \leq 2$, that is $p\leq
2/\lambda $. We conclude that $u_{m}^{^{\prime }}$ is uniformly bounded in $%
L^{p}\left( V^{-\gamma }\right) $, with $p=\min \{2,2/\lambda \}$.
Concerning the time derivative of $\phi _{m},$ we use the second equation of
(\ref{e:galerkin-abs}) and the previous uniform bounds. We have%
\begin{align*}
\Vert \phi _{m}^{^{\prime }}||_{-2}& \lesssim \left( \left\Vert \mu
_{m}\right\Vert _{L^{2}}+\Vert Q_{m}B_{1}(u_{m},\phi _{m})\Vert _{-2}\right)
\\
& \lesssim \left( \left\Vert \mu _{m}\right\Vert _{L^{2}}+\Vert u_{m}\Vert
_{-\theta _{2}}\left\Vert \phi _{m}\right\Vert _{2}\right) ,
\end{align*}%
provided that $\theta _{2}\geq -1$. From this we can conclude that $\phi
_{m}^{^{\prime }}$ is uniformly bounded in $L^{2}\left( W^{-2}\right) .$ By
the Aubin-Lions-Simon compactness criterion (see, e.g., \cite{T}), we can
now obtain strong convergence properties for our sequence $\left( u_{m},\phi
_{m}\right) $, as follows. There exists%
\begin{equation*}
\left( u,\phi \right) \in C(0,T;V^{-\gamma }\times W^{0})\cap L^{\infty
}(0,T;\mathcal{Y}_{\theta _{2}})
\end{equation*}%
such that, in addition to (\ref{e:exist-weak-conv}), we also have%
\begin{equation}
\begin{array}{l}
\left( u_{m},\phi _{m}\right) \rightarrow \left( u,\phi \right) 
\mbox{
strongly in }L^{2}(0,T;V^{s}\times W^{2-}) \\ 
\phi _{m}\rightarrow \phi \mbox{ strongly in }C(0,T;W^{1-}),%
\end{array}
\label{e:exist-conv-u}
\end{equation}%
for any $s<\theta -\theta _{2},$ where $W^{s-}$ denotes $W^{s-\delta },$ for
some sufficiently small $\delta \in (0,s].$

Now we will show that this limit $\left( u,\phi \right) $ indeed satisfies
the weak formulation (\ref{weak1})-(\ref{weak2}). To this end, let $\left(
w,\psi \right) \in C^{\infty }(0,T;\mathcal{V}\times \mathcal{W})$ be an
arbitrary vector-valued function with $\left( w,\psi \right) (T)=\left(
0,0\right) $, and let $\left( w_{m},\psi _{m}\right) \in C^{1}(0,T;%
\widetilde{V}_{m}\times W_{m})$ be such that $\left( w_{m},\psi _{m}\right)
(T)=\left( 0,0\right) $ and%
\begin{align*}
w_{m}& \rightarrow w\text{ strongly in }C^{1}(0,T;V^{\bar{\gamma}}), \\
\psi _{m}& \rightarrow \psi \text{ strongly in }C^{1}\left( 0,T;W^{2}\right)
.
\end{align*}%
We have%
\begin{align}
& -\int_{0}^{T}\langle u_{m}(t),w_{m}^{^{\prime }}(t)\rangle
dt+\int_{0}^{T}\langle A_{0}u_{m}(t),w_{m}(t)\rangle dt  \label{e:convk} \\
& +\int_{0}^{T}b_{0}(u_{m}(t),u_{m}(t),w_{m}(t))dt-\int_{0}^{T}\left\langle
R_{0}(\mu _{m},\phi _{m}),w_{m}\left( t\right) \right\rangle dt  \notag \\
& =\langle u_{m}(0),w_{m}(0)\rangle +\int_{0}^{T}\langle
g(t),w_{m}(t)\rangle dt,  \notag
\end{align}%
and%
\begin{align}
& -\int_{0}^{T}\langle \phi _{m}(t),\psi _{m}^{^{\prime }}(t)\rangle
dt+\int_{0}^{T}\langle \mu _{m}(t),\psi _{m}(t)\rangle
dt+\int_{0}^{T}b_{1}\left( u_{m}\left( t\right) ,\phi _{m}\left( t\right)
,\psi _{m}\left( t\right) \right) dt  \label{e:convk2} \\
& =\langle \phi _{m}(0),\psi _{m}(0)\rangle .  \notag
\end{align}%
We would like to show that each term in the above equations converge to the
corresponding terms in (\ref{weak1})-(\ref{weak2}). The convergence in the
nonlinear term $b_{0}$ is shown in \cite[Theorem 3.1]{HLT}. Here we only
show it for the Korteweg term on the left-hand side of equation (\ref%
{e:convk}). The convergence in the nonlinear term $b_{1}$\ of (\ref{e:convk2}%
) is analogous. We have 
\begin{equation}
\int_{0}^{T}\left\vert \left\langle R_{0}(A_{1}\phi _{m},\phi
_{m}),w_{m}\left( t\right) \right\rangle -\left\langle R_{0}(A_{1}\phi ,\phi
),w\left( t\right) \right\rangle \right\vert dt\leq J_{m}+JJ_{m}+JJJ_{m},
\label{e:b-conv}
\end{equation}%
where the terms $J_{m}$, $JJ_{m}$, and $JJJ_{m}$ are defined below. Firstly,
it holds that%
\begin{equation}
\begin{split}
J_{m}& =\int_{0}^{T}\left\vert \left\langle R_{0}(A_{1}\phi _{m}(t),\phi
_{m}(t),w_{m}(t)-w(t))\right\rangle \right\vert dt \\
& \lesssim \int_{0}^{T}\Vert A_{1}\phi _{m}(t)\Vert _{L^{2}}\Vert \phi
_{m}(t)\Vert _{1}\Vert w_{m}(t)-w(t)\Vert _{\bar{\gamma}}dt \\
& \leq \Vert \phi _{m}\Vert _{L^{2}(D\left( A_{1}\right) )}\Vert \phi
_{m}\Vert _{L^{2}(W^{1})}\Vert w_{m}-w\Vert _{C(V^{\bar{\gamma}})}.
\end{split}%
\end{equation}%
thus, we get $\lim_{m\rightarrow \infty }J_{m}=0$. For $JJ_{m}$ we have%
\begin{equation}
\begin{split}
JJ_{m}& =\int_{0}^{T}\left\vert \left\langle R_{0}(A_{1}\left( \phi
_{m}(t)-\phi \left( t\right) \right) ,\phi _{m}(t),w(t))\right\rangle
\right\vert dt \\
& \lesssim \int_{0}^{T}\Vert \phi _{m}(t)-\phi (t)\Vert _{2-s}\Vert \phi
_{m}(t)\Vert _{2}\Vert w(t)\Vert _{\bar{\gamma}}dt \\
& \leq \Vert \phi _{m}-\phi \Vert _{L^{2}(W^{1})}\Vert \phi _{m}\Vert
_{L^{2}(W^{2})}\Vert w\Vert _{C(V^{\bar{\gamma}})},
\end{split}%
\end{equation}%
so $\lim_{m\rightarrow \infty }J_{m}=0$ by (\ref{e:exist-conv-u}) for as
long as $s\in \left( 0,2\right) $. Similarly, we have%
\begin{equation}
\begin{split}
JJJ_{m}& =\int_{0}^{T}\left\vert \left\langle R_{0}(A_{1}\phi \left(
t\right) ,\phi _{m}(t)-\phi \left( t\right) ,w(t))\right\rangle \right\vert
dt \\
& \lesssim \Vert A_{1}\phi \Vert _{L^{2}(L^{2})}\Vert \phi _{m}-\phi \Vert
_{L^{2}(W^{1})}\Vert w\Vert _{C(V^{\bar{\gamma}})},
\end{split}%
\end{equation}%
so $\lim_{m\rightarrow \infty }JJJ_{m}=0$ by (\ref{e:exist-conv-u}). The
proof of the theorem is now finished.
\end{proof}

Keeping in mind the Examples \ref{x:spaces}, \ref{y:spaces}, \ref{z:spaces}
from Section \ref{s:prelim}, our theorem covers the following special cases
listed in Table \ref{t:spec2}.

\begin{remark}
Let $\theta +\theta _{1}>\frac{1}{2}$ and recall that $\theta \geq 0$ and $%
\theta _{2}\geq -1$. By \cite[Proposition 2.5]{HLT}, the trilinear form $%
b_{00},$ defined by (\ref{b01})-(\ref{b01bis}), fulfills the hypotheses
(ii)-(iv) of Theorem \ref{t:exist} for $-\gamma \leq \theta -\theta _{2}-1$
with $-\gamma <\min \{2\theta +2\theta _{1}-\frac{n+2}{2},\theta -\theta
_{2}+2\theta _{1},\theta +\theta _{2}-1\}$. Similarly, the trilinear form $%
b_{01}$ satisfies (ii)-(iv) for $-\gamma \leq \theta -\theta _{2}-1$ with $%
-\gamma <\min \{2\theta +2\theta _{1}-\frac{n+2}{2},\theta -\theta
_{2}+2\theta _{1}-1,\theta +\theta _{2}\}$. In particular, our result yields
the global existence of a weak solution for both the inviscid and viscous
Leray-Allen-Cahn-$\alpha $ models in two and three space dimensions, and for
all the other regularized models listed in Table \ref{t:spec}.\ As far as we
know, except for the 3D NS-AC-$\alpha $-model \cite{GM1}, any of these
results have not been reported previously.
\end{remark}

\subsection{Uniqueness and stability}

\label{ss:stab}

Now we shall provide sufficient conditions for uniqueness and continuous
dependence with respect to the initial data of weak solutions of the general
three-parameter family of regularized models. Recall that $\theta _{1}\in 
\mathbb{R}$ and $\theta \geq 0.$

\begin{theorem}
\label{t:stab} Let $\left( u_{i},\phi _{i}\right) \in L^{\infty }(0,T;%
\mathcal{Y}_{\theta _{2}}),$ $i=1,2,$ be two solutions in the sense of
Definition \ref{weak}, corresponding to the initial conditions $\left(
u_{i}(0),\phi _{i}\left( 0\right) \right) \in \mathcal{Y}_{\theta _{2}}$, $%
i=1,2.$ Assume the following.

(i) $b_{0}:V^{\sigma _{1}}\times V^{\theta -\theta _{2}}\times V^{\sigma
_{2}}\rightarrow \mathbb{R}$ is bounded for some $\sigma _{1}\leq \theta
-\theta _{2}$ and $\sigma _{2}\leq \theta +\theta _{2}$ with $\sigma
_{1}+\sigma _{2}\leq \theta $.

(ii) $b_{0}(v,w,Nw)=0$ for any $v\in V^{\sigma _{1}}$ and $w\in V^{\sigma
_{2}}.$

We have the following cases:

(a) \underline{\textbf{Case }$\mathbf{n=2}$}: $\theta +\theta _{2}\geq 1$
and $\theta _{2}\geq 0.$

(b) \underline{\textbf{Case }$\mathbf{n=3}$}: $\theta _{2}\geq 1$ and $%
\theta \geq 0.$

Then the following estimate holds%
\begin{align}
& \Vert u_{1}(t)-u_{2}(t)\Vert _{-\theta _{2}}^{2}+\left\Vert \phi
_{1}\left( t\right) -\phi _{2}\left( t\right) \right\Vert _{1}^{2}
\label{uniq_stab} \\
& +\int_{0}^{t}\left( \Vert u_{1}\left( s\right) -u_{2}\left( s\right) \Vert
_{\theta -\theta _{2}}^{2}+\left\Vert A_{1}\left( \phi _{1}\left( s\right)
-\phi _{2}\left( s\right) \right) \right\Vert _{L^{2}}^{2}\right) ds  \notag
\\
& \leq \rho \left( t\right) \left( \Vert u_{1}(0)-u_{2}(0)\Vert _{-\theta
_{2}}^{2}+\left\Vert \phi _{1}\left( 0\right) -\phi _{2}\left( 0\right)
\right\Vert _{1}^{2}\right) ,  \notag
\end{align}%
for $t\in \lbrack 0,T],$ for some positive continuous function $\rho :%
\mathbb{R}_{+}\rightarrow \mathbb{R}_{+},$ $\rho \left( 0\right) >0$, which
depends only on the initial data $\left( u_{i}\left( 0\right) ,\phi
_{i}\left( 0\right) \right) $ in $\mathcal{Y}_{\theta _{2}}$-norm.
\end{theorem}

\begin{proof}
Let $v=u_{1}-u_{2}$ and $\psi =\phi _{1}-\phi _{2}$. Then subtracting the
equations for $\left( u_{1},\phi _{1}\right) $ and $\left( u_{2},\phi
_{2}\right) $ we have%
\begin{align}
& \left\langle \partial _{t}v,w\right\rangle +\left\langle
A_{0}v,w\right\rangle +\left\langle B_{0}\left( v,u_{1}\right)
,w\right\rangle +\left\langle B_{0}\left( u_{2},v\right) ,w\right\rangle
\label{diffuniq1} \\
& =\left\langle R_{0}\left( A_{1}\phi _{2},\psi \right) ,w\right\rangle
+\left\langle R_{0}\left( A_{1}\psi ,\phi _{1}\right) ,w\right\rangle , 
\notag
\end{align}%
and%
\begin{align}
& \left\langle \partial _{t}\psi ,\eta \right\rangle +\left\langle A_{1}\psi
,\eta \right\rangle +\left\langle B_{1}\left( v,\phi _{1}\right) ,\eta
\right\rangle +\left\langle B_{1}\left( u_{2},\psi \right) ,\eta
\right\rangle  \label{diffuniq2} \\
& =\left\langle f\left( \phi _{1}\right) -f\left( \phi _{2}\right) ,\eta
\right\rangle .  \notag
\end{align}%
First, observe that by the assumptions on $\theta ,\theta _{2},$ according
to (\ref{uniqest2})-(\ref{uniqest7}) below, the weak solution $\left(
u_{i},\phi _{i}\right) $\ of (\ref{e:op}) enjoys additional regularity.
Indeed, by \ virtue of (\ref{1.8}) it follows that $R_{0}\left( A_{1}\phi
_{i},\phi _{i}\right) \in L^{2}\left( 0,T;V^{-\theta -\theta _{2}}\right) $
and $B_{1}\left( u_{i},\phi _{i}\right) \in L^{2}\left( 0,T;L^{2}\left(
\Omega \right) \right) $. Thus, recalling the assumption on $b_{0}$ (see
(i)-(ii) above and (\ref{uniqest2}) below), it is easy to see that every
pairing $\left\langle \cdot ,\cdot \right\rangle $\ in (\ref{diffuniq1})-(%
\ref{diffuniq2}) is well-defined as a functional on the corresponding spaces
for $w\in L^{2}\left( 0,T;V^{\theta +\theta _{2}}\right) $ and $\eta \in
L^{2}\left( 0,T;L^{2}\left( \Omega \right) \right) $, respectively. Thus, by
(\ref{1.8}) we can take $w=Nv$ and $\eta =A_{1}\psi $ into (\ref{diffuniq1}%
)-(\ref{diffuniq2}) to infer%
\begin{align}
& \frac{d}{dt}\left( \Vert v\Vert _{-\theta _{2}}^{2}+\left\Vert
A_{1}^{1/2}\psi \right\Vert _{L^{2}}^{2}\right) +2c_{A_{0}}\Vert v\Vert
_{\theta -\theta _{2}}^{2}+2\left\Vert A_{1}\psi \right\Vert _{L^{2}}^{2}
\label{uniqest1} \\
& \leq 2\left\vert b_{0}(v,u_{1},Nv)\right\vert +2\left\vert b_{1}\left(
v,\psi ,A_{1}\phi _{2}\right) \right\vert +2\left\vert b_{1}\left(
u_{2},\psi ,A_{1}\psi \right) \right\vert +2\left\vert \left\langle f\left(
\phi _{1}\right) -f\left( \phi _{2}\right) ,A_{1}\psi \right\rangle
\right\vert .  \notag
\end{align}%
The first term on the right-hand side of (\ref{uniqest1}) can be bounded as
follows:%
\begin{align}
\left\vert b_{0}(v,u_{1},Nv)\right\vert & \lesssim \Vert v\Vert _{\sigma
_{1}}\Vert u_{1}\Vert _{\theta -\theta _{2}}\Vert v\Vert _{\sigma
_{2}-2\theta _{2}}  \label{uniqest2} \\
& \lesssim \Vert v\Vert _{-\theta _{2}}^{2-\lambda _{1}-\lambda _{2}}\Vert
v\Vert _{\theta -\theta _{2}}^{\lambda _{1}+\lambda _{2}}\Vert u_{1}\Vert
_{\theta -\theta _{2}}  \notag \\
& \lesssim \delta ^{-\frac{\lambda _{1}+\lambda _{2}}{2-\lambda _{1}-\lambda
_{2}}}\Vert v\Vert _{-\theta _{2}}^{2}\Vert u_{1}\Vert _{\theta -\theta
_{2}}^{2}+\delta \Vert v\Vert _{\theta -\theta _{2}}^{2},  \notag
\end{align}%
for any $\delta >0$, where $\lambda _{1}=\frac{\sigma _{1}+\theta _{2}}{%
\theta }$ and $\lambda _{2}=\frac{\sigma _{2}-\theta _{2}}{\theta }$, and
where in the last step we applied Young's inequality. Exploiting the fact
that $\phi _{i}\in L^{\infty }\left( 0,\infty ;L^{\infty }\left( \Omega
\right) \right) $ with $\left\vert \phi \right\vert \leq 1$ a.e. on $\Omega
\times \left( 0,\infty \right) $, the last term in (\ref{uniqest1}) is easy:%
\begin{equation}
\left\vert \left\langle f\left( \phi _{1}\right) -f\left( \phi _{2}\right)
,A_{1}\psi \right\rangle \right\vert \lesssim \delta \left\Vert A_{1}\psi
\right\Vert _{L^{2}}^{2}+\delta ^{-1}\left\Vert \psi \right\Vert
_{L^{2}}^{2}.  \label{uniqest3}
\end{equation}%
In order to bound the second and third terms on the right-hand side of (\ref%
{uniqest1}), we will treat each cases (a) and (b) separately.

\emph{Case (a)}: $n=2.$ Let $p\in \left( 2,\infty \right) .$ We have%
\begin{align*}
\left\vert b_{1}\left( u_{2},\psi ,A_{1}\psi \right) \right\vert &
=\left\vert \left\langle B_{1}\left( Nu_{2},\psi \right) ,A_{1}\psi
\right\rangle \right\vert \\
& \leq \left\Vert Nu_{2}\right\Vert _{L^{2p/\left( p-2\right) }}\left\Vert
A_{1}\psi \right\Vert _{L^{2}}^{2-\frac{2}{p}}\left\Vert \nabla \psi
\right\Vert _{L^{2}}^{\frac{2}{p}} \\
& \lesssim \delta ^{-1/\left( p-1\right) }\left\Vert Nu_{2}\right\Vert
_{L^{2p/\left( p-2\right) }}^{p}\left\Vert \nabla \psi \right\Vert
_{L^{2}}^{2}+\delta \left\Vert A_{1}\psi \right\Vert _{L^{2}}^{2} \\
& \lesssim \delta ^{-1/\left( p-1\right) }\left\Vert Nu_{2}\right\Vert
_{1}^{2}\left\Vert Nu_{2}\right\Vert _{1}^{\left( p-2\right) }\left\Vert
\psi \right\Vert _{1}^{2}+\delta \left\Vert A_{1}\psi \right\Vert
_{L^{2}}^{2},
\end{align*}%
where in the last step we have used the following inequality (see Appendix,
Lemma \ref{GNS}),%
\begin{equation}
\left\Vert w\right\Vert _{L^{2p/\left( p-2\right) }}\lesssim \left\Vert
w\right\Vert _{1}^{\frac{n}{2}-\frac{n\left( p-2\right) }{2p}}\left\Vert
w\right\Vert _{0}^{\frac{n\left( p-2\right) }{2p}+1-\frac{n}{2}},\text{ }%
n\geq 2.  \label{interp1}
\end{equation}%
Exploiting now the fact that $N:V^{s}\rightarrow V^{s+2\theta _{2}}$ is
bounded (see (\ref{e:bdd-amn})), and since $V^{\theta -\theta _{2}}\subseteq
V^{1-2\theta _{2}}$ and $V^{-\theta _{2}}\subseteq V^{-2\theta _{2}},$ one
has%
\begin{equation}
\left\vert b_{1}\left( u_{2},\psi ,A_{1}\psi \right) \right\vert \leq \delta
\left\Vert A_{1}\psi \right\Vert _{L^{2}}^{2}+\delta ^{-1/\left( p-1\right)
}\left\Vert u_{2}\right\Vert _{\theta -\theta _{2}}^{2}\left\Vert
u_{2}\right\Vert _{-\theta _{2}}^{\left( p-2\right) }\left\Vert \psi
\right\Vert _{1}^{2}.  \label{uniqest4}
\end{equation}%
For the last term, we have%
\begin{align}
\left\vert b_{1}\left( v,\psi ,A_{1}\phi _{2}\right) \right\vert &
=\left\vert \left\langle B_{1}\left( Nv,\psi \right) ,A_{1}\phi
_{2}\right\rangle \right\vert  \label{uniqest5} \\
& \leq \left\Vert Nv\right\Vert _{L^{4}}\left\Vert \nabla \psi \right\Vert
_{L^{4}}\left\Vert A_{1}\phi _{2}\right\Vert _{L^{2}}  \notag \\
& \lesssim \left( \left\Vert Nv\right\Vert _{0}^{1/2}\left\Vert
Nv\right\Vert _{1}^{1/2}\right) \left( \left\Vert \nabla \psi \right\Vert
_{L^{2}}^{1/2}\left\Vert A_{1}\psi \right\Vert _{L^{2}}^{1/2}\right)
\left\Vert A_{1}\phi _{2}\right\Vert _{L^{2}}  \notag \\
& \lesssim \left\Vert v\right\Vert _{-\theta _{2}}^{1/2}\left\Vert
v\right\Vert _{\theta -\theta _{2}}^{1/2}\left\Vert \psi \right\Vert
_{1}^{1/2}\left\Vert A_{1}\psi \right\Vert _{L^{2}}^{1/2}\left\Vert
A_{1}\phi _{2}\right\Vert _{L^{2}}  \notag \\
& \lesssim \delta \left\Vert v\right\Vert _{\theta -\theta _{2}}^{2}+\delta
\left\Vert A_{1}\psi \right\Vert _{L^{2}}^{2}+\delta ^{-1}\left( \left\Vert
v\right\Vert _{-\theta _{2}}\left\Vert \psi \right\Vert _{1}\right)
\left\Vert A_{1}\phi _{2}\right\Vert _{L^{2}}^{2}.  \notag
\end{align}

\emph{Case (b)}: $n=3.$ Let $p\in \left( 2,6\right) $ be fixed but otherwise
arbitrary (a suitable value will be chosen below). Once again, we have%
\begin{align*}
& \left\vert b_{1}\left( u_{2},\psi ,A_{1}\psi \right) \right\vert \\
& \leq \left\Vert Nu_{2}\right\Vert _{L^{2p/\left( p-2\right) }}\left\Vert
A_{1}\psi \right\Vert _{L^{2}}^{\frac{5}{2}-\frac{3}{p}}\left\Vert \nabla
\psi \right\Vert _{L^{2}}^{\frac{3}{p}-\frac{1}{2}} \\
& \lesssim \delta ^{-\frac{5p-6}{6-p}}\left\Vert Nu_{2}\right\Vert
_{L^{2p/\left( p-2\right) }}^{\frac{4p}{6-p}}\left\Vert \nabla \psi
\right\Vert _{L^{2}}^{2}+\delta \left\Vert A_{1}\psi \right\Vert _{L^{2}}^{2}
\\
& \lesssim \delta ^{-\frac{5p-6}{6-p}}\left\Vert Nu_{2}\right\Vert _{1}^{%
\frac{12}{6-p}}\left\Vert Nu_{2}\right\Vert _{1}^{\frac{4\left( p-3\right) }{%
6-p}}\left\Vert \psi \right\Vert _{1}^{2}+\delta \left\Vert A_{1}\psi
\right\Vert _{L^{2}}^{2},
\end{align*}%
where we have once more used (\ref{interp1}). By virtue of (\ref{e:bdd-amn}%
), and since $V^{\theta -\theta _{2}}\subseteq V^{-2\theta _{2}}$ and $%
V^{-\theta _{2}}\subseteq V^{1-2\theta _{2}},$ we can now conclude%
\begin{equation}
\left\vert b_{1}\left( u_{2},\psi ,A_{1}\psi \right) \right\vert \lesssim
\delta ^{-\frac{5p-6}{6-p}}\left\Vert u_{2}\right\Vert _{-\theta _{2}}^{%
\frac{12}{6-p}}\left\Vert u_{2}\right\Vert _{\theta -\theta _{2}}^{\frac{%
4\left( p-3\right) }{6-p}}\left\Vert \psi \right\Vert _{1}^{2}+\delta
\left\Vert A_{1}\psi \right\Vert _{L^{2}}^{2},  \label{uniqest6}
\end{equation}%
for any $\delta >0$, provided that we choose a suitable $p\in \left(
2,6\right) $ such that $\frac{4\left( p-3\right) }{6-p}\leq 2$ (this is
easily the case for any fixed $p\in \lbrack 3,4]$). For the last term, since 
$\theta _{2}\geq 1$ we have%
\begin{align}
\left\vert b_{1}\left( v,\psi ,A_{1}\phi _{2}\right) \right\vert & \leq
\left\Vert Nv\right\Vert _{\theta _{2}}\left\Vert \nabla \psi A_{1}\phi
\right\Vert _{-\theta _{2}}  \label{uniqest7} \\
& \lesssim \left\Vert v\right\Vert _{-\theta _{2}}\left\Vert A_{1}\phi
_{2}\right\Vert _{L^{2}}\left\Vert \nabla \psi \right\Vert _{L^{3}}  \notag
\\
& \lesssim \left\Vert v\right\Vert _{-\theta _{2}}\left\Vert A_{1}\phi
_{2}\right\Vert _{L^{2}}\left\Vert A_{1}\psi \right\Vert
_{L^{2}}^{1/2}\left\Vert \nabla \psi \right\Vert _{L^{2}}^{1/2}  \notag \\
& \lesssim \delta \left\Vert A_{1}\psi \right\Vert _{L^{2}}^{2}+\delta
^{-1/3}\left\Vert v\right\Vert _{-\theta _{2}}^{4/3}\left\Vert A_{1}\phi
_{2}\right\Vert _{L^{2}}^{4/3}\left\Vert \psi \right\Vert _{1}^{2/3}  \notag
\\
& \lesssim \delta \left\Vert A_{1}\psi \right\Vert _{L^{2}}^{2}+\delta
^{-1/3}\left( \left\Vert \psi \right\Vert _{1}^{2}+\left\Vert v\right\Vert
_{-\theta _{2}}^{2}\left\Vert A_{1}\phi _{2}\right\Vert _{L^{2}}^{2}\right) .
\notag
\end{align}

Collecting all estimates from (\ref{uniqest2}) to (\ref{uniqest7}), and
choosing a sufficiently small $\delta \sim \min \left(
c_{A_{0}},c_{A_{1}}\right) >0$, we can now apply Gr{\"{o}}nwall's inequality
in (\ref{uniqest1}) to deduce 
\begin{equation}
\Vert v(t)\Vert _{-\theta _{2}}^{2}+\left\Vert \psi \left( t\right)
\right\Vert _{1}^{2}\leq \left( \Vert v(0)\Vert _{-\theta
_{2}}^{2}+\left\Vert \psi \left( 0\right) \right\Vert _{1}^{2}\right) \exp
\int_{0}^{t}\Theta \left( s\right) ds,
\end{equation}%
for a suitable function $\Theta \in L^{1}\left( 0,T\right) .$ Integrating (%
\ref{uniqest1}) over $\left( 0,t\right) $ gives the desired inequality (\ref%
{uniq_stab}). The proof is finished.
\end{proof}

To clarify these results at least in the case of the specific models listed
in Table \ref{t:spec}, the corresponding conditions and stability results
derived from Theorem \ref{t:stab} are given below. Recall that both
conditions (a) and (b) of Theorem \ref{t:stab} are in force according to
whether $n=2$ or $n=3,$ respectively.

\begin{remark}
Exploiting \cite[Proposition 2.5]{HLT}, the trilinear form $b_{00}$
satisfies the hypotheses of Theorem \ref{t:stab} provided $\theta +\theta
_{1}\geq \frac{1-k}{2}$, $\theta +2\theta _{1}\geq k$, $\theta +\theta
_{2}\geq \frac{1}{2}$, $2\theta +2\theta _{1}+\theta _{2}>\frac{n+2}{2}$,
and $3\theta +2\theta _{1}+2\theta _{2}\geq 2-k$, for some $k\in \{0,1\}$.
The trilinear form $b_{01}$ satisfies the hypotheses of Theorem \ref{t:stab}
for $\theta +2\theta _{1}\geq 1$, $\theta +\theta _{1}\geq \frac{1}{2},$ $%
\theta +\theta _{2}\geq 0$, $2\theta +2\theta _{1}+\theta _{2}>\frac{n+2}{2}$%
, and $3\theta +2\theta _{1}+2\theta _{2}\geq 1$. For instance, these
assumptions allow us to recover the stability and uniqueness of the weak
solutions for the 3D Navier-Stokes-Allen-Cahn-$\alpha $-model included in
Table \ref{t:spec} (see also Table \ref{t:spec2}). This result was reported
previously in \cite{GM1}.
\end{remark}

\subsection{Regularity of weak solutions}

\label{ss:reg}

In this subsection, we develop a regularity result on weak solutions for the
general family of regularized models constructed in Section \ref{ss:exist}.
Incidently, the result below also allows us to obtain globally well-defined
(unique) \emph{strong} solutions for our regularized models, which will
become important in Section \ref{s:longbehav}\ to the study of the
asymptotic behavior as time goes to infinity. As in Section \ref{ss:stab},
recall that $\theta \geq 0$ and $\theta _{1}\in \mathbb{R}.$ As before, due
to the coupling of the regularized NSE with the Allen-Cahn equation, we will
separately derive optimal estimates in dimensions $n=2,3$.

\begin{theorem}
\label{t:reg} Let%
\begin{equation*}
\left( u,\phi \right) \in L^{\infty }\left( 0,T;\mathcal{Y}_{\theta
_{2}}\right) \cap L^{2}\left( 0,T;V^{\theta -\theta _{2}}\times D\left(
A_{1}\right) \right)
\end{equation*}%
be a weak solution in the sense of Definition \ref{weak}. Under the
assumptions of Theorems \ref{t:exist} and \ref{t:stab}, for some $\theta -%
\frac{1}{2}\geq \beta >-\theta _{2}$ when $n=3$, and some $\theta \geq \beta
\geq \max \left( 1-2\theta _{2},-\theta _{2}\right) $ with $\beta \neq
-\theta _{2}$ when $n=2$, let the following conditions hold.

(i) $b_{0}:V^{\alpha }\times V^{\alpha }\times V^{\theta -\beta }\rightarrow 
\mathbb{R}$ is bounded, where $\alpha =\min \{\beta ,\theta -\theta _{2}\}$;

(ii) $b_{0}(v,w,Nw)=0$ for any $v,w\in \mathcal{V}$;

(iii) $u_{0}\in V^{\beta }$, $\phi _{0}\in D\left( A_{1}\right) $, and $g\in
L^{2}(0,T;V^{\beta -\theta })$.

Then we have 
\begin{equation}
\left( u,\phi \right) \in L^{\infty }(0,T;V^{\beta }\times D\left(
A_{1}\right) )\cap L^{2}(0,T;V^{\beta +\theta }\times D(A_{1}^{3/2})).
\label{e:reg}
\end{equation}
\end{theorem}

\begin{proof}
The following estimates will be also deduced by a formal argument. However,
even in this case, they can be rigorously justified working with a
sufficiently smooth approximating solution, see Theorem \ref{t:exist}.
Pairing the first equation of (\ref{e:op}) with $\Lambda ^{2\beta }u,$ the
second and third equations with $A_{1}^{2}\phi ,$ respectively, we deduce%
\begin{align}
& \left\langle \partial _{t}u,\Lambda ^{2\beta }u\right\rangle +\left\langle
A_{0}u,\Lambda ^{2\beta }u\right\rangle +b_{0}(u,u,\Lambda ^{2\beta }u)
\label{e:innpro} \\
& =\left\langle g+R_{0}\left( A_{1}\phi ,\phi \right) ,\Lambda ^{2\beta
}u\right\rangle ,\qquad a.e.\ \text{ in }\left( 0,T\right) ,  \notag
\end{align}%
and%
\begin{equation}
\left\langle \partial _{t}\phi ,A_{1}^{2}\phi \right\rangle +\left\langle
A_{1}\phi ,A_{1}^{2}\phi \right\rangle +b_{1}(u,\phi ,A_{1}^{2}\phi
)=-\left\langle f\left( \phi \right) ,A_{1}^{2}\phi \right\rangle ,\text{ }%
a.e.\text{ in }\left( 0,T\right) .  \label{e:innpro2}
\end{equation}%
By employing the boundedness of $b_{0}$ (see (i)) and the coercivity of $%
A_{0}$, we infer%
\begin{equation}
b_{0}(u,u,\Lambda ^{2\beta }u)\lesssim \delta ^{-1}\Vert u\Vert _{\theta
-\theta _{2}}^{2}\Vert u\Vert _{\beta }^{2}+\delta \left\Vert u\right\Vert
_{\beta +\theta }^{2},\qquad a.e.\ \text{ in }\left( 0,T\right) ,
\label{est10}
\end{equation}%
for any $\delta >0$, while for the first term on the right-hand side of (\ref%
{e:innpro}), we have%
\begin{equation}
\left\langle g,\Lambda ^{2\beta }u\right\rangle \lesssim \delta \left\Vert
u\right\Vert _{\beta +\theta }^{2}+\delta ^{-1}\left\Vert g\right\Vert
_{\beta -\theta }^{2}.  \label{est11}
\end{equation}%
Moreover, since $\phi \in \left[ -1,1\right] $ a.e. on $\Omega \times \left(
0,T\right) $, one has%
\begin{equation}
\left\vert \left\langle f\left( \phi \right) ,A_{1}^{2}\phi \right\rangle
\right\vert =\left\vert \left\langle A_{1}^{1/2}f\left( \phi \right)
,A_{1}^{3/2}\phi \right\rangle \right\vert \lesssim \delta \left\Vert
A_{1}^{3/2}\phi \right\Vert _{L^{2}}^{2}+\delta ^{-1}\left\Vert \phi
\right\Vert _{1}^{2}.  \label{est12}
\end{equation}%
To bound the remaining terms in (\ref{e:innpro})-(\ref{e:innpro2}), we
divide our proof according to the different assumptions we employed for $n=3$
and $n=2.$

\underline{\textbf{Case }$\mathbf{n=3}$}: We begin to estimate the term
involving the Korteweg force. Since $V^{1/2}\subset L^{3}$, we deduce%
\begin{align}
\left\vert \left\langle R_{0}\left( A_{1}\phi ,\phi \right) ,\Lambda
^{2\beta }u\right\rangle \right\vert & \leq \left\Vert \Lambda ^{2\beta
}u\right\Vert _{\theta -\beta }\left\Vert R_{0}\left( A_{1}\phi ,\phi
\right) \right\Vert _{\beta -\theta }  \label{est13} \\
& \lesssim \left\Vert u\right\Vert _{\theta +\beta }\left\Vert R_{0}\left(
A_{1}\phi ,\phi \right) \right\Vert _{-1/2}  \notag \\
& \lesssim \left\Vert u\right\Vert _{\theta +\beta }\left\Vert A_{1}\phi
\right\Vert _{L^{2}}\left\Vert \nabla \phi \right\Vert _{L^{6}}  \notag \\
& \lesssim \delta \left\Vert u\right\Vert _{\beta +\theta }^{2}+\delta
^{-1}\left( \left\Vert A_{1}\phi \right\Vert _{L^{2}}^{2}\right) \left\Vert
A_{1}\phi \right\Vert _{L^{2}}^{2}.  \notag
\end{align}%
In order to estimate the term $b_{1}$ in (\ref{e:innpro2}) we require the
following basic inequality:%
\begin{equation}
\left\Vert \nabla ^{s}\left( f\cdot g\right) \right\Vert _{L^{p}}\lesssim
\left( \left\Vert f\right\Vert _{L^{p_{1}}}\left\Vert g\right\Vert
_{W^{s,p_{2}}}+\left\Vert f\right\Vert _{W^{s,p_{3}}}\left\Vert g\right\Vert
_{L^{p_{4}}}\right) ,  \label{basic00}
\end{equation}%
for any $p_{i}\in \left( 1,\infty \right) $, $i=1,...,4,$ with $%
1/p_{1}+1/p_{2}=1/p_{3}+1/p_{4}=1/p$ and $s\in \mathbb{N}_{0}.$ On account
of (\ref{basic00}), we have%
\begin{align*}
b_{1}(u,\phi ,A_{1}^{2}\phi )& =\left\langle A_{1}^{1/2}B_{1}\left( u,\phi
\right) ,A_{1}^{3/2}\phi \right\rangle \leq \left\Vert B_{1}\left( u,\phi
\right) \right\Vert _{1}\left\Vert A_{1}^{3/2}\phi \right\Vert _{L^{2}} \\
& \lesssim \left( \left\Vert Nu\right\Vert _{L^{6}}\left\Vert \nabla \phi
\right\Vert _{W^{1,3}}+\left\Vert Nu\right\Vert _{1}\left\Vert \nabla \phi
\right\Vert _{L^{\infty }}\right) \left\Vert A_{1}^{3/2}\phi \right\Vert
_{L^{2}} \\
& =:I_{1}+I_{2}.
\end{align*}%
To estimate the $I_{1}$-term, we recall that $N:V^{1-2\theta
_{2}}\rightarrow V^{1}$ is bounded, and $V^{-\theta _{2}}\subseteq
V^{1-2\theta _{2}}.$ Indeed, exploiting the following inequality $\left\Vert
\cdot \right\Vert _{W^{1,3}}\lesssim \left\Vert \cdot \right\Vert
_{2}^{3/4}\left\Vert \cdot \right\Vert _{L^{2}}^{1/4}$ (see Lemma \ref{GNS}%
), we infer%
\begin{align}
I_{1}& \lesssim \left\Vert u\right\Vert _{-\theta _{2}}\left\Vert \phi
\right\Vert _{1}^{1/4}\left\Vert A_{1}^{3/2}\phi \right\Vert _{L^{2}}^{7/4}
\label{est14} \\
& \lesssim \delta \left\Vert A_{1}^{3/2}\phi \right\Vert _{L^{2}}^{2}+\delta
^{-7}\left\Vert u\right\Vert _{-\theta _{2}}^{8}\left\Vert \phi \right\Vert
_{1}^{2}.  \notag
\end{align}%
Similarly, by the 3D Agmon's inequality, we have%
\begin{align}
I_{2}& \lesssim \left\Vert u\right\Vert _{-\theta _{2}}\left\Vert A_{1}\phi
\right\Vert _{L^{2}}^{1/2}\left\Vert A_{1}^{3/2}\phi \right\Vert
_{L^{2}}^{3/2}  \label{est15} \\
& \lesssim \delta \left\Vert A_{1}^{3/2}\phi \right\Vert _{L^{2}}^{2}+\delta
^{-3}\left\Vert u\right\Vert _{-\theta _{2}}^{4}\left\Vert A_{1}\phi
\right\Vert _{L^{2}}^{2}.  \notag
\end{align}%
Collecting all the above estimates from (\ref{est10}) to (\ref{est15}), and
then adding together the relations (\ref{e:innpro})-(\ref{e:innpro}), for a
sufficiently small $\delta \sim \min \left( c_{A_{1}},c_{A_{0}}\right) >0$,
we deduce the following inequality%
\begin{align}
& \frac{d}{dt}\left( \left\Vert u\right\Vert _{\beta }^{2}+\left\Vert
A_{1}\phi \right\Vert _{L^{2}}^{2}\right) +\left\Vert A_{1}^{3/2}\phi
\right\Vert _{L^{2}}^{2}+c_{A_{0}}\left\Vert u\right\Vert _{\beta +\theta
}^{2}  \label{est15bis} \\
& \leq \Xi \left( t\right) \left( \left\Vert u\right\Vert _{\beta
}^{2}+\left\Vert A_{1}\phi \right\Vert _{L^{2}}^{2}\right) +C\left\Vert
g\right\Vert _{\beta -\theta }^{2},  \notag
\end{align}%
where we have set%
\begin{equation}
\Xi :=C_{\delta }\left( 1+\left\Vert u\right\Vert _{-\theta
_{2}}^{8}+\left\Vert \nabla \phi \right\Vert _{L^{2}}^{2}+\left\Vert
A_{1}\phi \right\Vert _{L^{2}}^{2}+\left\Vert u\right\Vert _{\theta -\theta
_{2}}^{2}\right) ,  \label{f15}
\end{equation}%
for some $C_{\delta }>0$. Notice that since $\left( u,\phi \right) $ is a
weak solution in the sense of Definition \ref{weak}, we have $\Xi \in
L^{1}\left( 0,T\right) $. Integrating (\ref{est15bis}) over $\left(
0,t\right) $, and using Gr\"{o}nwall's inequality, we conclude 
\begin{align}
& \left\Vert u\left( t\right) \right\Vert _{\beta }^{2}+\left\Vert A_{1}\phi
\left( t\right) \right\Vert _{L^{2}}^{2}  \label{est15tris} \\
& \lesssim \left( \int_{0}^{t}\Vert g\Vert _{\beta -\theta }^{2}+\Vert
u(0)\Vert _{\beta }^{2}+\left\Vert A_{1}\phi \left( 0\right) \right\Vert
_{L^{2}}^{2}\right) \exp \left( 2\int_{0}^{t}\Xi \left( s\right) ds\right)
,\qquad a.e.\ \text{ in }\left( 0,T\right) .  \notag
\end{align}%
Therefore, by assumption (iii) we have $\left( u,\phi \right) \in L^{\infty
}(0,T;V^{\beta }\times D\left( A_{1}\right) )$, which transfers to $\left(
u,\phi \right) \in L^{2}(0,T;V^{\theta +\beta }\times D(A_{1}^{3/2}))$ on
account of (\ref{est15bis})-(\ref{est15tris}).

\underline{\textbf{Case }$\mathbf{n=2}$}: Concerning the Korteweg term,
using the Agmon's inequality in two dimensions, we have%
\begin{align}
\left\vert \left\langle R_{0}\left( A_{1}\phi ,\phi \right) ,\Lambda
^{2\beta }u\right\rangle \right\vert & \leq \left\Vert \Lambda ^{2\beta
}u\right\Vert _{\theta -\beta }\left\Vert R_{0}\left( A_{1}\phi ,\phi
\right) \right\Vert _{\beta -\theta }  \label{est13bis} \\
& \lesssim \left\Vert u\right\Vert _{\theta +\beta }\left\Vert R_{0}\left(
A_{1}\phi ,\phi \right) \right\Vert _{0}  \notag \\
& \lesssim \delta \left\Vert u\right\Vert _{\beta +\theta }^{2}+\delta
^{-1}\left\Vert \nabla \phi \right\Vert _{L^{\infty }}^{2}\left\Vert
A_{1}\phi \right\Vert _{L^{2}}^{2}  \notag \\
& \lesssim \delta \left\Vert u\right\Vert _{\beta +\theta }^{2}+\delta
\left\Vert A_{1}^{3/2}\phi \right\Vert _{L^{2}}^{2}  \notag \\
& +\delta ^{-2}\left\Vert A_{1}\phi \right\Vert _{L^{2}}^{2}\left(
\left\Vert A_{1}\phi \right\Vert _{L^{2}}^{2}\left\Vert \phi \right\Vert
_{1}^{2}\right) ,  \notag
\end{align}%
for any $\delta >0$. Next, using (\ref{basic00}) we have%
\begin{align*}
b_{1}(u,\phi ,A_{1}^{2}\phi )& \leq \left\Vert B_{1}\left( u,\phi \right)
\right\Vert _{1}\left\Vert A_{1}^{3/2}\phi \right\Vert _{L^{2}} \\
& \lesssim \left( \left\Vert Nu\right\Vert _{W^{1,2}}\left\Vert \nabla \phi
\right\Vert _{L^{\infty }}+\left\Vert Nu\right\Vert _{L^{4}}\left\Vert
\nabla \phi \right\Vert _{W^{1,4}}\right) \left\Vert A_{1}^{3/2}\phi
\right\Vert _{L^{2}} \\
& =:I_{1}+I_{2}.
\end{align*}%
Now, since $N$ is bounded from $V^{-2\theta _{2}}\rightarrow V^{0}$ and from 
$V^{1-2\theta _{2}}\rightarrow V^{1}$, the second term%
\begin{align}
I_{2}& \lesssim \left\Vert Nu\right\Vert _{0}^{1/2}\left\Vert Nu\right\Vert
_{1}^{1/2}\left\Vert A_{1}\phi \right\Vert _{L^{2}}^{1/2}\left\Vert
A_{1}^{3/2}\phi \right\Vert _{L^{2}}^{3/2}  \label{est14bis} \\
& \lesssim \delta \left\Vert A_{1}^{3/2}\phi \right\Vert _{L^{2}}^{2}+\delta
^{-3}\left\Vert u\right\Vert _{-\theta _{2}}^{2}\left\Vert u\right\Vert
_{\theta -\theta _{2}}^{2}\left\Vert A_{1}\phi \right\Vert _{L^{2}}^{2}. 
\notag
\end{align}%
Similarly, by elementary inequalities and the facts that $N$ is bounded and $%
V^{\theta -\theta _{2}}\subseteq V^{1-2\theta _{2}},$ we have%
\begin{align}
I_{1}& \leq C\left\Vert \nabla \phi \right\Vert _{L^{2}}^{1/2}\left\Vert
A_{1}^{3/2}\phi \right\Vert _{L^{2}}^{3/2}\left\Vert Nu\right\Vert _{1}
\label{est15q} \\
& \lesssim \delta ^{-3}\left\Vert u\right\Vert _{\theta -\theta
_{2}}^{2}\left\Vert u\right\Vert _{\beta }^{2}\left\Vert \nabla \phi
\right\Vert _{L^{2}}^{2}+\delta \left\Vert A_{1}^{3/2}\phi \right\Vert
_{L^{2}}^{2},  \notag
\end{align}%
provided that $\beta \geq \max \left( 1-2\theta _{2},-\theta _{2}\right) .$
Collecting the above estimates (\ref{est13bis})-(\ref{est15q}) and insert
them into the right hand sides of (\ref{e:innpro})-(\ref{e:innpro}), for a
sufficiently small $\delta \sim \min \left( c_{A_{0}},c_{A_{1}}\right) $, we
infer%
\begin{align}
& \frac{d}{dt}\left( \left\Vert u\right\Vert _{\beta }^{2}+\left\Vert
A_{1}\phi \right\Vert _{L^{2}}^{2}\right) +\left\Vert A_{1}^{3/2}\phi
\right\Vert _{L^{2}}^{2}+c_{A_{0}}\left\Vert u\right\Vert _{\beta +\theta
}^{2}  \label{est16} \\
& \leq \Psi \left( t\right) \left( \left\Vert u\right\Vert _{\beta
}^{2}+\left\Vert A_{1}\phi \right\Vert _{L^{2}}^{2}\right) +C\left\Vert
g\right\Vert _{\beta -\theta }^{2},  \notag
\end{align}%
where we have set%
\begin{equation}
\Psi :=C_{\delta }\left[ 1+\Vert u\Vert _{\theta -\theta _{2}}^{2}\left(
\left\Vert \nabla \phi \right\Vert _{L^{2}}^{2}+\left\Vert u\right\Vert
_{-\theta _{2}}^{2}\right) +\left\Vert A_{1}\phi \right\Vert
_{L^{2}}^{2}\left\Vert \phi \right\Vert _{1}^{2}\right] ,  \label{f16}
\end{equation}%
for some $C_{\delta }>0$. We remark once again that $\Psi \in L^{1}\left(
0,T\right) $ for any weak solution $\left( u,\phi \right) $ to problem (\ref%
{e:op}). Thus, the application of Gronwall's inequality in (\ref{est16})
yields the desired conclusion. The proof of the theorem is finished.
\end{proof}

\begin{remark}
\label{reg-rem}Let $4\theta +4\theta _{1}+2\theta _{2}>n+2$, $2\theta
+2\theta _{1}\geq 1-k$, $\theta +2\theta _{2}\geq 1$, $3\theta +4\theta
_{1}\geq 1$, $\theta +2\theta _{1}\geq \ell $, and $3\theta +2\theta
_{1}+2\theta _{2}\geq 2-\ell $, for some $k,\ell \in \{0,1\}$. In addition
to the assumptions of Theorem \ref{t:reg}, let 
\begin{equation*}
\beta \in (\frac{n+2}{2}-2(\theta _{1}+\theta _{2})-\theta ,3\theta +2\theta
_{1}-\frac{n+2}{2})\cap \lbrack \frac{1-\ell }{2}-\theta _{1}-\theta
_{2},\min \{2\theta +\theta _{2}-1,2\theta -\theta _{2}+2\theta _{1}-k\}].
\end{equation*}%
Then, by \cite[Proposition 2.5]{HLT} the trilinear form $b_{00}$ satisfies
the hypotheses of the above theorem.
\end{remark}

\begin{remark}
Let $4\theta +4\theta _{1}+2\theta _{2}>n+2$, $\theta +2\theta _{2}\geq 0$,
and $\theta +2\theta _{1}\geq 1$. Let 
\begin{equation*}
\beta \in (\frac{n+2}{2}-2(\theta _{1}+\theta _{2})-\theta ,3\theta +2\theta
_{1}-\frac{n+2}{2})\cap \lbrack \frac{1}{2}-\theta _{1}-\theta _{2},\min
\{2\theta +\theta _{2},2\theta -\theta _{2}+2\theta _{1}-1\}].
\end{equation*}%
By \cite[Proposition 2.5]{HLT}, the trilinear form $b_{01}$ satisfies the
hypotheses of the above theorem.
\end{remark}

The corresponding conditions and results of Theorem \ref{t:reg} above are
listed in the Table \ref{t:spec-reg} below for the \emph{three-dimensional}
regularized models listed in Table \ref{t:spec}. For the NS-AC-$\alpha $%
-like model in Table \ref{t:spec}, the allowed values for $\beta $ are $%
\beta \leq 2\theta -\theta _{2}-1$ with $\beta <3\theta -\frac{5}{2}$ and $%
\beta \leq \theta -\frac{1}{2},$ provided that $\theta \geq \frac{1}{2}$ and 
$4\theta +2\theta _{2}>5$.{\small 
\begin{table}[th]
\caption{Regularity results for some special cases of the model (\protect\ref%
{e:pde}) in \emph{three} space dimensions. The table gives values of $%
\protect\beta $ for our recovered \emph{global} regularity results for the
regularized two-phase flow models, where $u_{0}\in V^{\protect\beta }$ and $%
\protect\phi _{0}\in D\left( A_{1}\right) $. (In the NS--AC-$\protect\alpha $%
-like case, see the text for the allowed values of $\protect\beta $.) }
\label{t:spec-reg}
\begin{center}
{\small $%
\begin{tabular}{||l||l||l||l||l||l||l||}
\hline\hline
3D Model & NSE-AC & Leray-AC-$\alpha $ & ML-AC-$\alpha $ & SBM-AC & NSV-AC & 
NS-AC-$\alpha $ \\ \hline\hline
Regularity & N/A & $\beta \in (0,\frac{1}{2}]$ & $\beta \in (-1,\frac{1}{2}]$
& $\beta \in (-1,\frac{1}{2}]$ & $\beta \in (-1,-\frac{1}{2}]$ & $\beta \in
(-1,0]$ \\ \hline\hline
\end{tabular}%
$  }
\end{center}
\end{table}
}

\section{Singular perturbations}

\label{s:pert}

In this section, following \cite{HLT} we will consider the situation where
the operators $A_{0}$ and $B_{0}$ in the general three-parameter family of
regularized models represented by problem (\ref{e:op}) have values from a
convergent (in a certain sense) sequence, and study the limiting behavior of
the corresponding sequence of solutions. As special cases we have inviscid
limits $\left( \nu =0\right) $ in the viscous equations and $\alpha
\rightarrow 0^{+}$ limits in the $\alpha $-models.

\subsection{Perturbations to the linear part of the flow component}

\label{ss:pert-lin}

Consider the problem%
\begin{equation}
\left\{ 
\begin{array}{l}
\partial _{t}u+A_{0}u+B_{0}(u,u)=g+R_{0}\left( A_{1}\phi ,\phi \right) , \\ 
\partial _{t}\phi +A_{1}\phi +B_{1}\left( u,\phi \right) +f\left( \phi
\right) =0%
\end{array}%
\right.  \label{e:pert-lin-0}
\end{equation}%
and its perturbation%
\begin{equation}
\left\{ 
\begin{array}{l}
\partial _{t}u_{i}+A_{0i}u_{i}+B_{0}(u_{i},u_{i})=g+R_{0}\left( A_{1}\phi
_{i},\phi _{i}\right) , \\ 
\partial _{t}\phi _{i}+A_{1}\phi _{i}+B_{1}\left( u_{i},\phi _{i}\right)
+f\left( \phi _{i}\right) =0,%
\end{array}%
\right.  \label{e:pert-lin}
\end{equation}%
for $i\in \mathbb{N}$, where $A_{0}$, $B_{0}$, $B_{1},$ $R_{0\text{ }}$and $%
N $ satisfy the assumptions stated in Section \ref{s:prelim}, and for $i\in 
\mathbb{N}$, $A_{0i}:V^{s}\rightarrow V^{s-2\varepsilon }$ is a bounded
linear operator satisfying 
\begin{equation}
\Vert A_{0i}v\Vert _{-\varepsilon -\theta _{2}}^{2}+\Vert v\Vert _{\theta
-\theta _{2}}^{2}\lesssim \left\langle A_{0i}v,Nv\right\rangle +\Vert v\Vert
_{-\theta _{2}}^{2},\quad v\in V^{\varepsilon -\theta _{2}}.
\label{e:pert-lin-ell}
\end{equation}%
Assuming that both problems (\ref{e:pert-lin-0}) and (\ref{e:pert-lin}) have
the same initial condition $\left( u_{0},\phi _{0}\right) $, and that $%
A_{0i}\rightarrow A_{0}$ in some topology, we are concerned with the
behavior of $\left( u_{i},\phi _{i}\right) $ as $i\rightarrow \infty $. We
will also assume that $\varepsilon \geq \theta $.

\begin{theorem}
\label{t:pert-lin} Assume the above setting, and in addition let the
following conditions hold.

\begin{itemize}
\item[i)] $\left( u_{0},\phi _{0}\right) \in \mathcal{Y}_{\theta _{2}}$,
with any $\theta _{2}\geq -1$, and $g\in L^{2}(0,T;V^{-\theta -\theta _{2}})$%
, $T>0$;

\item[ii)] $b_{0}(v,v,Nv)=0$ for any $v\in \mathcal{V}$;

\item[iii)] $b_{0}:V^{\sigma _{1}}\times V^{\sigma _{2}}\times V^{\gamma
}\rightarrow \mathbb{R}$ is bounded for some $\sigma _{j}\in \lbrack -\theta
_{2},\theta -\theta _{2}]$, $j=1,2$, and $\gamma \in \lbrack \varepsilon
+\theta _{2},\infty )\cap (\theta _{2},\infty )\cap \left( \frac{n}{2}%
,\infty \right) $;

\item[iv)] $b_{0}:V^{\bar{\sigma}_{1}}\times V^{\bar{\sigma}_{2}}\times V^{%
\bar{\gamma}}\rightarrow \mathbb{R}$ is bounded for some $\bar{\sigma}%
_{j}<\theta -\theta _{2}$, $j=1,2$, and $\bar{\gamma}\geq \gamma $;

\item[v)] $A_{0i}$ converge weakly to $A_{0}$ as $i\rightarrow \infty $.
\end{itemize}

Then, there exists a solution $\left( u,\phi \right) \in L^{\infty }(0,T;%
\mathcal{Y}_{\theta _{2}})\cap L^{2}(0,T;V^{\theta -\theta _{2}}\times
D\left( A_{1}\right) )$ to (\ref{e:pert-lin-0}) such that, up to a
subsequence,%
\begin{equation}
\left\{ 
\begin{array}{l}
\left( u_{i},\phi _{i}\right) \rightarrow \left( u,\phi \right) 
\mbox{
weak-star in }L^{\infty }(0,T;\mathcal{Y}_{\theta _{2}}), \\ 
\left( u_{i},\phi _{i}\right) \rightarrow \left( u,\phi \right) 
\mbox{
weakly in }L^{2}(0,T;V^{\theta -\theta _{2}}\times D\left( A_{1}\right) ),
\\ 
\left( u_{i},\phi _{i}\right) \rightarrow \left( u,\phi \right) 
\mbox{
strongly in }L^{2}(0,T;V^{s}\times W^{l}), \\ 
\phi _{i}\rightarrow \phi \mbox{
strongly in }C(0,T;W^{\zeta })\mbox{ for any }\zeta <1,%
\end{array}%
\right.  \label{e:pert-conv-u}
\end{equation}%
for any $s<\theta -\theta _{2},$ $l<2,$ as $i\rightarrow \infty $.
\end{theorem}

\begin{proof}
First, from Theorem~\ref{t:exist}, we know that for $i\in \mathbb{N}$ there
exists a solution%
\begin{equation*}
\left( u_{i},\phi _{i}\right) \in L^{\infty }(0,T;\mathcal{Y}_{\theta
_{2}})\cap L^{2}(0,T;V^{\varepsilon -\theta _{2}}\times D\left( A_{1}\right)
)
\end{equation*}%
to (\ref{e:pert-lin}). Moreover, $\phi _{i}\in L^{\infty }\left(
0,T;L^{\infty }\left( \Omega \right) \right) $ with $\phi _{i}\in \left[ -1,1%
\right] $ for all $i\in \mathbb{N}.$ Duality pairing of the first and second
equations of (\ref{e:pert-lin}) with $Nu_{i}$ and $A_{1}\phi _{i}+f\left(
\phi _{i}\right) $, respectively, arguing as in the proof of Theorem \ref%
{t:exist}, we deduce%
\begin{equation}
\begin{split}
& \frac{d}{dt}\mathcal{E}\left( u_{i}\left( t\right) ,\phi _{i}\left(
t\right) \right) +2\langle A_{0i}u_{i},Nu_{i}\rangle +\left\Vert A_{1}\phi
_{i}+f\left( \phi _{i}\right) \right\Vert _{L^{2}}^{2} \\
& \lesssim \delta ^{-1}\Vert g\Vert _{-\theta -\theta _{2}}^{2}+\delta \Vert
u_{i}\Vert _{\theta -\theta _{2}}^{2},
\end{split}
\label{e:pert-1}
\end{equation}%
Choosing $\delta >0$ small enough, then using (\ref{e:pert-lin-ell}) and
integrating over $\left( 0,t\right) $ we have 
\begin{equation}
\mathcal{E}\left( u_{i}\left( t\right) ,\phi _{i}\left( t\right) \right)
\lesssim C_{T}\left( 1+\mathcal{E}\left( u\left( 0\right) ,\phi \left(
0\right) \right) \right) .  \label{e:pert-2}
\end{equation}%
Moreover, integrating (\ref{e:pert-1}), and taking into account (\ref%
{e:pert-2}) and (\ref{e:pert-lin-ell}), we once again infer%
\begin{equation}
\Vert A_{0i}u_{i}\Vert _{L^{2}(0,t;V^{-\varepsilon -\theta
_{2}})}^{2}+\left\Vert A_{1}\phi _{i}+f\left( \phi _{i}\right) \right\Vert
_{L^{2}\left( 0,t;L^{2}\right) }^{2}+\Vert u_{i}\Vert _{L^{2}(0,t;V^{\theta
-\theta _{2}})}^{2}\leq C_{T},\qquad t\in \left( 0,T\right) .
\label{e:pert-3}
\end{equation}%
For any fixed $T>0$, this gives $u_{i}\in L^{\infty }(0,T;V^{-\theta
_{2}})\cap L^{2}(0,T;V^{\theta -\theta _{2}})$ and $\phi _{i}\in L^{2}\left(
0,T;D\left( A_{1}\right) \right) $, respectively, with uniformly (in $i$)
bounded norms. On the other hand, we have%
\begin{align}
\Vert u_{i}^{^{\prime }}\Vert _{-\gamma }& \leq \Vert A_{0i}u_{i}\Vert
_{-\gamma }+\Vert B_{0}(u_{i},u_{i})\Vert _{-\gamma }+\Vert g+R_{0}\left(
A_{1}\phi _{i},\phi _{i}\right) \Vert _{-\gamma }, \\
||\phi _{i}^{^{\prime }}||_{-2}& \leq \left\Vert A_{1}\phi _{i}+f\left( \phi
_{i}\right) \right\Vert _{L^{2}}^{2}+\left\Vert B_{1}\left( u_{i},\phi
_{i}\right) \right\Vert _{-2}.  \notag
\end{align}%
By estimating these terms as in the proof of Theorem \ref{t:exist}, and
taking into account (\ref{e:pert-2})-(\ref{e:pert-3}), we conclude that $%
(u_{i}^{^{\prime }},\phi _{i}^{^{\prime }})$ is uniformly bounded in $%
L^{2}(0,T;V^{-\gamma }\times W^{-2})$. Passing now to a subsequence owing to
compactness arguments, we infer the existence of $\left( u,\phi \right) $
satisfying (\ref{e:pert-conv-u}). Now taking into account the weak
convergence of $A_{0i}$ to $A_{0}$, the rest of the proof proceeds similarly
to that of Theorem \ref{t:exist}.
\end{proof}

For example, setting $\varepsilon =1$, with $\theta =0$ and $\theta _{2}=1$,
and checking all the requirements (i)-(v) of Theorem~\ref{t:pert-lin}, the
viscous solutions to the 3D SBM-AC model converge to the inviscid solution
as the viscosity $\nu $ tends to zero. Recall that the global existence of a
weak solution to the inviscid 3D SBM-AC model (see Table \ref{t:spec}) is
also contained in Theorem~\ref{t:exist}. Similarly, setting $\varepsilon =0$%
, with $\theta =0$ and $\theta _{2}=1$, the viscous solutions to the 3D
Leray--AC-$\alpha $ model converge to the inviscid solution as the viscosity
tends to zero. This result gives another proof of the global existence of a
weak solution for the inviscid ($\nu =0$) 3D Leray--AC-$\alpha $ model. Note
that as in \cite[Section 4, Theorem 4.1]{HLT} where the inviscid problem for
single-fluids was investigated, Theorem~\ref{t:pert-lin} establishes
analogue global existence results for homogeneous two-phase Allen-Cahn flows
without any viscosity terms. In particular, by Theorem \ref{t:pert-lin} we
recover the convergence results as $\nu \rightarrow 0$ for the global weak
solutions of NSE-AC system in two space dimensions, which were previously
reported in \cite{X, ZGH}.

\subsection{Perturbations involving the nonlinear part of the flow component}

\label{ss:pert-nonlin}

We employ the same assumptions as in \cite[Section 4]{HLT}. For $i\in 
\mathbb{N}$, let $A_{0i}:V^{s}\rightarrow V^{s-2\varepsilon }$ and $%
N_{i}:V^{s}\rightarrow V^{s+2\varepsilon _{2}}$ be bounded linear operators,
satisfying 
\begin{equation}
\Vert v\Vert _{\theta +\theta _{2}}^{2}\lesssim \left\langle
A_{0i}N_{i}^{-1}v,v\right\rangle +\Vert v\Vert _{\theta _{2}}^{2},\qquad
v\in V^{\theta +\theta _{2}},  \label{e:pert-nl-ell}
\end{equation}%
and 
\begin{equation}
\Vert v\Vert _{\theta _{2}}^{2}\lesssim \left\langle
N_{i}^{-1}v,v\right\rangle ,\qquad v\in V^{\theta _{2}},
\label{e:pert-nl-nell}
\end{equation}%
where we also assumed that $N_{i}$ is invertible. In this subsection, we
continue with perturbations of (\ref{e:pert-lin-0}) of the form%
\begin{equation}
\left\{ 
\begin{array}{l}
\partial _{t}u_{i}+A_{0i}u_{i}+B_{0i}(u_{i},u_{i})=g+R_{0}\left( A_{1}\phi
_{i},\phi _{i}\right) , \\ 
\partial _{t}\phi _{i}+A_{1}\phi _{i}+B_{1}\left( u_{i},\phi _{i}\right)
+f\left( \phi _{i}\right) =0,\text{ }(i\in \mathbb{N}),%
\end{array}%
\right.  \label{e:pert-nonlin}
\end{equation}%
where $B_{0i}$ is some bilinear map. Again assuming that both problems (\ref%
{e:pert-lin-0}) and (\ref{e:pert-nonlin}) have the same initial condition $%
\left( u_{0},\phi _{0}\right) $, and that $A_{0i}\rightarrow A_{0}$ and $%
B_{0i}\rightarrow B_{0}$ in some topology, we are concerned with the
behavior of $\left( u_{i},\phi _{i}\right) $ as $i\rightarrow \infty $. For
reference, define the trilinear form $b_{0i}(u,v,w)=\left\langle
B_{0i}\left( u,v\right) ,w\right\rangle $.

\begin{theorem}
\label{t:pert-nonlin} Assume the above setting, and in addition let the
following conditions hold.

\begin{itemize}
\item[i)] $\left( u_{0},\phi _{0}\right) \in Y_{\theta _{2}}$ with any $%
\theta _{2}\geq -1$, and $g\in L^{2}(0,T;V^{-\theta -\theta _{2}})$, $T>0$;

\item[ii)] $b_{0i}(v,v,N_{i}v)=0$ for any $v\in \mathcal{V}$;

\item[iii)] $b_{0i}:V^{\sigma }\times V^{\sigma }\times V^{\gamma
}\rightarrow \mathbb{R}$ is uniformly bounded for some $\sigma \in \lbrack
-\theta _{2},\theta +\theta _{2}-2\varepsilon _{2})$, and $\gamma \in
\lbrack \theta +\varepsilon _{2},\infty )\cap (\varepsilon _{2},\infty )\cap
\left( \frac{n}{2},\infty \right) $;

\item[iv)] $A_{0i}:V^{\theta -\theta _{2}}\rightarrow V^{-\gamma }$ is
uniformly bounded and converges weakly to $A_{0}$;

\item[v)] $N_{i}^{-1}:V^{s+2\theta _{2}}\rightarrow V^{s+2\theta
_{2}-2\varepsilon _{2}}$ is uniformly bounded;

\item[vi)] $N_{i}^{-1}N:V^{\theta -\theta _{2}}\rightarrow V^{\theta +\theta
_{2}-2\varepsilon _{2}}$ converges strongly to the identity map;

\item[vii)] For any $v\in V^{\theta -\theta _{2}}$, $B_{0i}(v,v)$ converges
weakly to $B_{0}(v,v)$.
\end{itemize}

Then, there exists a solution%
\begin{equation*}
\left( u,\phi \right) \in L^{\infty }(0,T;\mathcal{Y}_{\theta _{2}})\cap
L^{2}(0,T;V^{\theta -\theta _{2}}\times D\left( A_{1}\right) )
\end{equation*}%
to (\ref{e:pert-lin-0}) such that, up to a subsequence, $%
y_{i}=N^{-1}N_{i}u_{i}$ and $\phi _{i}$ satisfy%
\begin{equation}
\left\{ 
\begin{array}{l}
\left( y_{i},\phi _{i}\right) \rightarrow \left( u,\phi \right) 
\mbox{
weak-star in }L^{\infty }(0,T;\mathcal{Y}_{\theta _{2}}), \\ 
\left( y_{i},\phi _{i}\right) \rightarrow \left( u,\phi \right) 
\mbox{
weakly in }L^{2}(0,T;V^{\theta -\theta _{2}}\times D\left( A_{1}\right) ),
\\ 
\left( y_{i},\phi _{i}\right) \rightarrow \left( u,\phi \right) 
\mbox{
strongly in }L^{2}(0,T;V^{s}\times W^{l}), \\ 
\phi _{i}\rightarrow \phi \mbox{
strongly in }C(0,T;W^{\zeta }),%
\end{array}%
\right.  \label{e:pertn-conv-u}
\end{equation}%
for any $s<\theta -\theta _{2},$ $l<2$ and $\zeta <1,$as $i\rightarrow
\infty $.
\end{theorem}

\begin{proof}
Recall that by Theorem \ref{t:exist}, we know that for $i\in \mathbb{N}$
there exists a solution to (\ref{e:pert-nonlin}) with the following
properties:%
\begin{align}
u_{i}& \in L^{\infty }(0,T;V^{-\varepsilon _{2}})\cap
L^{2}(0,T;V^{\varepsilon -\varepsilon _{2}}),  \label{e:pertn-0} \\
\phi _{i}& \in L^{\infty }\left( 0,T;L^{\infty }\left( \Omega \right)
\right) \text{ with }\phi _{i}\in \left[ -1,1\right] ,  \notag \\
\phi _{i}& \in L^{\infty }\left( 0,T;W^{1}\right) \cap L^{2}\left(
0,T;D\left( A_{1}\right) \right) .  \notag
\end{align}%
Pairing now the first and second equations of (\ref{e:pert-nonlin}) with $%
v_{i}:=N_{i}u_{i}$ and $\psi _{i}=A_{1}\phi _{i}+f\left( \phi _{i}\right) $,
respectively, after standard transformations, we have 
\begin{align}
& \frac{d}{dt}\left( \left\langle N_{i}^{-1}u_{i},u_{i}\right\rangle
+\left\Vert A_{1}^{1/2}\phi _{i}\right\Vert _{L^{2}}^{2}+2\int_{\Omega
}F\left( \phi _{i}\right) dx\right) +2\left\langle
A_{0i}N_{i}^{-1}u_{i},u_{i}\right\rangle  \label{e:pertn-1} \\
& =\left\langle g,u_{i}\right\rangle \lesssim \delta ^{-1}\Vert g\Vert
_{\theta -\theta _{2}}^{2}+\delta \Vert v_{i}\Vert _{\theta +\theta
_{2}}^{2}.  \notag
\end{align}%
Choosing $\delta \in \left( 0,1\right) $ small enough, then using (\ref%
{e:pert-nl-ell}), by Gr\"{o}nwall's inequality and (\ref{e:pert-nl-nell}) we
have 
\begin{equation}
\left( \Vert v_{i}(t)\Vert _{\theta _{2}}^{2}+\left\Vert \phi _{i}\left(
t\right) \right\Vert _{1}\right) \lesssim C_{T},\text{ }t\in \left(
0,T\right) .  \label{e:pertn-2}
\end{equation}%
Moreover, integrating (\ref{e:pertn-1}), and taking into account (\ref%
{e:pertn-0}), (\ref{e:pertn-2}) and (\ref{e:pert-nl-ell}), we infer that for
any fixed $T>0$,%
\begin{align*}
v_{i}& =N_{i}u_{i}\in L^{\infty }(0,T;V^{\theta _{2}})\cap
L^{2}(0,T;V^{\theta +\theta _{2}}), \\
\phi _{i}& \in L^{\infty }\left( 0,T;W^{1}\right) \cap L^{2}\left(
0,T;D\left( A_{1}\right) \right) ,
\end{align*}%
with uniformly (in $i$) bounded norms. On the other hand, we again have 
\begin{equation}
\Vert u_{i}^{^{\prime }}\Vert _{-\gamma }\leq \Vert A_{0i}u_{i}\Vert
_{-\gamma }+\Vert B_{0i}(u_{i},u_{i})\Vert _{-\gamma }+\Vert g+R_{0}\left(
A_{i}\phi _{i},\phi _{i}\right) \Vert _{-\gamma }
\end{equation}%
and%
\begin{equation*}
||\phi _{i}^{^{\prime }}||_{-2}\leq \left\Vert A_{1}\phi _{i}+f\left( \phi
_{i}\right) \right\Vert _{L^{2}}^{2}+\left\Vert B_{1}\left( u_{i},\phi
_{i}\right) \right\Vert _{-2}.
\end{equation*}%
By estimating the right hand sides as in the proof of Theorem \ref{t:exist},
we conclude that $\left( u_{i}^{^{\prime }},\phi _{i}^{^{\prime }}\right) $
is uniformly bounded in $L^{2}(0,T;V^{-\gamma }\times W^{-2})$; moreover, $%
v_{i}^{^{\prime }}=N_{i}u_{i}^{^{\prime }}$ is uniformly bounded in $%
L^{2}(0,T;V^{-\gamma })$. Compactness arguments as in the proof of Theorem %
\ref{t:exist}, allows us to pass to a subsequence, and, thus, infer the
existence of%
\begin{align*}
v& =Nu\in L^{\infty }(0,T;V^{\theta _{2}})\cap L^{2}(0,T;V^{\theta +\theta
_{2}}), \\
\phi & \in L^{\infty }\left( 0,T;W^{1}\cap L^{\infty }\left( \Omega \right)
\right) \cap L^{2}\left( 0,T;D\left( A_{1}\right) \right)
\end{align*}%
satisfying 
\begin{equation}
\begin{split}
& v_{i}\rightarrow v\mbox{ weak-star in }L^{\infty }(0,T;V^{\theta _{2}}), \\
& v_{i}\rightarrow v\mbox{ weakly in }L^{2}(0,T;V^{\theta +\theta _{2}}), \\
& v_{i}\rightarrow v\mbox{ strongly in }L^{2}(0,T;V^{s}),
\end{split}%
\end{equation}%
for any $s<\theta +\theta _{2}$, and%
\begin{equation*}
\phi _{i}\rightarrow \phi \mbox{ weakly in }L^{2}\left( 0,T;D\left(
A_{1}\right) \right) ,\text{ }\phi _{i}\rightarrow \phi \mbox{ weak-star in }%
L^{\infty }(0,T;W^{1}\cap L^{\infty }\left( \Omega \right) ),
\end{equation*}%
as $i\rightarrow \infty $. Define $u=N^{-1}v$ and $y_{i}=N^{-1}v_{i}$, and
note that these families satisfy (\ref{e:pertn-conv-u}).

We can now argue as in the proof of Theorem \ref{t:exist} to show that
indeed the limit $\left( u,\phi \right) $ satisfies the problem (\ref%
{e:pert-lin-0}). The procedure to pass to the limit is standard owing to the
following identities%
\begin{equation*}
\left\{ 
\begin{array}{l}
u_{i}-u=N_{i}^{-1}Ny_{i}-u=N_{i}^{-1}N(y_{i}-u)+(N_{i}^{-1}N-I)u, \\ 
A_{0i}u_{i}-A_{0}u=A_{0i}(u_{i}-u)+(A_{0i}-A_{0})u, \\ 
B_{0i}(u_{i},u_{i})-B_{0}(u,u)=B_{0i}(u_{i},u_{i}-u)+B_{0i}(u_{i}-u,u)+B_{0i}(u,u)-B_{0}(u,u),%
\end{array}%
\right.
\end{equation*}%
the uniform boundedness and weak convergence of $A_{0i}\rightarrow A_{0},$ $%
B_{0i}\rightarrow B_{0}$, and assumptions (iv)-(vii), see \cite[Theorem 4.2]%
{HLT}. The proof is finished.
\end{proof}

For example, setting $\varepsilon =\varepsilon _{2}=1$, with $\theta =1$ and 
$\theta _{2}=0$, and checking all the requirements (i)-(vii) of Theorem~\ref%
{t:pert-nonlin}, the weak solutions to the 3D NS--AC-$\alpha $ model
converge to a weak solution of the 3D NSE-AC model as the parameter $\alpha
\rightarrow 0$. This result was previously reported in~\cite{GM1}.

\section{Longtime behavior}

\label{s:longbehav}

In this section we establish the existence of global and exponential
attractors for the general three-parameter family of regularized models.
Moreover, assuming the potential $F$ to be real analytic and under
appropriate conditions on the external forces $g$ in (\ref{e:op}), we also
demonstrate that each trajectory converges to a single equilibrium, and find
a convergence rate estimate. We recall that by Theorem \ref{t:exist}, there
exists a weak solution%
\begin{equation*}
\left( u,\phi \right) \in L_{\mathrm{loc}}^{\infty }(0,\infty ;\mathcal{Y}%
_{\theta _{2}})\cap L_{\mathrm{loc}}^{2}(0,\infty ;V^{\theta -\theta
_{2}}\times D\left( A_{1}\right) )
\end{equation*}%
to (\ref{weak1})-(\ref{weak2}) with any given initial data $\left( u(0),\phi
\left( 0\right) \right) \in \mathcal{Y}_{\theta _{2}}$. By Theorem \ref%
{t:stab} the weak solution is unique and depends continuously on the initial
data in a Lipschitz way. Therefore, we have a continuous (nonlinear)\
semigroup%
\begin{equation}
\begin{array}{ll}
S_{\theta _{2}}(t):\mathcal{Y}_{\theta _{2}}\rightarrow \mathcal{Y}_{\theta
_{2}},\text{ }t\geq 0, &  \\ 
\left( u_{0},\phi _{0}\right) \mapsto \left( u\left( t\right) ,\phi \left(
t\right) \right) . & 
\end{array}
\label{semigroup}
\end{equation}%
Also for the sake of reference below, recall the following definition for
the space of translation bounded functions%
\begin{equation*}
L_{tb}^{2}\left( \mathbb{R}_{+};X\right) :=\left\{ g\in L_{loc}^{2}\left( 
\mathbb{R}_{+};X\right) :\left\Vert g\right\Vert _{L_{tb}^{2}\left( \mathbb{R%
}_{+};X\right) }^{2}:=\sup_{t\geq 0}\int_{t}^{t+1}\left\Vert g\left(
s\right) \right\Vert _{X}^{2}ds<\infty \right\} ,
\end{equation*}%
where $X$ is a given Banach space.

\subsection{Global and exponential attractors in the case $\protect\theta >0$%
}

\label{ss:attr-diss}

The following proposition establishes the existence of an absorbing ball in $%
\mathcal{Y}_{\theta _{2}}$ in both cases $\theta >0$ and $\theta =0$.
Moreover, with additional conditions in the case when $\theta >0$, we show
not only the existence of an absorbing ball in the space $V^{\beta }\times
D\left( A_{1}\right) $, but also that any weak solution with initial
condition in $\mathcal{Y}_{\theta _{2}}$ acquires additional smoothness in
an infinitesimal time. For our first result it suffices to have $\theta \geq
0.$

\begin{proposition}
\label{t:attr-exist}Let $\left( u,\phi \right) \in L_{\mathrm{loc}}^{\infty
}(0,\infty ;\mathcal{Y}_{\theta _{2}})\cap L_{\mathrm{loc}}^{2}(0,\infty
;V^{\theta -\theta _{2}}\times D\left( A_{1}\right) )$ be a weak solution in
the sense of Definition \ref{weak} with $\left( u(0),\phi \left( 0\right)
\right) \in \mathcal{Y}_{\theta _{2}}$. In addition, let the following
conditions hold.

(i) $\left\langle A_{0}v,Nv\right\rangle \geq c\Vert v\Vert _{\theta -\theta
_{2}}^{2}$ for any $v\in V^{\theta -\theta _{2}}$, with a constant $c>0$;

(ii) $g\in L_{tb}^{2}\left( \mathbb{R}_{+};V^{-\theta -\theta _{2}}\right) $;

Then for some constant $k>0$ independent of time and initial conditions, we
have 
\begin{align}
& \Vert u(t)\Vert _{-\theta _{2}}^{2}+\left\Vert \phi \left( t\right)
\right\Vert _{1}^{2}+\Vert \left( u,\phi \right) \Vert
_{L^{2}(t,t+1;V^{\theta -\theta _{2}}\times D\left( A_{1}\right) )}^{2}
\label{energy_fbis} \\
& \lesssim e^{-kt}\left( \Vert u(0)\Vert _{-\theta _{2}}^{2}+\left\Vert \phi
\left( 0\right) \right\Vert _{1}^{2}\right) +C_{\ast },  \notag
\end{align}%
for all $t\geq 0,$ for some $C_{\ast }>0$ independent of time and the
initial data.
\end{proposition}

\begin{proof}
Let us now set%
\begin{align}
E\left( t\right) :& =\Vert u(t)\Vert _{-\theta _{2}}^{2}+\left\Vert \phi
\left( t\right) \right\Vert _{1}^{2}+2\left\langle F\left( \phi \left(
t\right) \right) ,1\right\rangle _{L^{2}}+\left\Vert \phi \left( t\right)
\right\Vert _{L^{2}}^{2}+c_{E},  \label{energy} \\
\Theta \left( t\right) :& =-2\left\Vert u\left( t\right) \right\Vert
_{\theta -\theta _{2}}^{2}+\kappa \left\Vert u\left( t\right) \right\Vert
_{-\theta _{2}}^{2}-2\left\Vert \mu \left( t\right) \right\Vert _{L^{2}}^{2}
\notag \\
& -\left( 2-\kappa \right) \left\Vert A_{1}^{1/2}\phi \left( t\right)
\right\Vert _{L^{2}}^{2}  \notag \\
& +2\left[ \kappa \left\langle F\left( \phi \left( t\right) \right) -f\left(
\phi \left( t\right) \right) \phi \left( t\right) ,1\right\rangle
_{L^{2}}-\left( 1-\kappa \right) \left\langle f\left( \phi \left( t\right)
\right) \phi \left( t\right) ,1\right\rangle _{L^{2}}\right]  \notag \\
& +\kappa \left\Vert \phi \left( t\right) \right\Vert _{L^{2}}^{2}+2\kappa
C_{F}\text{vol}\left( \Omega \right) ,  \notag
\end{align}%
where $\kappa \in \left( 0,1\right) $. In \eqref{energy} we have set $%
c_{E}=2C_{F}$vol$\left( \Omega \right) >0$, with $C_{F}$ taken large enough
in order to ensure that $E\left( t\right) $ is nonnegative (note that $F$ is
bounded from below). On account of this choice and recalling Proposition \ref%
{maxp}, we can find $C>0$ such that 
\begin{equation}
\Vert u(t)\Vert _{-\theta _{2}}^{2}+\left\Vert \phi \left( t\right)
\right\Vert _{1}^{2}\leq E\left( t\right) \leq C\left( 1+\Vert u(t)\Vert
_{-\theta _{2}}^{2}+\left\Vert \phi \left( t\right) \right\Vert
_{1}^{2}\right) .  \label{3.9}
\end{equation}%
Following an argument from \cite[Proposition 3.1]{GG1}, on account of (i) we
have%
\begin{equation}
\frac{d}{dt}E\left( t\right) +\kappa E\left( t\right) \leq \Theta \left(
t\right) +\delta ^{-1}\Vert g\left( t\right) \Vert _{-\theta -\theta
_{2}}^{2}+\delta \Vert N\Vert _{-\theta _{2};\theta _{2}}^{2}\Vert u\left(
t\right) \Vert _{\theta -\theta _{2}}^{2},  \label{3.4}
\end{equation}%
for any $\delta >0$. Observe preliminarily that, owing to the first
assumption of (\ref{condf}), we have 
\begin{align}
& \left\vert f\left( y\right) \right\vert \left( 1+\left\vert y\right\vert
\right) \leq 2f\left( y\right) y+c_{f},  \label{3.6} \\
& F\left( y\right) -f\left( y\right) y\leq c_{f}^{\prime }\left\vert
y\right\vert ^{2}+c_{f}^{\prime \prime },  \label{3.7}
\end{align}%
for any $y\in \mathbb{R}$. Here $c_{f},$ $c_{f}^{\prime }$ and $%
c_{f}^{\prime \prime }$ are positive, sufficiently large constants that
depend on $f$ only. From (\ref{3.6})-(\ref{3.7}) and elementary
inequalities, it follows%
\begin{align*}
\Theta \left( t\right) & \leq -\left( 1-\kappa c_{\Omega }\left\vert \Omega
\right\vert \right) \left\Vert \boldsymbol{u}\left( t\right) \right\Vert
_{-\theta _{2}}^{2}-2\left\Vert \mu \left( t\right) \right\Vert
_{L^{2}}^{2}-\left( 2-\kappa \right) \left\Vert \nabla \phi \left( t\right)
\right\Vert _{L^{2}}^{2} \\
& -\left( 2-\kappa (1+2c_{f}^{\prime }\right) \left\Vert \phi \left(
t\right) \right\Vert _{L^{2}}^{2}-\left( 1-\kappa \right) \left\langle
\left\vert f\left( \phi \left( t\right) \right) \right\vert ,1+\left\vert
\phi \left( t\right) \right\vert \right\rangle _{L^{2}}+C,
\end{align*}%
where $c_{\Omega }$ depends on the shape of $\Omega ,$ but not on its size,
and $C>0$ depends on $\kappa $, $c_{f}$ and $c_{f}^{\prime \prime }$ at
most, but it is independent of time and the initial data. It is thus
possible to adjust sufficiently small $\kappa \in \left( 0,1\right) $ and $%
\delta >0,$ in order to have%
\begin{align}
& \frac{d}{dt}E\left( t\right) +\kappa E\left( t\right) +C\left( \left\Vert 
\boldsymbol{u}\left( t\right) \right\Vert _{\theta -\theta
_{2}}^{2}+\left\Vert \nabla \phi \left( t\right) \right\Vert
_{L^{2}}^{2}+\left\Vert \phi \left( t\right) \right\Vert _{L^{2}}^{2}\right)
+2\left\Vert \mu \left( t\right) \right\Vert _{L^{2}}^{2}  \label{diss} \\
& +C\left\langle \left\vert f\left( \phi \left( t\right) \right) \right\vert
,1+\left\vert \phi \left( t\right) \right\vert \right\rangle _{L^{2}}\leq
C\left( 1+\Vert g\left( t\right) \Vert _{-\theta -\theta _{2}}^{2}\right) . 
\notag
\end{align}%
Then, observing that assumption (\ref{condf}) also implies that $\left\vert
F\left( y\right) \right\vert \leq \left\vert f\left( y\right) \right\vert
\left( 1+\left\vert y\right\vert \right) +c_{f},$ for some positive constant 
$c_{f}$ and all $y\in \mathbb{R}$, and applying Gronwall's inequality (see
Appendix, Lemma \ref{Gineq}), we deduce%
\begin{align}
& E\left( t\right) +\int_{t}^{t+1}\left( C\left\Vert \boldsymbol{u}\left(
t\right) \right\Vert _{\theta -\theta _{2}}^{2}+2\left\Vert \mu \left(
s\right) \right\Vert _{L^{2}}^{2}+C\left\langle \left\vert F\left( y\left(
s\right) \right) \right\vert ,1\right\rangle _{L^{2}}\right) ds
\label{3.9bis} \\
& \leq E\left( 0\right) e^{-\kappa t}+C(1+\left\Vert g\right\Vert
_{L_{tb}^{2}\left( \mathbb{R}_{+};V^{-\theta -\theta _{2}}\right) }^{2}), 
\notag
\end{align}%
for all $t\geq 0$. Taking (\ref{3.9}) into account, by Proposition \ref{maxp}
we immediately obtain (\ref{energy_fbis}). This completes the proof.
\end{proof}

Consequently, for any $\theta \geq 0$ we have the following proposition.

\begin{proposition}
\label{t:attr-exist2}For every $R>0$, there exists $C_{\ast }=C_{\ast
}\left( R\right) >0$, independent of time, such that, for any $\varphi
_{0}:=\left( u_{0},\phi _{0}\right) \in \mathcal{B}_{\mathcal{Y}_{\theta
_{2}}}\left( R\right) ,$%
\begin{equation}
\sup_{t\geq 0}\left\Vert S_{\theta _{2}}\left( t\right) \varphi
_{0}\right\Vert _{\mathcal{Y}_{\theta _{2}}}+\int_{t}^{t+1}\left( \left\Vert
u\left( s\right) \right\Vert _{\theta -\theta _{2}}^{2}+\left\Vert A_{1}\phi
\left( s\right) \right\Vert _{L^{2}}^{2}\right) ds\leq C_{\ast },
\label{energy_ftris}
\end{equation}%
where $\mathcal{B}_{\mathcal{Y}_{\theta _{2}}}\left( R\right) $ denotes the
ball in $\mathcal{Y}_{\theta _{2}}$ of radius $R,$ centered at $0$.
\end{proposition}

The first main result of this subsection is the following.

\begin{theorem}
\label{t:attr-existb}Let the assumptions of Theorems \ref{t:exist} and \ref%
{t:stab} be satisfied for some $\theta >0.$ In addition, for some%
\begin{equation}
\beta \in (-\theta _{2},\min (\theta -\frac{1}{2},\theta -\theta _{2})],
\label{cond-beta}
\end{equation}%
when $n=3$, and%
\begin{equation}
\beta \in \lbrack \max \left( 1-2\theta _{2},-\theta _{2}\right) ,\min
(\theta ,\theta -\theta _{2})],\text{ }\beta \neq -\theta _{2},
\label{cond-beta-bis}
\end{equation}%
when $n=2,$ provided that the above intervals are nonempty, let the
following conditions hold.

(i) $b_{0}:V^{\alpha }\times V^{\alpha }\times V^{\theta -\beta }\rightarrow 
\mathbb{R}$ is bounded, where $\alpha =\min \{\beta ,\theta -\theta _{2}\}$;

(ii) $g\in V^{\beta -\theta }$ is time independent.

Then, there exists a compact attractor $\mathcal{A}\Subset \mathcal{Y}%
_{\theta _{2}}$ for the system (\ref{weak1})-(\ref{weak2}) which attracts
the bounded sets of $\mathcal{Y}_{\theta _{2}}$. Moreover, $\mathcal{A}$ is
connected and it is the maximal bounded attractor in $V^{\beta }\times
D\left( A_{1}\right) $.
\end{theorem}

\begin{proof}
By Propositions \ref{t:attr-exist} and \ref{maxp}, there is a ball $\mathcal{%
B}$ in $\mathcal{Y}_{\theta _{2}}$ which is absorbing in $\mathcal{Y}%
_{\theta _{2}}$, meaning that for any bounded set $U\subset \mathcal{Y}%
_{\theta _{2}}$ there exists $t_{0}=t_{0}(\left\Vert U\right\Vert _{\mathcal{%
Y}_{\theta _{2}}})>0$ such that $S_{\theta _{2}}(t)U\subset \mathcal{B}$ for
all $t\geq t_{0}$. Moreover, by Theorem \ref{t:reg} and application of the
uniform Gronwall's lemma \cite[Lemma III.1.1]{T}\ in (\ref{est15bis}) and (%
\ref{est16}), by virtue of (\ref{energy_fbis}) (cf. also (\ref{energy_ftris}%
)), we infer the existence of a new time $t_{1}=t_{0}+1$ such that%
\begin{equation}
\sup_{t\geq t_{1}}\left( \left\Vert u\left( t\right) \right\Vert _{\beta
}^{2}+\left\Vert A_{1}\left( \phi \left( t\right) \right) \right\Vert
_{L^{2}}^{2}\right) \leq C,  \label{comp-abs}
\end{equation}%
for some positive constant $C$ independent of time and the initial data.
Moreover, integration over $\left( t,t+1\right) $ of the previous
inequalities (\ref{est15bis}), (\ref{est16}) yields%
\begin{equation}
\sup_{t\geq t_{1}}\int_{t}^{t+1}\left( \left\Vert u\left( s\right)
\right\Vert _{\beta +\theta }^{2}+||A_{1}^{3/2}\left( \phi \left( s\right)
\right) ||_{L^{2}}^{2}\right) \leq C,  \label{comp-abs-bis}
\end{equation}%
owing once again to (\ref{comp-abs}). Thus, for any bounded set $U\subset 
\mathcal{Y}_{\theta _{2}},$ we have that $\cup _{t\geq t_{1}}S_{\theta
_{2}}(t)U$ is relatively compact in $\mathcal{Y}_{\theta _{2}}$, when $%
\mathcal{Y}_{\theta _{2}}$ is endowed with the metric topology of $%
V^{-\theta _{2}}\times W^{1}$. Finally, applying~\cite[Theorem I.1.1]{T} we
have that the set $\mathcal{A}=\cap _{s\geq 0}\overline{\cup _{t\geq
s}S_{\theta _{2}}(t)\mathcal{B}}$ is a compact attractor for $S_{\theta
_{2}} $, and the rest of the result is immediate.
\end{proof}

\begin{remark}
\label{attr-rem}

(i) All the special cases listed in Table \ref{t:spec} (except for the 3D
NSE-AC and 3D NSV-AC) satisfy the conditions of Theorem \ref{t:attr-existb}
when the space dimension is $n=3$, cf. also Table \ref{t:spec-reg}. For
instance, the uniform Gronwall's lemma cannot be applied to (\ref{est15bis}%
), (\ref{est16}), and so we cannot infer that (\ref{comp-abs}) is satisfied
in the case$\ \theta =0$. The NSV-AC system ($\theta =0,$ $\theta
_{1}=\theta _{2}=1$) will be handled in the next subsection. Note that an
absorbing set $\mathcal{B}$ in $\mathcal{Y}_{\theta _{2}}$ for problem (\ref%
{weak1})-(\ref{weak2}) can be constructed for all $\theta \geq 0$ on account
of Proposition \ref{t:attr-exist}. Indeed, there exists $R_{0}>0$
independent of time and initial data such that the ball $\mathcal{B}:=%
\mathcal{B}_{\mathcal{Y}_{\theta _{2}}}\left( R_{0}\right) $ is absorbing
for $S_{\theta _{2}}\left( t\right) $ on $\mathcal{Y}_{\theta _{2}}$.

(ii) We emphasize that estimate (\ref{comp-abs}) is also satisfied provided
that the external force $g$ is time dependent and $g\in L_{tb}^{2}\left( 
\mathbb{R}_{+};V^{\beta -\theta }\right) $. On account of this observation,
one can generalize the notion of global attractor and replace it by the
notion of pullback attractor, for example. One can still study the set of
all complete bounded trajectories, that is, trajectories which are bounded
for all $t\in \mathbb{R}.$ All the results that we have presented in this
section are still true in that case.

(iii) Note that (\ref{comp-abs}) also implies that the dynamical system $%
\left( S_{\theta _{2}}\left( t\right) ,\mathcal{Y}_{\theta _{2}}\right) $
possesses a compact absorbing set $\mathcal{B}_{\beta }$ which is contained
in $V^{\beta }\times D\left( A_{1}\right) $.
\end{remark}

Next we show the existence of exponential attractors for our regularized
family of models (\ref{e:op}) when $\theta >0$. It turns out that in order
to successfully construct an exponential attractor for problem (\ref{e:op}),
we need to derive a compact absorbing set with a higher degree of smoothness
than the set $\mathcal{B}_{\beta }\subset V^{\beta }\times D\left(
A_{1}\right) $. As in \cite{GG1}, this feature is due to the \emph{strong}
coupling of the regularized NSE with the Allen-Cahn equation. Indeed, as we
will see below it does not seem possible to get the exponential attractor
directly on the smooth set $\mathcal{B}_{\beta }$. The next lemma is
concerned with this issue.

\begin{lemma}
\label{l:smooth-compset}Let the assumptions of Theorem \ref{t:attr-existb}
be satisfied. In addition, for the same $\beta $ as in (\ref{cond-beta})-(%
\ref{cond-beta-bis}), and some value $\gamma \in \mathbb{R}$ satisfying%
\begin{equation*}
1-\theta -2\theta _{2}\leq \gamma \leq \min \left( \theta -\frac{n}{6}%
,\theta -1,\beta -\theta \right) ,
\end{equation*}%
provided that the above interval is nonempty, let the following conditions
hold.

(i) $\partial _{t}g\in L_{tb}^{2}\left( \mathbb{R}_{+};V^{\gamma -\theta
}\right) ,$ $g\in L^{\infty }\left( \mathbb{R}_{+};V^{\gamma }\right) ;$

(ii) $b_{0}:V^{\alpha }\times V^{\alpha }\times V^{\theta -\gamma
}\rightarrow \mathbb{R}$ is bounded, where $\alpha =\min \left( \gamma
,\beta +\theta \right) .$

Then, there exists a time $t_{2}\geq t_{1}$ such that%
\begin{equation}
\sup_{t\geq t_{2}}\left( \left\Vert u\left( t\right) \right\Vert _{2\theta
+\gamma }^{2}+||A_{1}^{3/2}\left( \phi \left( t\right) \right)
||_{L^{2}}^{2}\right) \leq C,  \label{comp-abs2}
\end{equation}%
for some positive constant $C$ independent of time and the initial data.
\end{lemma}

\begin{proof}
The following arguments are formal for the sake of simplicity, but they can
be rigorously justified within the Galerkin scheme used in the proof of
Theorem \ref{t:exist}. We first observe that, using the apriori bounds (\ref%
{comp-abs}), (\ref{comp-abs-bis}), and arguing exactly as in the proof of
Theorem \ref{t:reg}, (\ref{basic00})-(\ref{est15q}), from the second
equation of (\ref{e:op}) by comparison, we have%
\begin{equation}
\sup_{t\geq t_{1}}\left( \left\Vert \partial _{t}\phi \left( t\right)
\right\Vert _{L^{2}}^{2}+\int_{t}^{t+1}\left\Vert \partial _{t}\phi \left(
s\right) \right\Vert _{1}^{2}ds\right) \leq C.  \label{est-h1}
\end{equation}%
Similarly, by comparison in the first equation of (\ref{e:op}), we also have%
\begin{equation}
\sup_{t\geq t_{1}}\left( \left\Vert \partial _{t}u\left( t\right)
\right\Vert _{\beta -2\theta }^{2}+\int_{t}^{t+1}\left\Vert \partial
_{t}u\left( s\right) \right\Vert _{\beta -\theta }^{2}ds\right) \leq C,
\label{est-h2}
\end{equation}%
owing to the assumptions on $A_{0},$ $g\in L_{tb}^{2}\left( \mathbb{R}%
_{+};V^{\beta -\theta }\right) $ and boundednes of the form $b_{0},$ see
assumptions (i)-(ii) of Theorem \ref{t:attr-existb}.

Set $v:=\partial _{t}u,$ $\psi :=\partial _{t}\phi $, and observe that $%
\left( v,\psi \right) $ solves the following system%
\begin{equation}
\left\{ 
\begin{array}{l}
\partial _{t}v+A_{0}v+B_{0}(v,u)+B_{0}\left( u,v\right) =R_{0}\left(
A_{1}\psi ,\phi \right) +R_{0}\left( A_{1}\phi ,\psi \right) +\partial _{t}g,
\\ 
\partial _{t}\psi +B_{1}\left( v,\phi \right) +B_{1}\left( u,\psi \right)
+A_{1}\psi =-f^{^{\prime }}\left( \phi \right) \psi .%
\end{array}%
\right.  \label{e:time-deriv-eq}
\end{equation}%
We repeat the arguments from the proof of Theorem \ref{t:reg}, and whenever
necessary show the new estimates. Pairing the first and second equations of (%
\ref{e:time-deriv-eq}) with $\Lambda ^{2\gamma }v$ and $A_{1}\psi ,$
respectively, then adding the identities together, we deduce%
\begin{align}
& \frac{1}{2}\frac{d}{dt}\left( \left\Vert v\right\Vert _{\gamma
}^{2}+\left\Vert A_{1}^{1/2}\psi \right\Vert _{L^{2}}^{2}\right)
+\left\langle A_{0}v,\Lambda ^{2\gamma }v\right\rangle +\left\Vert A_{1}\psi
\right\Vert _{L^{2}}^{2}  \label{est-h3} \\
& =\left\langle \partial _{t}g,\Lambda ^{2\gamma }v\right\rangle
-b_{0}\left( v,u,\Lambda ^{2\gamma }v\right) -b_{0}\left( u,v,\Lambda
^{2\gamma }v\right)  \notag \\
& +\left\langle R_{0}\left( A_{1}\psi ,\phi \right) ,\Lambda ^{2\gamma
}v\right\rangle +\left\langle R_{0}\left( A_{1}\phi ,\psi \right) ,\Lambda
^{2\gamma }v\right\rangle -\left\langle f^{^{\prime }}\left( \phi \right)
\psi ,A_{1}\psi \right\rangle  \notag \\
& -b_{1}\left( v,\phi ,A_{1}\psi \right) -b_{1}\left( u,\psi ,A_{1}\psi
\right) .  \notag
\end{align}%
Some of the estimates on the right-hand side of (\ref{est-h3}) are easy. We
have%
\begin{equation}
\left\langle \partial _{t}g,\Lambda ^{2\gamma }v\right\rangle \lesssim
\delta ^{-1}\left\Vert \partial _{t}g\right\Vert _{\gamma -\theta
}^{2}+\delta \left\Vert v\right\Vert _{\gamma +\theta }^{2},\text{ for any }%
\delta >0  \label{est-h4}
\end{equation}%
and%
\begin{equation}
\left\langle f^{^{\prime }}\left( \phi \right) \psi ,A_{1}\psi \right\rangle
\leq \delta \left\Vert A_{1}\psi \right\Vert _{L^{2}}^{2}+C\delta
^{-1}\left\Vert \psi \right\Vert _{1}^{2}.  \label{est-h5}
\end{equation}%
To bound the last term we argue verbatim as in the proof of Theorem \ref%
{t:stab}, (\ref{uniqest4})-(\ref{uniqest7}). We have%
\begin{equation}
b_{1}\left( u,\psi ,A_{1}\psi \right) \leq \delta \left\Vert A_{1}\psi
\right\Vert _{L^{2}}^{2}+C\delta ^{-1/\left( p-1\right) }\left\Vert
u\right\Vert _{\theta -\theta _{2}}^{2}\left\Vert u\right\Vert _{-\theta
_{2}}^{\left( p-2\right) }\left\Vert \psi \right\Vert _{1}^{2},
\label{est-h5bis}
\end{equation}%
in two space dimensions, and%
\begin{equation}
b_{1}\left( u,\psi ,A_{1}\psi \right) \leq C\delta ^{-7}\left\Vert
u\right\Vert _{-\theta _{2}}^{6}\left\Vert u\right\Vert _{\theta -\theta
_{2}}^{2}\left\Vert \psi \right\Vert _{1}^{2}+\delta \left\Vert A_{1}\psi
\right\Vert _{L^{2}}^{2},  \label{est-h5tris}
\end{equation}%
when $n=3$. In order to bound the last remaining terms in (\ref{est-h3}), we
argue as follows. Using (\ref{basic00}), the fact that $V^{\frac{n}{6}%
}\subset L^{3}$ and $\gamma -\theta \leq -n/6$, we have%
\begin{align}
\left\langle R_{0}\left( A_{1}\phi ,\psi \right) ,\Lambda ^{2\gamma
}v\right\rangle & \lesssim \left\Vert \Lambda ^{2\gamma }v\right\Vert
_{\theta -\gamma }\left\Vert R_{0}\left( A_{1}\phi ,\psi \right) \right\Vert
_{\gamma -\theta }  \label{est-h8} \\
& \lesssim \left\Vert v\right\Vert _{\theta +\gamma }\left\Vert A_{1}\phi
\right\Vert _{L^{6}}\left\Vert \nabla \psi \right\Vert _{L^{2}}  \notag \\
& \lesssim \delta \left\Vert v\right\Vert _{\theta +\gamma }^{2}+\delta
^{-1}\left\Vert A_{1}^{3/2}\phi \right\Vert _{L^{2}}^{2}\left\Vert \psi
\right\Vert _{1}^{2}.  \notag
\end{align}%
Next, from standard interpolation inequalities and (\ref{basic00}), we infer%
\begin{align}
\left\langle R_{0}\left( A_{1}\psi ,\phi \right) ,\Lambda ^{2\gamma
}v\right\rangle & \lesssim \left\Vert \Lambda ^{2\gamma }v\right\Vert
_{\theta -\gamma }\left\Vert R_{0}\left( A_{1}\psi ,\phi \right) \right\Vert
_{\gamma -\theta }  \label{est-h9} \\
& \lesssim \left\Vert v\right\Vert _{\theta +\gamma }\left\Vert \psi
\right\Vert _{1}\left\Vert A_{1}^{3/2}\phi \right\Vert _{L^{2}}  \notag \\
& \leq \delta \left\Vert v\right\Vert _{\theta +\gamma }^{2}+C\delta
^{-1}\left\Vert \psi \right\Vert _{1}^{2}\left\Vert A_{1}^{3/2}\phi
\right\Vert _{L^{2}}^{2},  \notag
\end{align}%
thanks to the inequality%
\begin{equation}
\left\Vert R_{0}\left( A_{1}\psi ,\phi \right) \right\Vert _{\gamma -\theta
}\leq \left\Vert R_{0}\left( A_{1}\psi ,\phi \right) \right\Vert
_{-1}\lesssim \left\Vert \psi \right\Vert _{1}\left\Vert A_{1}^{3/2}\phi
\right\Vert _{L^{2}}.  \label{est-h10}
\end{equation}%
Finally, we can argue as above in (\ref{est-h8})-(\ref{est-h10}) to obtain,
owing to the boundedness of the map $N:V^{s}\rightarrow V^{s+2\theta _{2}},$%
\begin{align}
b_{1}\left( v,\phi ,A_{1}\psi \right) & =\left\langle R_{0}\left( A_{1}\psi
,\phi \right) ,Nv\right\rangle  \label{est-h11} \\
& \lesssim \left\Vert Nv\right\Vert _{\gamma +\theta +2\theta
_{2}}\left\Vert R_{0}\left( A_{1}\psi ,\phi \right) \right\Vert _{-\left(
\gamma +\theta +2\theta _{2}\right) }  \notag \\
& \lesssim \left\Vert v\right\Vert _{\gamma +\theta }\left\Vert R_{0}\left(
A_{1}\psi ,\phi \right) \right\Vert _{-1}  \notag \\
& \leq \delta \left\Vert v\right\Vert _{\gamma +\theta }^{2}+C\delta
^{-1}\left\Vert \psi \right\Vert _{1}^{2}\left\Vert A_{1}^{3/2}\phi
\right\Vert _{L^{2}}^{2}.  \notag
\end{align}%
The assumption (ii) on $b_{0}$ is enough to bound the following terms:%
\begin{align}
& b_{0}\left( v,u,\Lambda ^{2\gamma }v\right) +b_{0}\left( u,v,\Lambda
^{2\gamma }v\right)  \label{est-h12} \\
& \lesssim \left\Vert v\right\Vert _{\gamma }\left\Vert u\right\Vert _{\beta
+\theta }\left\Vert \Lambda ^{2\gamma }v\right\Vert _{\theta -\gamma } 
\notag \\
& \leq \delta \left\Vert v\right\Vert _{\gamma +\theta }^{2}+C\delta
^{-1}\left\Vert v\right\Vert _{\gamma }^{2}\left\Vert u\right\Vert _{\beta
+\theta }^{2}.  \notag
\end{align}%
Setting%
\begin{equation*}
Z\left( t\right) :=\left\Vert v\left( t\right) \right\Vert _{\gamma
}^{2}+\left\Vert \psi \left( t\right) \right\Vert _{1}^{2},
\end{equation*}%
and inserting all the previous estimates (\ref{est-h4})-(\ref{est-h11}) into
(\ref{est-h3}), then using the coercitivity assumption on $A_{0}$ and
choosing a sufficiently small $\delta \sim \min \left(
c_{A_{0}},c_{A_{1}}\right) >0,$ we obtain%
\begin{equation}
\frac{d}{dt}Z\left( t\right) +\left\Vert v\left( t\right) \right\Vert
_{\gamma +\theta }^{2}+\left\Vert A_{1}\psi \left( t\right) \right\Vert
_{L^{2}}^{2}\leq \Delta \left( t\right) Z\left( t\right) +C\left\Vert
\partial _{t}g\right\Vert _{\gamma -\theta }^{2},  \label{est-h13}
\end{equation}%
where the function $\Delta \in L^{1}\left( t,t+1\right) $ is%
\begin{equation*}
\Delta :=C_{\delta }\left( 1+\left\Vert A_{1}^{3/2}\phi \right\Vert
_{L^{2}}^{2}+\left\Vert u\right\Vert _{-\theta _{2}}^{6}\left\Vert
u\right\Vert _{\theta -\theta _{2}}^{2}+\left\Vert u\right\Vert _{\beta
+\theta }^{2}\right) ,
\end{equation*}%
for some $C_{\delta }>0$. Observe now that, on account of (\ref{comp-abs})-(%
\ref{comp-abs-bis}) and (\ref{est-h1})-(\ref{est-h2}), we have%
\begin{equation*}
\sup_{t\geq t_{1}}\int_{t}^{t+1}\Delta \left( s\right) ds\leq C,\text{ }%
\sup_{t\geq t_{1}}\int_{t}^{t+1}Z\left( s\right) ds\leq C.
\end{equation*}%
Hence, we can exploit the uniform Gronwall's lemma \cite[Lemma III.1.1]{T}
in (\ref{est-h13}) to infer the existence of a new time $t_{2}\geq t_{1}$
such that%
\begin{equation}
\sup_{t\geq t_{2}}\left\Vert \partial _{t}u\left( t\right) \right\Vert
_{\gamma }+\left\Vert \partial _{t}\phi \left( t\right) \right\Vert _{1}\leq
C,  \label{est-h14}
\end{equation}%
for some positive constant $C$ independent of time and the initial data.
Comparison in the second equation of (\ref{e:op}) yields the desired bound
for $\phi $ in (\ref{comp-abs2}) owing once more to (\ref{comp-abs}) and (%
\ref{est-h14}). Finally, from the first equation of (\ref{e:op}), we have%
\begin{equation}
\left\Vert u\right\Vert _{2\theta +\gamma }\approx \left\Vert
A_{0}u\right\Vert _{\gamma }\lesssim \left\Vert g\right\Vert _{\gamma
}+\left\Vert B_{0}\left( u,u\right) \right\Vert _{\beta -\theta }+\left\Vert
R_{0}\left( A_{1}\phi ,\phi \right) \right\Vert _{\beta -\theta },
\label{est-h15}
\end{equation}%
since $\gamma \leq \beta -\theta .$ Owing to the assumptions on $\beta ,$
and on $b_{0}$ in the statement of Theorem \ref{t:attr-existb}, each term on
the right-hand side is essentially bounded for times greater than $t_{2}.$
Estimate (\ref{est-h15}) entails the desired estimate for $u$ from (\ref%
{comp-abs2}), and the proof is finished.
\end{proof}

The next lemma is concerned with the Lipschitz-in-time regularity of the
semigroup $S_{\theta _{2}}(t).$

\begin{lemma}
\label{time_reg-positive}Let the assumptions of Lemma \ref{l:smooth-compset}
be satisfied. For any $R>0,$ there exists a time $t_{\ast }=t_{\ast }\left(
R\right) >0$, such that%
\begin{equation}
\left\Vert S_{\theta _{2}}(t)\varphi _{0}-S_{\theta _{2}}(\tilde{t})\varphi
_{0}\right\Vert _{V^{\gamma }\times W^{1}}\leq C|t-\tilde{t}|,
\label{time_regularity-posit}
\end{equation}%
for all $t,\tilde{t}\in \lbrack t_{\ast },\infty )$ and any $\varphi
_{0}=\left( u_{0},\varphi _{0}\right) \in \mathcal{B}_{\theta ,\gamma
}\left( R\right) \subset V^{2\theta +\gamma }\times D(A_{1}^{3/2}).$ Here $%
\mathcal{B}_{\theta ,\gamma }\left( R\right) $ denotes any ball in $%
V^{2\theta +\gamma }\times D(A_{1}^{3/2})$ of radius $R>0$, centered at $0.$
\end{lemma}

\begin{proof}
The claim (\ref{time_regularity-posit}) follows from the basic equality%
\begin{equation*}
S_{\theta _{2}}(t)\varphi _{0}-S_{\theta _{2}}(\tilde{t})\varphi
_{0}=\int_{t}^{\widetilde{t}}\partial _{y}\left( S_{\theta _{2}}(y)\varphi
_{0}\right) dy
\end{equation*}%
and estimate (\ref{est-h14}).
\end{proof}

With the essential Lemma \ref{l:smooth-compset}, the next result states the
validity of the smoothing property for the semigroup $S_{\theta _{2}}\left(
t\right) $ in the case $\theta >0$.

\begin{lemma}
\label{smooth-property-diss}Let the assumptions of Lemma \ref%
{l:smooth-compset} be satisfied. Indicate by $\left( u_{i},\phi _{i}\right) $
the solution to problem (\ref{e:op}) which corresponds to the initial data $%
\left( u_{i}\left( 0\right) ,\phi _{i}\left( 0\right) \right) \in \mathcal{B}
$, where $i=1,2$. Then the following estimate holds:%
\begin{align}
\left\Vert u_{1}\left( t\right) -u_{2}\left( t\right) \right\Vert & _{\beta
}^{2}+\left\Vert A_{1}\left( \phi _{1}\left( t\right) -\phi _{2}\left(
t\right) \right) \right\Vert _{L^{2}}^{2}  \label{3.59} \\
& \leq C\frac{\overline{t}+1}{\overline{t}}e^{Ct}\left( \left\Vert
u_{1}\left( 0\right) -u_{2}\left( 0\right) \right\Vert _{-\theta
_{2}}^{2}+\left\Vert \phi _{1}\left( 0\right) -\phi _{2}\left( 0\right)
\right\Vert _{1}^{2}\right) ,  \notag
\end{align}%
for all $\overline{t}:=t-t_{2}>0$, and some positive constant $C$ which only
depends on $\mathcal{B}$.
\end{lemma}

\begin{proof}
First, recall that by (\ref{comp-abs2}), there is a compact absorbing set
for $S_{\theta _{2}}\left( t\right) $ contained in $V^{2\theta +\gamma
}\times D(A_{1}^{3/2})$.

As before, let $v=u_{1}-u_{2}$ and $\psi =\phi _{1}-\phi _{2}$ and recall
that $\left( v,\psi \right) $ solves equations (\ref{diffuniq1})-(\ref%
{diffuniq2}). Taking $w=\Lambda ^{2\beta }v$ and $\eta =A_{1}^{2}\psi $ into
(\ref{diffuniq1}) and (\ref{diffuniq2}), respectively, we infer%
\begin{align}
& \frac{1}{2}\frac{d}{dt}\left( \Vert v\Vert _{\beta }^{2}+\left\Vert
A_{1}\psi \right\Vert _{L^{2}}^{2}\right) +c_{A_{0}}\Vert v\Vert _{\theta
+\beta }^{2}+\left\Vert A_{1}^{3/2}\psi \right\Vert _{L^{2}}^{2}
\label{diffhigher0} \\
& \leq b_{0}(v,u_{1},\Lambda ^{2\beta }v)+b_{0}(u_{2},v,\Lambda ^{2\beta
}v)+\left\langle R_{0}\left( A_{1}\phi _{2},\psi \right) ,\Lambda ^{2\beta
}v\right\rangle  \notag \\
& +\left\langle R_{0}\left( A_{1}\psi ,\phi _{1}\right) ,\Lambda ^{2\beta
}v\right\rangle +\left\langle f\left( \phi _{1}\right) -f\left( \phi
_{2}\right) ,A_{1}^{2}\psi \right\rangle +b_{1}\left( v,\phi
_{1},A_{1}^{2}\psi \right)  \notag \\
& +b_{1}\left( u_{2},\psi ,A_{1}^{2}\psi \right)  \notag \\
& =:I_{1}+...+I_{7}.  \notag
\end{align}%
The first two terms are bounded exactly as in (\ref{est10}). We have%
\begin{equation}
I_{1}\leq C\delta ^{-1}\Vert u_{1}\Vert _{\theta -\theta _{2}}^{2}\Vert
v\Vert _{\beta }^{2}+\delta \left\Vert v\right\Vert _{\beta +\theta }^{2},%
\text{ }I_{2}\leq C\delta ^{-1}\Vert u_{2}\Vert _{\theta -\theta
_{2}}^{2}\Vert v\Vert _{\beta }^{2}+\delta \left\Vert v\right\Vert _{\beta
+\theta }^{2},  \label{diffuniq3}
\end{equation}%
for any $\delta >0$. Concerning the terms $I_{3},I_{4},$ we can argue in a
similar fashion as in the estimates (\ref{est13}), (\ref{est13bis}) to deduce%
\begin{equation}
I_{3}\leq \delta \left\Vert v\right\Vert _{\beta +\theta }^{2}+C\delta
^{-1}\left( \left\Vert A_{1}\phi _{2}\right\Vert _{L^{2}}^{2}\right)
\left\Vert A_{1}\psi \right\Vert _{L^{2}}^{2}
\end{equation}%
when $n=3$, and in two space dimensions,%
\begin{equation}
I_{3}\leq \delta \left\Vert v\right\Vert _{\beta +\theta }^{2}+\delta
\left\Vert A_{1}^{3/2}\psi \right\Vert _{L^{2}}^{2}+C\delta ^{-2}\left\Vert
A_{1}\phi _{2}\right\Vert _{L^{2}}^{4}\left\Vert \nabla \psi \right\Vert
_{L^{2}}^{2}.
\end{equation}%
Analogously, by elementary Sobolev inequalities, we have%
\begin{align}
I_{4}& \leq \delta \left\Vert v\right\Vert _{\beta +\theta }^{2}+C\delta
^{-1}\left( \left\Vert A_{1}\phi _{1}\right\Vert _{L^{2}}^{2}\right)
\left\Vert A_{1}\psi \right\Vert _{L^{2}}^{2},\text{ } \\
I_{4}& \leq \delta \left\Vert v\right\Vert _{\beta +\theta }^{2}+\delta
\left\Vert A_{1}^{3/2}\psi \right\Vert _{L^{2}}^{2}+C\delta ^{-2}\left\Vert
A_{1}\psi \right\Vert _{L^{2}}^{2}\left( \left\Vert A_{1}\phi
_{1}\right\Vert _{L^{2}}^{4}\right) ,  \notag
\end{align}%
in both two and three space dimensions, respectively. Moreover, it easy to
see that%
\begin{align}
I_{5}& =\left\langle A_{1}^{1/2}f\left( \phi _{1}\right) -f\left( \phi
_{2}\right) ,A_{1}^{3/2}\psi \right\rangle \\
& \leq \delta \left\Vert A_{1}^{3/2}\psi \right\Vert _{L^{2}}^{2}+C\delta
^{-1}\left( 1+\left\Vert A_{1}\phi _{1}\right\Vert _{L^{2}}^{2}+\left\Vert
A_{1}\phi _{1}\right\Vert _{L^{2}}^{2}\right) \left\Vert \psi \right\Vert
_{1}^{2}.  \notag
\end{align}%
In order to bound the last two terms on the right-hand side of (\ref%
{diffhigher0}), we repeat the same estimates derived in the proof of Theorem %
\ref{t:reg}, see (\ref{basic00})-(\ref{est15}) and (\ref{est14bis})-(\ref%
{est15q}). First, \ from (\ref{est14})-(\ref{est15}) when $n=3$,%
\begin{align}
I_{7}& =\left\langle A_{1}^{1/2}B_{1}\left( u_{2},\psi \right)
,A_{1}^{3/2}\psi \right\rangle \\
& \leq \delta \left\Vert A_{1}^{3/2}\psi \right\Vert _{L^{2}}^{2}+C\delta
^{-7}\left\Vert u_{2}\right\Vert _{-\theta _{2}}^{8}\left\Vert \psi
\right\Vert _{1}^{2}+C\delta ^{-3}\left\Vert u_{2}\right\Vert _{-\theta
_{2}}^{4}\left\Vert A_{1}\psi \right\Vert _{L^{2}}^{2}  \notag
\end{align}%
and, from (\ref{est14bis})-(\ref{est15q}) when $n=2$, we have%
\begin{equation*}
I_{7}\leq \delta \left\Vert A_{1}^{3/2}\psi \right\Vert _{L^{2}}^{2}+C\delta
^{-3}\left( \left\Vert u_{2}\right\Vert _{-\theta _{2}}^{2}\left\Vert
u_{2}\right\Vert _{\theta -\theta _{2}}^{2}\left\Vert A_{1}\psi \right\Vert
_{L^{2}}^{2}+\left\Vert u_{2}\right\Vert _{\theta -\theta
_{2}}^{2}\left\Vert u_{2}\right\Vert _{\beta }^{2}\left\Vert \nabla \psi
\right\Vert _{L^{2}}^{2}\right) .
\end{equation*}

Finally, the bound for the term $I_{6}$ is similar. As in the proof of
Theorem \ref{t:reg}, we again have%
\begin{align}
I_{6}& =\left\langle A_{1}^{1/2}B_{1}\left( v,\phi _{1}\right)
,A_{1}^{3/2}\psi \right\rangle \lesssim \left\Vert B_{1}\left( v,\phi
_{1}\right) \right\Vert _{1}\left\Vert A_{1}^{3/2}\psi \right\Vert _{L^{2}}
\\
& \leq \delta \left\Vert A_{1}^{3/2}\psi \right\Vert _{L^{2}}^{2}+C\delta
^{-1}\left\Vert v\right\Vert _{-\theta _{2}}^{2}\left\Vert A_{1}\phi
_{1}\right\Vert _{L^{2}}\left\Vert A_{1}^{3/2}\phi _{1}\right\Vert _{L^{2}},
\notag
\end{align}%
when $n=3,$ and%
\begin{equation}
I_{6}\leq \delta \left\Vert A_{1}^{3/2}\psi \right\Vert _{L^{2}}^{2}+C\delta
^{-1}\left\Vert v\right\Vert _{\beta }^{2}\left\Vert A_{1}\phi
_{1}\right\Vert _{L^{2}}\left\Vert A_{1}^{3/2}\phi _{1}\right\Vert _{L^{2}},
\label{diffuniq4}
\end{equation}%
in two space dimensions ($n=2$). Let us now set%
\begin{equation*}
X\left( t\right) :=\left\Vert v\left( t\right) \right\Vert _{\beta
}^{2}+\left\Vert A_{1}\psi \left( t\right) \right\Vert _{L^{2}}^{2}.
\end{equation*}%
Collecting all the above estimates from (\ref{diffuniq3}) to (\ref{diffuniq4}%
), and inserting them into the right-hand side of (\ref{diffhigher0}), and
choosing $\delta \sim \min \left( c_{A_{0}},c_{A_{1}}\right) >0$
sufficiently small, we deduce the following inequality%
\begin{equation}
\frac{d}{dt}X\left( t\right) \leq \Xi \left( t\right) X\left( t\right) ,%
\text{ for all }t\geq t_{2},  \label{diffuniq5}
\end{equation}%
with the obvious definition for $\Xi $. We emphasize that in both dimensions 
$n=2,3,$ the function $\Xi \in L^{\infty }\left( t_{2},\infty \right) $,
thanks now to the uniform estimates (\ref{comp-abs}), (\ref{comp-abs-bis})
and (\ref{comp-abs2}). Multiplying now both sides of this inequality by $%
\overline{t}:=t-t_{2}$ and integrating the resulting relation over $\left(
t_{2},t\right) ,$ we get 
\begin{equation*}
\overline{t}X\left( t\right) \leq C\int_{t_{2}}^{t}\left( s-t_{2}+1\right)
X\left( s\right) ds,\quad \text{for all }t>t_{2},
\end{equation*}%
which entails (\ref{3.59}), owing to Theorem \ref{t:stab}, (\ref{uniq_stab}%
). The proof is complete.
\end{proof}

The second main result of this subsection is concerned with the existence of
exponential attractors for problem (\ref{e:op}) in the case $\theta >0$.

\begin{theorem}
\label{expo-thm}Let the assumptions of Lemma \ref{l:smooth-compset} be
satisfied, and assume $g$ is time independent. Then $\left( S_{\theta _{2}},%
\mathcal{Y}_{\theta _{2}}\right) $ possesses an exponential attractor $%
\mathcal{M}_{\theta _{2}}\subset \mathcal{Y}_{\theta _{2}}$ which is bounded
in $V^{2\theta +\gamma }\times D(A_{1}^{3/2})$. Thus, by definition, we have

(a) $\mathcal{M}_{\theta _{2}}$ is compact and semi-invariant with respect $%
S_{\theta _{2}}\left( t\right) ,$ that is,%
\begin{equation*}
S_{\theta _{2}}\left( t\right) \left( \mathcal{M}_{\theta _{2}}\right)
\subseteq \mathcal{M}_{\theta _{2}},\quad \forall \,t\geq 0.
\end{equation*}

(b) The fractal dimension $\dim _{F}\left( \mathcal{M}_{\theta _{2}},%
\mathcal{Y}_{\theta _{2}}\right) $ of $\mathcal{M}_{\theta _{2}}$ is finite
and an upper bound can be computed explicitly.

(c) $\mathcal{M}_{\theta _{2}}$ attracts exponentially fast any bounded
subset $B$ of $\mathcal{Y}_{\theta _{2}}$, that is, there exist a positive
nondecreasing function $Q$ and a constant $\rho >0$ such that 
\begin{equation*}
dist_{\mathcal{Y}_{\theta _{2}}}\left( S_{\theta _{2}}\left( t\right) B,%
\mathcal{M}_{\theta _{2}}\right) \leq Q(\Vert B\Vert _{\mathcal{Y}_{\theta
_{2}}})e^{-\rho t},\quad \forall \,t\geq 0.
\end{equation*}%
Here $dist_{\mathcal{Y}_{\theta _{2}}}$ denotes the Hausdorff semi-distance
between sets in $\mathcal{Y}_{\theta _{2}}$ and $\Vert B\Vert _{\mathcal{Y}%
_{\theta _{2}}}$ stands for the size of $B$ in $\mathcal{Y}_{\theta _{2}}.$
Both $Q$ and $\rho $ can be explicitly calculated.
\end{theorem}

\begin{proof}
Using Theorem \ref{t:stab}, Theorem \ref{t:attr-existb} and Lemma \ref%
{l:smooth-compset}, we can find a bounded subset $X_{0}$ of $V^{2\theta
+\gamma }\times D(A_{1}^{3/2})$ and a time $t^{\sharp }>0$ such that,
setting $\Sigma =S_{\theta _{2}}(t^{\sharp })$, the mapping $\Sigma
:X_{0}\rightarrow X_{0}$ enjoys the smoothing property (\ref{3.59}).
Therefore Theorem \ref{t4.5} applies to $\Sigma $ and there exists a compact
set $\mathcal{M}_{\theta _{2}}^{\ast }\in X_{0}$ with finite fractal
dimension (with respect to the metric topology of $V^{\gamma }\times W^{1}$)
that satisfies (\ref{3.57}) and (\ref{3.58}). Hence, setting 
\begin{equation*}
\mathcal{M}_{\theta _{2}}\mathcal{=}\cup _{t\in \left[ t^{\sharp
},2t^{\sharp }\right] }S_{\theta _{2}}\left( t\right) \mathcal{M}_{\theta
_{2}}^{\ast },
\end{equation*}%
we deduce that (a) and (c) are fulfilled, while (b) is a consequence of
Theorem \ref{t:stab} and (\ref{time_regularity-posit}). The attraction
property (c) in the metric topology of $V^{-\theta _{2}}\times W^{1}$ is a
standard corollary of the aforementioned properties, standard interpolation
inequalities and the fact that $\mathcal{M}_{\theta _{2}}$ is bounded in $%
V^{2\theta +\gamma }\times D(A_{1}^{3/2})$.
\end{proof}

As a consequence of the above theorem, we have the following.

\begin{corollary}
Under the assumptions of Theorem \ref{expo-thm}, the global attractor $%
\mathcal{A}$ is bounded in $V^{2\theta +\gamma }\times D(A_{1}^{3/2})$, and $%
\mathcal{A}$ has finite fractal dimension.
\end{corollary}

Note that Theorem \ref{t:attr-existb} provides many examples where the
conditions of the Lemma \ref{l:smooth-compset} are satisfied (with $\beta
=\theta -\theta _{2}$ and $\gamma =\beta -\theta =-\theta _{2}$). For
example, setting $\theta _{1}=\theta _{2}=0$, with $\theta =1$, and checking
all the requirements (i)--(ii) of Theorem \ref{t:attr-existb}, and
conditions (i)-(ii) of Lemma \ref{l:smooth-compset}, the 2D NSE-AC system
possesses an exponential attractor bounded in $V^{2}\times W^{3}$. This
result was previously reported in \cite{GG1}. Another example covered by the
assumptions of Theorem \ref{expo-thm} is the 3D NS-AC-$\alpha $ system ($%
\theta _{1}=\theta _{2}=\theta =1$) treated previously in \cite{Mej2}.

\begin{remark}
The above results can be used to recover estimates on the dimension of the
global attractor $\mathcal{A}$ for the generalized model (\ref{e:op}),
through the application of the classicial machinery previously used for the
special case of the 2D NSE-AC system, see \cite[Section 4]{GG1}. This is a
somewhat long calculation that we do not include here.
\end{remark}

\subsection{Global and exponential attractors in the case $\protect\theta =0$%
}

\label{ss:attr-nondiss}

In this subsection, we consider non-dissipative systems, which are
represented in our generalized model (\ref{e:op}) when $\theta =0$. This is
the case of the Navier-Stokes-Voight equation which can be seen as an
inviscid regularization of the Navier-Stokes equation. The parabolic
character of the NSV equation is lost; indeed, the semigroup generated by
the NSV equation is only asymptotically compact \cite{KT1, KT2}. In this
sense, NSV behaves more like a damped hyperbolic system, and the same is
true for the NSV-AC system and all other systems when $\theta =0$ (cf. also
Remark \ref{attr-rem}).

Our first task is to prove that the evolution system under consideration
posseses a global attractor, bounded in the energy phase space $\mathcal{Y}%
_{\theta _{2}}$, under rather general conditions on $\theta _{1}\in \mathbb{R%
},$ $\theta _{2}\geq 1$. The analogue of Theorem \ref{t:attr-existb}\ in the
case $\theta =0$ is as follows.

\begin{theorem}
\label{smooth}Let the assumptions of Theorems \ref{t:exist} and \ref{t:stab}
be satisfied for $\theta =0$. For $s\in (\frac{1}{2},1]$, define%
\begin{equation*}
\mathbb{I}_{s}:=(-\infty ,0]\cap (-\infty ,2s-\frac{n}{2})\cap \lbrack
2s-1-2\theta _{2},\infty )\cap (\frac{n}{2}+2s-2-2\theta _{2},\infty )
\end{equation*}%
such that $\mathbb{I}_{s}\neq \varnothing $. In addition, for some%
\begin{equation}
\beta \in \mathbb{J}_{s}:=[\max \left( 1-2\theta _{2},-\theta _{2}\right)
,0]\cap \mathbb{I}_{s},\text{ with }\beta \neq -\theta _{2},
\label{interv-beta1}
\end{equation}%
when $n=2,$ and for some%
\begin{equation}
\beta \in \mathbb{J}_{s}:=(-\theta _{2},-\frac{1}{2}]\cap \mathbb{I}_{s}
\label{interv-beta2}
\end{equation}%
when $n=3$, provided that $\mathbb{J}_{s}\neq \varnothing $, let the
following conditions hold.

(i) $\left\langle A_{0}v,Nv\right\rangle \geq c_{A_{0}}\Vert v\Vert
_{-\theta _{2}}^{2}$ for any $v\in V^{-\theta _{2}}$, with a constant $%
c_{A_{0}}>0$;

(ii) $\left\langle A_{0}v,\left( I-\Delta \right) ^{\beta }v\right\rangle
\geq c_{A_{0}}\Vert v\Vert _{\beta }^{2},$ for any $v\in V^{\beta },$ for
some $c_{A_{0}}>0$;

(iii) $b_{0}(v,w,Nw)=b_{1}(v,\phi ,\phi )=0,$ for any $v,w\in \mathcal{V}$, $%
\phi \in \mathcal{W}$;

(iv) $b_{0}:V^{-\theta _{2}}\times V^{\beta }\times V^{-\beta }\rightarrow 
\mathbb{R}$ is bounded;

(v) $g\in V^{\beta }$ is time independent.

Then, there exists a compact attractor $\mathcal{A}_{0}\Subset \mathcal{Y}%
_{\theta _{2}}$, for the system (\ref{weak1})-(\ref{weak2}) which attracts
the bounded sets of $\mathcal{Y}_{\theta _{2}}$. Moreover, $\mathcal{A}_{0}$
is connected and it is the maximal bounded attractor in $\mathcal{Y}_{\theta
_{2}}$.
\end{theorem}

\begin{proof}
For instance, we follow our argument devised in \cite{GM2} for the
Navier-Stokes-Voight equation with memory. On account of (\ref{semigroup})
and Remark \ref{attr-rem}, let us take a fixed $\varphi _{0}:=\left(
u_{0},\phi _{0}\right) \in \mathcal{B}$ ($\mathcal{B}$ is the absorbing set
derived from Proposition \ref{t:attr-exist}) and consider the corresponding
trajectory $\left( u\left( t\right) ,\phi \left( t\right) \right) =S_{\theta
_{2}}\left( t\right) \varphi _{0}$. Recall that, by Proposition \ref%
{t:attr-exist2}, we have%
\begin{equation}
\left\{ 
\begin{array}{l}
\sup_{\varphi _{0}\in \mathcal{B}}\left\Vert S_{\theta _{2}}\left( t\right)
\varphi _{0}\right\Vert _{\mathcal{Y}_{\theta _{2}}}^{2}\leq C, \\ 
\sup_{\varphi _{0}\in \mathcal{B}}\int_{s}^{t}\left\Vert A_{1}\phi \left(
s\right) \right\Vert _{L^{2}}^{2}ds\leq C\left( 1+t-s\right) ,%
\end{array}%
\right.  \label{decomp-prelim}
\end{equation}%
for all $t>s\geq 0,$ with a constant $C>0$ independent of time and initial
data. Moreover, in this proof the generic constant $C>0$ depends only on $%
\mathcal{B},$ $g$ (see Remark \ref{attr-rem}, (i)) and other physical
parameters of the problem. We split a given trajectory $\left( u\left(
t\right) ,\phi \left( t\right) \right) $ as follows:%
\begin{equation}
\left( u\left( t\right) ,\phi \left( t\right) \right) =\left( u_{d}\left(
t\right) ,\phi _{d}\left( t\right) \right) +\left( u_{c}\left( t\right)
,\phi _{c}\left( t\right) \right) ,  \label{decomp}
\end{equation}%
where%
\begin{equation}
\left\{ 
\begin{array}{l}
\partial _{t}u_{d}+A_{0}u_{d}+B_{0}(u,u_{d})=R_{0}\left( A_{1}\phi _{d},\phi
_{d}\right) , \\ 
\partial _{t}\phi _{d}+B_{1}\left( u_{d},\phi _{d}\right) +A_{1}\phi _{d}=0,
\\ 
u_{d}\left( 0\right) =u_{0},\phi _{d}\left( 0\right) =\phi _{0},%
\end{array}%
\right.  \label{d1}
\end{equation}%
and%
\begin{equation}
\left\{ 
\begin{array}{l}
\partial _{t}u_{c}+A_{0}u_{c}+B_{0}(u,u_{c})=R_{0}\left( A_{1}\phi ,\phi
_{c}\right) +R_{0}\left( A_{1}\phi _{c},\phi _{d}\right) +g, \\ 
\partial _{t}\phi _{c}+B_{1}\left( u,\phi _{c}\right) +B_{1}\left(
u_{c},\phi _{d}\right) +A_{1}\phi _{c}=-f\left( \phi \right) , \\ 
u_{c}\left( 0\right) =0,\phi _{c}\left( 0\right) =0.%
\end{array}%
\right.  \label{c1}
\end{equation}%
Let us first show that $\left( u_{d}\left( t\right) ,\phi _{d}\left(
t\right) \right) $ decays exponentially to zero with respect to the norm of $%
V^{-\theta _{2}}\times W^{1}$. We begin by noting that, since $\left( u,\phi
\right) \in L^{\infty }\left( 0,\infty ;\mathcal{Y}_{\theta _{2}}\right) $
we can easily adapt the proof of Proposition \ref{t:attr-exist} (see also (%
\ref{decomp-prelim})), to find%
\begin{equation}
\sup_{\varphi _{0}\in \mathcal{B}}\left\Vert \left( u_{d}\left( t\right)
,\phi _{d}\left( t\right) \right) \right\Vert _{\mathcal{Y}_{\theta
_{2}}}^{2}+\int_{0}^{\infty }\left( \Vert A_{1}\phi _{d}\left( s\right)
\Vert _{L^{2}}^{2}+\left\Vert u_{d}\left( s\right) \right\Vert _{-\theta
_{2}}^{2}\right) ds\leq C,\text{ for all }t\geq 0,  \label{energyd1}
\end{equation}%
which implies on account of (\ref{energy_fbis}) and (\ref{decomp}), that%
\begin{equation}
\left\{ 
\begin{array}{l}
\sup_{\varphi _{0}\in \mathcal{B}}\left\Vert \left( u_{c}\left( t\right)
,\phi _{c}\left( t\right) \right) \right\Vert _{\mathcal{Y}_{\theta
_{2}}}^{2}\leq C, \\ 
\sup_{\varphi _{0}\in \mathcal{B}}\int_{s}^{t}\Vert A_{1}\phi _{c}\left(
s\right) \Vert _{L^{2}}^{2}ds\leq C\left( 1+t-s\right) ,\text{ for all }%
t>s\geq 0,%
\end{array}%
\right.  \label{energyc1}
\end{equation}%
where $C>0$ is obviously independent of time and initial data (cf. also
Proposition \ref{t:attr-exist2} and (\ref{decomp-prelim})).

Let us now consider the functional%
\begin{equation*}
E_{d}\left( t\right) :=\left\Vert u_{d}\left( t\right) \right\Vert _{-\theta
_{2}}^{2}+\left\Vert A_{1}^{1/2}\phi _{d}\left( t\right) \right\Vert
_{L^{2}}^{2}.
\end{equation*}%
We now pair the first and second equations of (\ref{d1}) with $Nu_{d}$ and $%
A_{1}\phi _{d},$ respectively. Adding together the resulting relations, on
account of assumptions (i)-(ii), we easily derive%
\begin{equation}
\frac{d}{dt}E_{d}\left( t\right) +2\left( c_{A_{0}}\left\Vert u_{d}\left(
t\right) \right\Vert _{-\theta _{2}}^{2}+\left\Vert A_{1}\phi _{d}\left(
t\right) \right\Vert _{L^{2}}^{2}\right) \leq 0.  \label{dest1}
\end{equation}%
Thus, applying a suitable Gronwall's inequality (see Appendix, Lemma \ref%
{Gineq}) to (\ref{dest1}), we obtain%
\begin{equation}
E_{d}\left( t\right) \leq E_{d}\left( 0\right) e^{-\min
\{c_{A_{1}},c_{A_{0}}\}t},\quad \text{\ for all }t\geq 0.  \label{dest2}
\end{equation}%
This estimate gives the desired exponential decay of $\left( u_{d},\phi
_{d}\right) $ in the norm of $V^{-\theta _{2}}\times W^{1}$.

Let us now obtain a bound for $\Vert \left( u_{c},\phi _{c}\right) \Vert _{%
\mathbb{V}_{\beta ,s}}$. To this end, we argue as in the proof of Theorem %
\ref{t:reg}, namely, we pair the first and second equations of (\ref{c1})
with $\Lambda ^{2\beta }u_{c}$ and $A_{1}^{2s}\phi _{c},$ respectively, for $%
s\in (\frac{1}{2},1]$, and then we add the resulting relationships. The
analogue of (\ref{est15bis})-(\ref{est16}) for (\ref{c1}) is the following
identity:%
\begin{align}
& \frac{d}{dt}\left( \left\Vert u_{c}\right\Vert _{\beta }^{2}+\left\Vert
A_{1}^{s}\phi _{c}\right\Vert _{L^{2}}^{2}\right) +2\left\Vert A_{1}^{\left(
2s+1\right) /2}\phi _{c}\right\Vert _{L^{2}}^{2}+2\left\langle
A_{0}u_{c},\Lambda ^{2\beta }u_{c}\right\rangle  \label{ide-comp} \\
& =-2b_{0}\left( u,u_{c},\Lambda ^{2\beta }u_{c}\right) -2\left\langle
A_{1}^{\left( 2s-1\right) /2}f\left( \phi \right) ,A_{1}^{\left( 2s+1\right)
/2}\phi _{c}\right\rangle +2\left\langle g,\Lambda ^{2\beta
}u_{c}\right\rangle  \notag \\
& +2\left\langle R_{0}\left( A_{1}\phi ,\phi _{c}\right) ,\Lambda ^{2\beta
}u_{c}\right\rangle +2\left\langle R_{0}\left( A_{1}\phi _{c},\phi
_{d}\right) ,\Lambda ^{2\beta }u_{c}\right\rangle  \notag \\
& -2\left\langle A_{1}^{\left( 2s-1\right) /2}B_{1}\left( u,\phi _{c}\right)
,A_{1}^{\left( 2s+1\right) /2}\phi _{c}\right\rangle -2\left\langle
A_{1}^{\left( 2s-1\right) /2}B_{1}\left( u_{c},\phi _{d}\right)
,A_{1}^{\left( 2s+1\right) /2}\phi _{c}\right\rangle .  \notag
\end{align}%
Note that, using the fact that $\left( u,\phi \right) \in L^{\infty }\left(
0,\infty ;\mathcal{Y}_{\theta _{2}}\right) $, we can once again find a
positive constant $C_{\delta }\sim 1/\delta $ such that%
\begin{equation}
\left\langle A_{1}^{\left( 2s-1\right) /2}f\left( \phi \right)
,A_{1}^{\left( 2s+1\right) /2}\phi _{c}\right\rangle \leq \delta \left\Vert
A_{1}^{\left( 2s+1\right) /2}\phi _{c}\right\Vert _{L^{2}}^{2}+C\delta
^{-1}\left\Vert \phi \right\Vert _{1}^{2},\text{ for any }\delta >0.
\label{cest0}
\end{equation}%
Moreover, by assumptions (iv)-(v) we have%
\begin{equation}
\begin{array}{l}
2\left\langle B_{0}\left( u,u_{c}\right) ,\Lambda ^{2\beta
}u_{c}\right\rangle \lesssim \left\Vert u\right\Vert _{-\theta
_{2}}\left\Vert u_{c}\right\Vert _{\beta }\left\Vert \Lambda ^{2\beta
}u_{c}\right\Vert _{-\beta }\lesssim \delta ^{-1}\left\Vert u\right\Vert
_{-\theta _{2}}^{2}\left\Vert u_{c}\right\Vert _{\beta }^{2}+\delta
\left\Vert u_{c}\right\Vert _{\beta }^{2}, \\ 
2\left\langle g,\Lambda ^{2\beta }u_{c}\right\rangle \lesssim \delta
^{-1}\left\Vert g\right\Vert _{\beta }^{2}+\delta \left\Vert
u_{c}\right\Vert _{\beta }^{2},%
\end{array}
\label{cest1}
\end{equation}%
for any $\delta >0$. In order to bound the last four terms on the right-hand
side of (\ref{ide-comp}), we argue as follows. By Lemma \ref{l:hole}, the
bilinear mapping $B_{1}:V^{\theta _{2}}\times V^{2s-\varepsilon }\rightarrow
W^{2s-1}$ is continuous provided that $\theta _{2}\geq 1$, for any $s\in (%
\frac{1}{2},1],$ $n=2,3$, and some $\varepsilon \in (0,\frac{1}{2}).$
Moreover, by the Sobolev inequality $\left\Vert \cdot \right\Vert
_{2s+1-\varepsilon }\lesssim \left\Vert \cdot \right\Vert
_{2s+1}^{1-\varepsilon }\left\Vert \cdot \right\Vert _{2s}^{\varepsilon }$,
we have%
\begin{align}
2\left\langle A_{1}^{\left( 2s-1\right) /2}B_{1}\left( u,\phi _{c}\right)
,A_{1}^{\left( 2s+1\right) /2}\phi _{c}\right\rangle & \lesssim \left\Vert
A_{1}^{\left( 2s+1\right) /2}\phi _{c}\right\Vert _{L^{2}}\left\Vert
B_{1}\left( u,\phi _{c}\right) \right\Vert _{2s-1}  \label{cest3} \\
& \lesssim \left\Vert Nu\right\Vert _{\theta _{2}}\left\Vert \nabla \phi
_{c}\right\Vert _{2s-\varepsilon }\left\Vert A_{1}^{\left( 2s+1\right)
/2}\phi _{c}\right\Vert _{L^{2}}  \notag \\
& \lesssim \left\Vert u\right\Vert _{-\theta _{2}}\left\Vert \phi
_{c}\right\Vert _{2s+1}^{2-\varepsilon }\left\Vert \phi _{c}\right\Vert
_{2s}^{\varepsilon }  \notag \\
& \leq \delta \left\Vert A_{1}^{\left( 2s+1\right) /2}\phi _{c}\right\Vert
_{L^{2}}^{2}+C\delta ^{1-2/\varepsilon }\left\Vert u\right\Vert _{-\theta
_{2}}^{\frac{2}{\varepsilon }}\left\Vert A_{1}^{s}\phi _{c}\right\Vert
_{L^{2}}^{2}.  \notag
\end{align}%
Similarly, the mapping $R_{0}$ is continuous from $W^{2s-1}\times
V^{1}\rightarrow V^{\beta },$ and from $W^{0}\times V^{2s}\rightarrow
V^{\beta }$, for as long as $\beta \in \mathbb{I}_{s}\neq \varnothing $ and $%
s\in (\frac{1}{2},1].$ Henceforth, we have the following bounds%
\begin{align}
2\left\langle R_{0}\left( A_{1}\phi ,\phi _{c}\right) ,\Lambda ^{2\beta
}u_{c}\right\rangle & \lesssim \left\Vert \Lambda ^{2\beta }u_{c}\right\Vert
_{-\beta }\left\Vert R_{0}\left( A_{1}\phi ,\phi _{c}\right) \right\Vert
_{\beta }  \label{cest2} \\
& \leq \delta \left\Vert A_{1}^{\left( 2s+1\right) /2}\phi _{c}\right\Vert
_{L^{2}}^{2}+C\delta ^{-1}\left\Vert u_{c}\right\Vert _{\beta
}^{2}\left\Vert A_{1}\phi \right\Vert _{L^{2}}^{2}  \notag
\end{align}%
and%
\begin{align}
2\left\langle R_{0}\left( A_{1}\phi _{c},\phi _{d}\right) ,\Lambda ^{2\beta
}u_{c}\right\rangle & \lesssim \left\Vert u_{c}\right\Vert _{\beta
}\left\Vert R_{0}\left( A_{1}\phi _{c},\phi _{d}\right) \right\Vert _{\beta }
\label{cest2bis} \\
& \leq C\delta ^{-1}\left\Vert u_{c}\right\Vert _{\beta }^{2}\left\Vert
A_{1}\phi _{d}\right\Vert _{L^{2}}^{2}+\delta \left\Vert A_{1}^{\left(
2s+1\right) /2}\phi _{c}\right\Vert _{L^{2}}^{2}.  \notag
\end{align}%
Finally, for the last term we exploit Lemma \ref{l:hole} once more, and
observe that for $\beta \in \mathbb{I}_{s}$, the continuity of the mapping $%
B_{1}:V^{\beta +2\theta _{2}}\times V^{1}\rightarrow W^{2s-1}$ and the
boundedness of the map $N:V^{\beta }\rightarrow V^{\beta +2\theta _{2}}$,
yields%
\begin{align}
\left\langle A_{1}^{\left( 2s-1\right) /2}B_{1}\left( u_{c},\phi _{d}\right)
,A_{1}^{\left( 2s+1\right) /2}\phi _{c}\right\rangle & \lesssim \left\Vert
B_{1}\left( u_{c},\phi _{d}\right) \right\Vert _{2s-1}\left\Vert
A_{1}^{\left( 2s+1\right) /2}\phi _{c}\right\Vert _{L^{2}}  \label{cest4} \\
& \lesssim \left\Vert Nu_{c}\right\Vert _{\beta +2\theta _{2}}\left\Vert
\nabla \phi _{d}\right\Vert _{1}\left\Vert A_{1}^{\left( 2s+1\right) /2}\phi
_{c}\right\Vert _{L^{2}}  \notag \\
& \leq \delta \left\Vert A_{1}^{\left( 2s+1\right) /2}\phi _{c}\right\Vert
_{L^{2}}^{2}+C\delta ^{-1}\left\Vert u_{c}\right\Vert _{\beta
}^{2}\left\Vert A_{1}\phi _{d}\right\Vert _{L^{2}}^{2}.  \notag
\end{align}%
Therefore, by setting%
\begin{align*}
E_{c}\left( t\right) & :=\left\Vert u_{c}\left( t\right) \right\Vert _{\beta
}^{2}+\left\Vert A_{1}^{s}\phi _{c}\left( t\right) \right\Vert _{L^{2}}^{2},
\\
\Xi _{c}\left( t\right) & :=C_{\delta }\left( 1+\left\Vert u\left( t\right)
\right\Vert _{-\theta _{2}}^{2/\varepsilon }+\left\Vert A_{1}\phi \left(
t\right) \right\Vert _{L^{2}}^{2}+\left\Vert A_{1}\phi _{d}\left( t\right)
\right\Vert _{L^{2}}^{2}\right) ,
\end{align*}%
on account of (\ref{ide-comp}), we can choose a sufficiently small $\delta
=(1/8)\min \left( c_{A_{1}},c_{A_{0}}\right) >0$ to deduce%
\begin{align}
& \frac{d}{dt}E_{c}\left( t\right) +C\left( ||A_{1}^{\left( 2s+1\right)
/2}\phi _{c}\left( t\right) ||_{L^{2}}^{2}+\left\Vert u_{c}\left( t\right)
\right\Vert _{\beta }^{2}\right)  \label{ineq8bis} \\
& \leq \Xi _{c}\left( t\right) E_{c}\left( t\right) +C_{\delta }\left(
\left\Vert g\right\Vert _{\beta }^{2}+\left\Vert \phi \left( t\right)
\right\Vert _{1}^{2}\right) ,  \notag
\end{align}%
for all $t\geq 0$. Next, integrate this relation over $\left( 0,t\right) $
and note that $E_{c}\left( 0\right) =0$. Hence, exploiting (\ref{energy_fbis}%
), (\ref{energy_ftris}), (\ref{decomp-prelim}) (\ref{energyd1})-(\ref%
{energyc1}), from the application of Gronwall's lemma we obtain 
\begin{equation}
\left\Vert u_{c}\left( t\right) \right\Vert _{\beta }^{2}+\left\Vert
A_{1}^{s}\phi _{c}\left( t\right) \right\Vert _{L^{2}}^{2}\lesssim
e^{C\left( 1+t\right) },\quad \text{\ for all }t\geq 0.  \label{3.8}
\end{equation}%
Finally, integrating (\ref{ineq8bis}) once more between $0$ and $t$, owing
to (\ref{3.8}) we also find the estimate%
\begin{equation}
\sup_{\varphi _{0}\in \mathcal{B}}\int_{0}^{t}\left\Vert A_{1}^{\left(
2s+1\right) /2}\phi _{c}\left( s\right) \right\Vert _{L^{2}}^{2}ds\lesssim
e^{C\left( 1+t\right) }\text{,}\quad \text{\ for all }t\geq 0.  \label{3.8n}
\end{equation}%
In particular, for every fixed $T>0$, we have found a compact subset $%
V^{\beta }\times D\left( A_{1}^{s}\right) \subset V^{-\theta _{2}}\times
W^{1}$ such that the mapping $S_{\theta _{2}}^{c}\left( t\right) \varphi
_{0}:=\left( u_{c}\left( t\right) ,\phi _{c}\left( t\right) \right) $
satisfies%
\begin{equation}
\bigcup\limits_{\varphi _{0}\in \mathcal{B}}S_{\theta _{2}}^{c}\left(
t\right) \varphi _{0}\subset V^{\beta }\times D\left( A_{1}^{s}\right) \text{%
, for all }t\in \left[ 0,T\right] .  \label{3.8w}
\end{equation}%
This fact, together with the exponential decay (\ref{dest2}), implies that $%
S_{\theta _{2}}\left( t\right) :\mathcal{Y}_{\theta _{2}}\rightarrow 
\mathcal{Y}_{\theta _{2}}$ is asymptotically smooth for $t\geq 0$. The
existence of the global attractor follows by means of standard methods in
the theory of dynamical systems (see, for instance, \cite[Theorem 3.4.6]%
{Hale}). The proof is finished.
\end{proof}

Although Theorem \ref{smooth} yields the global attractor, no conclusion can
be drawn at this stage about its \emph{optimal} regularity. Ideally, one
would like to directly check that $\mathcal{A}_{0}$ is bounded in $V^{\beta
}\times D\left( A_{1}\right) $ as in the statement of Theorem \ref%
{t:attr-existb}. On the contrary, the proof of the above theorem seems to
suggest that this is generally a much harder task, if not out of reach in
just one step. In order to show that $S_{\theta _{2}}\left( t\right) $
enjoys a stronger dissipativity property, we shall employ another semigroup
decomposition, which is much more complicated than the one in (\ref{d1})-(%
\ref{c1}). This step is also crucial in order to demonstrate the existence
of an exponential attractor below.

Our second result establishes the existence of a bounded exponentially
attracting $\mathcal{C}_{\beta ,s}$\ set in $\mathbb{V}_{\beta ,s}:=V^{\beta
}\times D\left( A_{1}^{s}\right) $.

\begin{theorem}
\label{t:glob-nondis}Let the assumptions of Theorems \ref{t:exist} and \ref%
{t:stab} be satisfied for $\theta =0,$ and assume that the conditions (i),
(iii) of Theorem \ref{smooth} are also satisfied. Fix $s\in (\frac{1}{2}%
,1]\cap (0,\frac{\theta _{2}}{2}-\frac{n}{4}+1)$ and consider the nonempty
interval%
\begin{equation*}
\mathbb{K}_{s}:=\mathbb{J}_{s}\cap \left( -\infty ,\frac{1}{2}-\frac{\theta
_{2}}{2}-\frac{n}{4}\right) \cap (-\infty ,s-\frac{1}{2}-\frac{\theta _{2}}{2%
}].
\end{equation*}%
Suppose now that the conditions (ii), (v) of Theorem \ref{smooth} hold for
some $\beta \in \mathbb{K}_{s}$, and the slightly stronger condition:

(vi) $b_{0}:V^{-\theta _{2}}\times V^{-\theta _{2}}\times V^{-\beta
}\rightarrow \mathbb{R}$ is bounded;

There exists $R_{1}>0,$ and a closed ball $\mathcal{B}_{\mathbb{V}_{\beta
,s}}\left( R_{1}\right) \subset \mathbb{V}_{\beta ,s}\cap \mathcal{Y}%
_{\theta _{2}}$ which attracts $\mathcal{B}$ exponentially fast, that is,%
\begin{equation}
dist_{V^{-\theta _{2}}\times W^{1}}\left( S_{\theta _{2}}\left( t\right) 
\mathcal{B},\mathcal{B}_{\mathbb{V}_{\beta ,s}}\left( R_{1}\right) \right)
\leq Ce^{-\rho t},\text{ \ for all }t\geq 0,  \label{exp_attracting}
\end{equation}%
for some positive constants $C$ and $\rho $ independent of time. Here $%
dist_{V^{-\theta _{2}}\times W^{1}}$ denotes the non-symmetric Hausdorff
distance in $V^{-\theta _{2}}\times W^{1}$.
\end{theorem}

\begin{proof}
The main steps require nothing more than what is already contained in the
proof of Theorem \ref{smooth}. Once again, we will rely on the semigroup
decomposition developed there, and the corresponding estimates. First, we
employ another semigroup decomposition and adopt a strategy devised in \cite%
{CP}. For $\varphi _{0}\in \mathcal{B}$, we write $\varphi _{0}=\overline{%
\varphi }_{0}+\overline{\varphi }_{1}$ with $\overline{\varphi }_{0}\in 
\mathcal{B}$ and $\overline{\varphi }_{1}\in V^{\beta }\times D\left(
A_{1}^{s}\right) $, and consider%
\begin{equation}
S_{\theta _{2}}\left( t\right) \varphi _{0}=U_{\varphi _{0}}\left( t\right) 
\overline{\varphi }_{0}+V_{\varphi _{0}}\left( t\right) \overline{\varphi }%
_{1}.  \label{new-decomp}
\end{equation}%
We define $U_{\varphi _{0}}\left( t\right) \overline{\varphi }_{0}=\overline{%
\varphi }_{d}\left( t\right) $ and $V_{\varphi _{0}}\left( t\right) 
\overline{\varphi }_{1}=\overline{\varphi }_{c}\left( t\right) ,$ where $%
\overline{\varphi }_{d}\left( t\right) =\left( \overline{u}_{d}\left(
t\right) ,\overline{\phi }_{d}\left( t\right) \right) $ and $\overline{%
\varphi }_{c}\left( t\right) =\left( \overline{u}_{c}\left( t\right) ,%
\overline{\phi }_{c}\left( t\right) \right) $ solve the Cauchy problems:%
\begin{equation}
\left\{ 
\begin{array}{l}
\partial _{t}\overline{u}_{d}+A_{0}\overline{u}_{d}+B_{0}(u,\overline{u}%
_{d})=R_{0}\left( A_{1}\overline{\phi }_{d},\phi _{c}\right) , \\ 
\partial _{t}\overline{\phi }_{d}+B_{1}\left( \overline{u}_{d},\phi
_{c}\right) +A_{1}\overline{\phi }_{d}=0, \\ 
\overline{u}_{d}\left( 0\right) =\overline{u}_{0},\overline{\phi }_{d}\left(
0\right) =\overline{\phi }_{0}%
\end{array}%
\right.  \label{d2}
\end{equation}%
and%
\begin{equation}
\left\{ 
\begin{array}{l}
\partial _{t}\overline{u}_{c}+A_{0}\overline{u}_{c}+B_{0}(u,\overline{u}%
_{c})=R_{0}\left( A_{1}\phi ,\phi _{d}\right) +R_{0}\left( A_{1}\overline{%
\phi }_{c},\phi _{c}\right) +g, \\ 
\partial _{t}\overline{\phi }_{c}+B_{1}\left( u,\phi _{d}\right)
+B_{1}\left( \overline{u}_{c},\phi _{c}\right) +A_{1}\overline{\phi }%
_{c}=-f\left( \phi \right) , \\ 
\overline{u}_{c}\left( 0\right) =\overline{u}_{1},\overline{\phi }_{c}\left(
0\right) =\overline{\phi }_{1}.%
\end{array}%
\right.  \label{c2}
\end{equation}%
Concerning the mapping $U_{\varphi _{0}}$, for every $\overline{\varphi }%
_{0}\in \mathcal{B}$ one argues exactly as in (\ref{dest1})-(\ref{dest2}) to
deduce%
\begin{equation}
\left\Vert U_{\varphi _{0}}\left( t\right) \overline{\varphi }%
_{0}\right\Vert _{\mathcal{Y}_{\theta _{2}}}^{2}\leq e^{-\rho t}\left\Vert 
\overline{\varphi }_{0}\right\Vert _{\mathcal{Y}_{\theta _{2}}}^{2}\leq
R_{0}e^{-\rho t},  \label{expo-stability}
\end{equation}%
for all $t\geq 0$, for some $\rho >0.$ The final step is to deduce an energy
inequality for $\overline{\varphi }_{c}\left( t\right) .$ To this end, we
set 
\begin{align*}
\Lambda _{1}\left( t\right) & :=\left\Vert \overline{u}_{c}\left( t\right)
\right\Vert _{\beta }^{2}+\left\Vert A_{1}^{s}\overline{\phi }_{c}\left(
t\right) \right\Vert _{L^{2}}^{2}, \\
\Lambda _{2}\left( t\right) & :=C\left\Vert u\left( t\right) \right\Vert
_{-\theta _{2}}^{2}\left( \left\Vert \overline{u}_{c}\left( t\right)
\right\Vert _{-\theta _{2}}^{2}+\left\Vert \phi _{d}\left( t\right)
\right\Vert _{2}^{2}\right) +C\left\Vert g\left( t\right) \right\Vert
_{\beta }^{2}+\left\Vert \phi \left( t\right) \right\Vert _{1}^{2} \\
& +C\left\Vert \overline{u}_{c}\left( t\right) \right\Vert _{-\theta
_{2}}^{2}\left\Vert \phi _{c}\left( t\right) \right\Vert
_{2s+1}^{2}+C\left\Vert \overline{u}_{c}\left( t\right) \right\Vert
_{-\theta _{2}}\left\Vert A_{1}\phi \left( t\right) \right\Vert
_{L^{2}}\left\Vert A_{1}\phi _{d}\left( t\right) \right\Vert _{L^{2}}
\end{align*}%
and then observe, on account of (\ref{new-decomp}), (\ref{expo-stability})
and (\ref{decomp-prelim}), (\ref{energyd1}), (\ref{3.8n}), that%
\begin{equation}
\sup_{\varphi _{0}\in \mathcal{B}}\int_{0}^{t}\Lambda _{2}\left( s\right)
ds\lesssim e^{C\left( 1+t\right) }\text{,}\quad \text{\ for all }t\geq 0.
\label{c3}
\end{equation}%
Once again we pair the first and second equations of (\ref{c2}) with $%
\Lambda ^{2\beta }\overline{u}_{c}$ and $A_{1}^{2s}\overline{\phi }_{c},$
respectively. Adding the resulting relationships, we get the following
identity:%
\begin{align}
& \frac{d\Lambda _{1}}{dt}+2\left\Vert A_{1}^{\left( 2s+1\right) /2}%
\overline{\phi }_{c}\right\Vert _{L^{2}}^{2}+2\left\langle A_{0}\overline{u}%
_{c},\Lambda ^{2\beta }\overline{u}_{c}\right\rangle  \label{eqq10} \\
& =-2b_{0}\left( u,\overline{u}_{c},\Lambda ^{2\beta }\overline{u}%
_{c}\right) -2\left\langle A_{1}^{\left( 2s-1\right) /2}f\left( \phi \right)
,A_{1}^{\left( 2s+1\right) /2}\overline{\phi }_{c}\right\rangle
+2\left\langle g,\Lambda ^{2\beta }\overline{u}_{c}\right\rangle  \notag \\
& +2\left\langle R_{0}\left( A_{1}\phi ,\phi _{d}\right) ,\Lambda ^{2\beta }%
\overline{u}_{c}\right\rangle +2\left\langle R_{0}\left( A_{1}\overline{\phi 
}_{c},\phi _{c}\right) ,\Lambda ^{2\beta }\overline{u}_{c}\right\rangle 
\notag \\
& -2\left\langle A_{1}^{\left( 2s-1\right) /2}B_{1}\left( u,\phi _{d}\right)
,A_{1}^{\left( 2s+1\right) /2}\overline{\phi }_{c}\right\rangle
-2\left\langle A_{1}^{\left( 2s-1\right) /2}B_{1}\left( \overline{u}%
_{c},\phi _{c}\right) ,A_{1}^{\left( 2s+1\right) /2}\overline{\phi }%
_{c}\right\rangle .  \notag
\end{align}%
As in the proof of Theorem \ref{smooth}, for every $\delta >0$ we have%
\begin{equation}
\left\langle A_{1}^{\left( 2s-1\right) /2}f\left( \phi \right)
,A_{1}^{\left( 2s+1\right) /2}\overline{\phi }_{c}\right\rangle \lesssim
\delta \left\Vert \overline{\phi }_{c}\right\Vert _{2s+1}^{2}+\delta
^{-1}\left\Vert \phi \right\Vert _{1}^{2}  \label{eqq11}
\end{equation}%
and%
\begin{equation}
2\left\langle B_{0}\left( u,\overline{u}_{c}\right) -g,\Lambda ^{2\beta }%
\overline{u}_{c}\right\rangle \lesssim \delta ^{-1}\left( \left\Vert
u\right\Vert _{-\theta _{2}}^{2}\left\Vert \overline{u}_{c}\right\Vert
_{-\theta _{2}}^{2}+\left\Vert g\right\Vert _{\beta }^{2}\right) +\delta
\left\Vert \overline{u}_{c}\right\Vert _{\beta }^{2},  \label{eqq12}
\end{equation}%
owing to assumption (vi). Exactly as in (\ref{cest3}), it follows%
\begin{align}
2\left\langle A_{1}^{\left( 2s-1\right) /2}B_{1}\left( \overline{u}_{c},\phi
_{c}\right) ,A_{1}^{\left( 2s+1\right) /2}\overline{\phi }_{c}\right\rangle
& \lesssim \left\Vert N\overline{u}_{c}\right\Vert _{\theta _{2}}\left\Vert
\nabla \phi _{c}\right\Vert _{2s-\varepsilon }\left\Vert \overline{\phi }%
_{c}\right\Vert _{2s+1}  \label{eqq13} \\
& \lesssim \delta \left\Vert \overline{\phi }_{c}\right\Vert
_{2s+1}^{2}+\delta ^{-1}\left\Vert \overline{u}_{c}\right\Vert _{-\theta
_{2}}^{2}\left\Vert \phi _{c}\right\Vert _{2s+1}^{2}.  \notag
\end{align}%
Next, for $2s<\theta _{2}-n/2+2,$ the boundedness of the mapping $%
B_{1}:V^{\theta _{2}}\times V^{1}\rightarrow W^{2s-1}$ yields%
\begin{align}
2\left\langle A_{1}^{\left( 2s-1\right) /2}B_{1}\left( u,\phi _{d}\right)
,A_{1}^{\left( 2s+1\right) /2}\overline{\phi }_{c}\right\rangle & \lesssim
\left\Vert Nu\right\Vert _{\theta _{2}}\left\Vert \nabla \phi
_{d}\right\Vert _{1}\left\Vert \overline{\phi }_{c}\right\Vert _{2s+1}
\label{eqq14} \\
& \lesssim \delta \left\Vert \overline{\phi }_{c}\right\Vert
_{2s+1}^{2}+\delta ^{-1}\left\Vert u\right\Vert _{-\theta
_{2}}^{2}\left\Vert \phi _{d}\right\Vert _{2}^{2}.  \notag
\end{align}%
For the remaining two terms, we use the boundedness of the mapping $R_{0}$
from $W^{2s-1}\times V^{2s}\rightarrow V^{2\beta +\theta _{2}}$ and $%
W^{0}\times V^{1}\rightarrow V^{2\beta +\theta _{2}},$ for every $\beta \in 
\mathbb{K}_{s}$ to derive%
\begin{align}
& 2\left\langle R_{0}\left( A_{1}\phi ,\phi _{d}\right) ,\Lambda ^{2\beta }%
\overline{u}_{c}\right\rangle +2\left\langle R_{0}\left( A_{1}\overline{\phi 
}_{c},\phi _{c}\right) ,\Lambda ^{2\beta }\overline{u}_{c}\right\rangle
\label{eqq15} \\
& \lesssim \delta \left\Vert \overline{\phi }_{c}\right\Vert
_{2s+1}^{2}+\delta ^{-1}\left\Vert \overline{u}_{c}\right\Vert _{-\theta
_{2}}^{2}\left\Vert \phi _{c}\right\Vert _{2s+1}^{2}+\left\Vert \overline{u}%
_{c}\right\Vert _{-\theta _{2}}\left\Vert A_{1}\phi \right\Vert
_{L^{2}}\left\Vert A_{1}\phi _{d}\right\Vert _{L^{2}}.  \notag
\end{align}%
Finally, taking $\delta =\frac{1}{4}\min \left( c_{A_{1}},c_{A_{0}}\right)
>0 $ sufficiently small, the basic energy identity (\ref{eqq10}) takes the
form of an inequality%
\begin{equation*}
\frac{d}{dt}\Lambda _{1}\left( t\right) +4\delta \Lambda _{1}\left( t\right)
\lesssim \Lambda _{2}\left( t\right) ,\text{ for }t\geq 0.
\end{equation*}%
In view of (\ref{c3}), the classical Gronwall lemma gives%
\begin{equation}
\left\Vert V_{\varphi _{0}}\left( t\right) \overline{\varphi }%
_{1}\right\Vert _{\mathbb{V}_{\beta ,s}}^{2}=\Lambda _{1}\left( t\right)
\leq e^{-4\delta t}\left\Vert \overline{\varphi }_{1}\right\Vert _{\mathbb{V}%
_{\beta ,s}}^{2}+J\left( t\right) ,  \label{map-stability}
\end{equation}%
for some positive continuous\ function $J:\mathbb{R}_{+}\rightarrow \mathbb{R%
}_{+},$ $J\left( t\right) \sim e^{Ct}$ as $t\rightarrow \infty .$ The thesis
(\ref{exp_attracting}) then follows from (\ref{expo-stability}), (\ref%
{map-stability}) and the application of \cite[Theorem 3.1]{CP}, which
finishes the proof.
\end{proof}

As a byproduct of the previous result, we also obtain a regularity result
for the global attractor $\mathcal{A}_{0}$ in the case $\theta =0$.

\begin{corollary}
\label{co-b}The global attractor $\mathcal{A}_{0}$ is bounded in $V^{\beta
}\times D\left( A_{1}^{s}\right) $, for some $\beta \in \mathbb{K}_{s}$ and $%
s\in (\frac{1}{2},1]\cap (0,\frac{\theta _{2}}{2}-\frac{n}{4}+1).$
\end{corollary}

\begin{remark}
For example, setting $\theta =0$ and $\theta _{2}=\theta _{1}=1$, and
checking all the requirements of Theorem~\ref{t:glob-nondis}, the 3D NSV-AC
admits a global attractor in the sense of Corollary \ref{co-b}, provided
that $s\in (0.5,0.75)$. As far as we know this result was not reported
anywhere else. On the other hand, once the regularity in $\mathcal{C}_{\beta
,s}:=\mathcal{B}_{\mathbb{V}_{\beta ,s}}\left( R_{1}\right) $ is
established, it may be also possible to exploit the semigroup decomposition (%
\ref{d1})-(\ref{c1}) and a bootstrap argument to show that, for $\left(
u,\phi \right) \in \mathcal{A}_{0}$, it also holds $\phi \in D\left(
A_{1}\right) $. We omit the details in order to avoid further technicalities.
\end{remark}

\begin{remark}
Note that estimate (\ref{exp_attracting}) is also satisfied provided that
the external force $g$ is time dependent and $g\in L_{tb}^{2}\left( \mathbb{R%
}_{+};V^{\beta }\right) .$ Even in the case $\theta =0$, one can generalize
the notion of global attractor and replace it by the notion of pullback
attractor, see Remark \ref{attr-rem}, (ii).
\end{remark}

In the second part, we show the existence of exponential attractors for our
regularized family of models (\ref{e:op}) when $\theta =0$. First, we note
the following straight-forward proposition.

\begin{proposition}
\label{t:to-H1} Let the assumptions of Theorem \ref{t:glob-nondis}. There
exists a time $t_{3}>0$ such that $S_{\theta _{2}}\left( t\right) \mathcal{C}%
_{\beta ,s}\subset \mathcal{C}_{\beta ,s},$ for all $t\geq t_{3}.$ Moreover,
the following estimate holds:%
\begin{equation}
\sup_{\varphi _{0}\in \mathcal{C}_{\beta ,s}}\left\Vert \left( u\left(
t\right) ,\phi \left( t\right) \right) \right\Vert _{\mathbb{V}_{\beta
,s}}^{2}+\int_{t}^{t+1}\left\Vert \phi \left( s\right) \right\Vert
_{2s+1}^{2}ds\lesssim R_{1},  \label{higher-est}
\end{equation}%
for all $t\geq 0$, for some $\beta \in \mathbb{K}_{s}$ and $s\in (\frac{1}{2}%
,1]\cap (0,\frac{\theta _{2}}{2}-\frac{n}{4}+1).$
\end{proposition}

The next lemma is concerned with the H\"{o}lder-in-time regularity of the
semigroup $S_{\theta _{2}}(t).$

\begin{lemma}
\label{time_reg} Let the assumptions of Proposition \ref{t:to-H1} be
satisfied. There exists a constant $C>0$ such that%
\begin{equation}
\left\Vert S_{\theta _{2}}(t)\varphi _{0}-S_{\theta _{2}}(\tilde{t})\varphi
_{0}\right\Vert _{V^{-\theta _{2}}\times W^{2s-1}}\leq C\left\vert t-%
\widetilde{t}\right\vert ^{1/2},  \label{time_regularity}
\end{equation}%
for all $t,\tilde{t}\in \left[ 0,\infty \right) $ and any $\varphi
_{0}=\left( u_{0},\varphi _{0}\right) \in \mathcal{C}_{\beta ,s}\subset
V^{-\theta _{2}}\times W^{2s-1}$.
\end{lemma}

\begin{proof}
The following basic equality plays an essential role:%
\begin{equation}
S_{\theta _{2}}(t)\varphi _{0}-S_{\theta _{2}}(\tilde{t})\varphi
_{0}=:\left( u\left( t\right) -u\left( \widetilde{t}\right) ,\phi \left(
t\right) -\phi \left( \widetilde{t}\right) \right) =\int_{t}^{\widetilde{t}%
}\partial _{y}\left( S_{\theta _{2}}(y)\varphi _{0}\right) dy.  \label{diff}
\end{equation}%
We claim that the bounds below hold for any strong solution $\left( u,\phi
\right) $ of problem (\ref{e:op}). More precisely, we have%
\begin{equation}
\int_{0}^{\infty }\left\Vert \partial _{t}u\left( s\right) \right\Vert
_{-\theta _{2}}^{2}ds\leq C\left( R_{1}\right) ,\text{ }\int_{0}^{\infty
}\left\Vert \partial _{t}\phi \left( s\right) \right\Vert _{2s-1}^{2}ds\leq
C\left( R_{1}\right) .  \label{timebis}
\end{equation}%
As usual in order to rigorously justify (\ref{timebis}) we can appeal once
more to the approximation scheme used in the proof of Theorem \ref{t:exist}.
From (\ref{e:op}), and by the boundedness of $A_{0}:V^{-\theta
_{2}}\rightarrow V^{-\theta _{2}}$, we have%
\begin{equation}
\left\{ 
\begin{array}{l}
\left\Vert \partial _{t}u\right\Vert _{-\theta _{2}}\lesssim \left\Vert
u\right\Vert _{-\theta _{2}}+\left\Vert B_{0}\left( u,u\right) \right\Vert
_{-\theta _{2}}+\left\Vert R_{0}\left( A_{1}\phi ,\phi \right) \right\Vert
_{-\theta _{2}}+\left\Vert g\right\Vert _{-\theta _{2}}, \\ 
\left\Vert \partial _{t}\phi \right\Vert _{2s-1}\leq \left\Vert A_{1}\phi
+f\left( \phi \right) \right\Vert _{2s-1}+\left\Vert B_{1}\left( u,\phi
\right) \right\Vert _{2s-1}.%
\end{array}%
\right.  \label{time2}
\end{equation}%
In order to bound all terms on the right-hand side of (\ref{time2}), we
appeal to the same estimates derived in the proof of Theorem \ref{t:exist}
and Theorem \ref{smooth}. Indeed, exactly as in (\ref{exist-1})-(\ref{exist3}%
), and exploiting the fact that $R_{0}:W^{2s-2}\times V^{2s}\rightarrow
V^{-1}$ is bounded for any $s\in (\frac{1}{2},1]$, we have%
\begin{equation*}
\left\Vert B_{0}\left( u,u\right) \right\Vert _{-\theta _{2}}\lesssim \Vert
u\Vert _{-\theta _{2}}^{\lambda _{1}+\lambda _{2}},\text{ }\left\Vert
R_{0}\left( A_{1}\phi ,\phi \right) \right\Vert _{-\theta _{2}}\lesssim
\left\Vert R_{0}\left( A_{1}\phi ,\phi \right) \right\Vert _{-1}\lesssim
\left\Vert \phi \right\Vert _{2s}||\phi ||_{2s+1}
\end{equation*}%
where we recall that $\lambda =\lambda _{1}+\lambda _{2}=1$ if $\theta =0$.
Thus, we get thanks to (\ref{higher-est}) that $\partial _{t}u\in
L^{2}\left( 0,\infty ;V^{-\theta _{2}}\right) $ and the first estimate of (%
\ref{timebis}) holds. On the other hand, recalling that for $\beta \in 
\mathbb{J}_{s},$ the map $B_{1}:V^{\beta +2\theta _{2}}\times
V^{1}\rightarrow W^{2s-1}$ is continuous, we can also bound the $W^{2s-1}$%
-norm of $B_{1}\left( u,\phi \right) $ in the last inequality of (\ref{time2}%
). Hence, owing to (\ref{higher-est}) we obtain $\partial _{t}\phi \in
L^{2}(0,\infty ;W^{2s-1})$ such that the second inequality of (\ref{timebis}%
) is satisfied. Thanks to (\ref{timebis}), from (\ref{diff}) and by
application of H\"{o}lder's inequality, we infer%
\begin{equation}
\left\Vert u\left( t\right) -u\left( \widetilde{t}\right) \right\Vert
_{-\theta _{2}}\leq C|t-\tilde{t}|^{1/2},\text{ }\left\Vert \phi \left(
t\right) -\phi \left( \widetilde{t}\right) \right\Vert _{2s-1}\leq C|t-%
\tilde{t}|^{1/2}.  \label{time3}
\end{equation}%
Thus, we immediately arrive at the inequality (\ref{time_regularity}) owing
to (\ref{time3}). The proof is finished.
\end{proof}

The second main result of this subsection is the following.

\begin{theorem}
\label{expo-thm2}Let the assumptions of Theorem \ref{t:glob-nondis} be
satisfied. Then $\left( S_{\theta _{2}},\mathcal{Y}_{\theta _{2}}\right) $
possesses an exponential attractor $\mathcal{M}_{0}\subset \mathcal{Y}%
_{\theta _{2}}$ which is bounded in $\mathbb{V}_{\beta ,s}$. Thus, by
definition, we have

(a) $\mathcal{M}_{0}$ is compact and semi-invariant with respect $S_{\theta
_{2}}\left( t\right) ,$ that is,%
\begin{equation*}
S_{\theta _{2}}\left( t\right) \left( \mathcal{M}_{0}\right) \subseteq 
\mathcal{M}_{0},\quad \forall \,t\geq 0.
\end{equation*}

(b) The fractal dimension $\dim _{F}\left( \mathcal{M}_{0},\mathcal{Y}%
_{\theta _{2}}\right) $ of $\mathcal{M}_{0}$ is finite and an upper bound
can be computed explicitly.

(c) $\mathcal{M}_{0}$ attracts exponentially fast any bounded subset $B$ of $%
\mathcal{Y}_{\theta _{2}}$, that is, there exist a positive nondecreasing
function $Q$ and a constant $\rho >0$ such that 
\begin{equation*}
dist_{\mathcal{Y}_{\theta _{2}}}\left( S_{\theta _{2}}\left( t\right) B,%
\mathcal{M}_{0}\right) \leq Q(\Vert B\Vert _{\mathcal{Y}_{\theta
_{2}}})e^{-\rho t},\quad \forall \,t\geq 0.
\end{equation*}
\end{theorem}

\begin{proof}
\emph{Step 1} (The smoothing property). For $\varphi _{0i}=\left(
u_{0i},\phi _{0i}\right) \in \mathcal{C}_{\beta ,s},$ let $\varphi
_{i}=\left( u_{i},\phi _{i}\right) =S_{\theta _{2}}\left( t\right) \varphi
_{0i},$ $i=1,2$, be the corresponding solutions. We decompose $\varphi
\left( t\right) :=\left( u\left( t\right) ,\phi \left( t\right) \right)
=\varphi _{1}\left( t\right) -\varphi _{2}\left( t\right) ,$ such that%
\begin{equation}
\varphi \left( t\right) =\left( v\left( t\right) ,\rho \left( t\right)
\right) +\left( \omega \left( t\right) ,\psi \left( t\right) \right) ,
\label{decomp-diff}
\end{equation}%
where $\left( v,\rho \right) $ solves%
\begin{equation}
\left\{ 
\begin{array}{l}
\partial _{t}v+A_{0}v+B_{0}(u_{2},v)=R_{0}\left( A_{1}\rho ,\phi _{1}\right)
, \\ 
\partial _{t}\rho +A_{1}\rho +B_{1}\left( v,\phi _{1}\right) =0, \\ 
v\left( 0\right) =u_{01}-u_{02},\rho \left( 0\right) =\phi _{01}-\phi _{02},%
\end{array}%
\right.  \label{conc-diff1}
\end{equation}%
and $\left( \omega ,\psi \right) $ solves%
\begin{equation}
\left\{ 
\begin{array}{l}
\partial _{t}\omega +A_{0}\omega +B_{0}(u,u_{1})+B_{0}\left( u_{2},\omega
\right) =R_{0}\left( A_{1}\psi ,\phi _{1}\right) +R_{0}\left( A_{1}\phi
_{2},\phi \right) , \\ 
\partial _{t}\psi +B_{1}\left( \omega ,\phi _{1}\right) +B_{1}\left(
u_{2},\phi \right) +A_{1}\psi =-\left( f\left( \phi _{1}\right) -f\left(
\phi _{2}\right) \right) ,%
\end{array}%
\right.  \label{conc-diff2}
\end{equation}%
supplemented with null initial data. It is apparent that upon pairing the
first and second equations of (\ref{conc-diff1}) with $Nv$ and $A_{1}\rho ,$
respectively, the norm $\left\Vert \left( v,\rho \right) \right\Vert
_{V^{-\theta _{2}}\times W^{1}}$ is exponentially decaying to zero. More
precisely, we easily get%
\begin{equation}
\left\Vert \left( v\left( t\right) ,\rho \left( t\right) \right) \right\Vert
_{V^{-\theta _{2}}\times W^{1}}^{2}\leq e^{-\kappa t}\left\Vert \varphi
_{01}-\varphi _{02}\right\Vert _{V^{-\theta _{2}}\times W^{1}}^{2},
\label{eq35}
\end{equation}%
for all $t\geq 0,$ for some positive constant $\kappa $ independent of time.
Concerning the other component in (\ref{decomp-diff}), we pair the first and
second equations of (\ref{conc-diff2}) with $\Lambda ^{2\beta }\omega $ and $%
A_{1}^{2s}\psi ,$ respectively, for $s\in (\frac{1}{2},1]$, and then we add
the resulting relationships. We have%
\begin{align}
& \frac{d}{dt}\left( \left\Vert \omega \right\Vert _{\beta }^{2}+\left\Vert
A_{1}^{s}\psi \right\Vert _{L^{2}}^{2}\right) +2\left\Vert A_{1}^{\left(
2s+1\right) /2}\psi \right\Vert _{L^{2}}^{2}+2\left\langle A_{0}\omega
,\Lambda ^{2\beta }\omega \right\rangle  \label{eq36} \\
& =-2b_{0}\left( u,u_{1},\Lambda ^{2\beta }\omega \right) -2b_{0}\left(
u_{2},\omega ,\Lambda ^{2\beta }\omega \right)  \notag \\
& -2\left\langle A_{1}^{\left( 2s-1\right) /2}\left( f\left( \phi
_{1}\right) -f\left( \phi _{2}\right) \right) ,A_{1}^{\left( 2s+1\right)
/2}\psi \right\rangle  \notag \\
& +2\left\langle R_{0}\left( A_{1}\psi ,\phi _{1}\right) ,\Lambda ^{2\beta
}\omega \right\rangle +2\left\langle R_{0}\left( A_{1}\phi _{2},\phi \right)
,\Lambda ^{2\beta }\omega \right\rangle  \notag \\
& -2\left\langle A_{1}^{\left( 2s-1\right) /2}B_{1}\left( \omega ,\phi
_{1}\right) ,A_{1}^{\left( 2s+1\right) /2}\psi \right\rangle -2\left\langle
A_{1}^{\left( 2s-1\right) /2}B_{1}\left( u_{2},\phi \right) ,A_{1}^{\left(
2s+1\right) /2}\psi \right\rangle .  \notag
\end{align}%
At this point, we basically repeat the proof of Theorem \ref{smooth} up to
the estimate (\ref{cest4}). Note that, using the fact that any trajectory $%
\varphi _{i},$ $i=1,2,$ satisfies (\ref{higher-est}) we can once again find
a constant $C>0$ such that%
\begin{equation}
\left\langle A_{1}^{\left( 2s-1\right) /2}\left( f\left( \phi _{1}\right)
-f\left( \phi _{2}\right) \right) ,A_{1}^{\left( 2s+1\right) /2}\psi
\right\rangle \leq \delta ||A_{1}^{\left( 2s+1\right) /2}\psi
||_{L^{2}}^{2}+C\delta ^{-1}\left\Vert \phi \right\Vert _{1}^{2},\text{ for
any }\delta >0.  \label{eq37}
\end{equation}%
Moreover, exactly as in (\ref{cest1}) it follows that%
\begin{align}
2b_{0}\left( u,u_{1},\Lambda ^{2\beta }\omega \right) & \lesssim \delta
^{-1}\left\Vert u\right\Vert _{-\theta _{2}}^{2}\left\Vert u_{1}\right\Vert
_{\beta }^{2}+\delta \left\Vert \omega \right\Vert _{\beta }^{2},
\label{eq38} \\
2b_{0}\left( u_{2},\omega ,\Lambda ^{2\beta }\omega \right) & \lesssim
\delta ^{-1}\left\Vert u_{2}\right\Vert _{-\theta _{2}}^{2}\left\Vert \omega
\right\Vert _{\beta }^{2}+\delta \left\Vert \omega \right\Vert _{\beta }^{2}.
\notag
\end{align}%
Owing to (\ref{cest2})-(\ref{cest2bis}), we infer%
\begin{equation}
2\left\langle R_{0}\left( A_{1}\psi ,\phi _{1}\right) ,\Lambda ^{2\beta
}\omega \right\rangle \lesssim \delta \left\Vert \omega \right\Vert _{\beta
}^{2}\left\Vert A_{1}\phi _{1}\right\Vert _{L^{2}}^{2}+\delta
^{-1}\left\Vert A_{1}^{\left( 2s+1\right) /2}\psi \right\Vert _{L^{2}}^{2}
\label{eq39}
\end{equation}%
and, respectively,%
\begin{equation}
2\left\langle R_{0}\left( A_{1}\phi _{2},\phi \right) ,\Lambda ^{2\beta
}\omega \right\rangle \lesssim \left\Vert \omega \right\Vert _{\beta
}^{2}\left\Vert A_{1}^{\left( 2s+1\right) /2}\phi _{2}\right\Vert
_{L^{2}}^{2}+\left\Vert A_{1}\phi \right\Vert _{L^{2}}^{2}.  \label{eq39bis}
\end{equation}%
Finally, by virtue of (\ref{cest3})-(\ref{cest4}) the following inequalities
also hold:%
\begin{equation}
2\left\langle A_{1}^{\left( 2s-1\right) /2}B_{1}\left( \omega ,\phi
_{1}\right) ,A_{1}^{\left( 2s+1\right) /2}\psi \right\rangle \lesssim \delta
\left\Vert A_{1}^{\left( 2s+1\right) /2}\psi \right\Vert _{L^{2}}^{2}+\delta
^{-1}\left\Vert \omega \right\Vert _{\beta }^{2}\left\Vert A_{1}\phi
_{1}\right\Vert _{L^{2}}^{2},  \label{eq40}
\end{equation}%
\begin{equation}
2\left\langle A_{1}^{\left( 2s-1\right) /2}B_{1}\left( u_{2},\phi \right)
,A_{1}^{\left( 2s+1\right) /2}\psi \right\rangle \lesssim \delta \left\Vert
A_{1}^{\left( 2s+1\right) /2}\psi \right\Vert _{L^{2}}^{2}+\delta
^{-1}\left\Vert u_{2}\right\Vert _{\beta }^{2}\left\Vert A_{1}\phi
\right\Vert _{L^{2}}^{2}.  \label{eq41}
\end{equation}%
Let us now set%
\begin{align*}
\mathcal{P}\left( t\right) & :=\left\Vert \omega \left( t\right) \right\Vert
_{\beta }^{2}+\left\Vert A_{1}^{s}\psi \left( t\right) \right\Vert
_{L^{2}}^{2}, \\
\mathcal{N}_{1}\left( t\right) & :=C_{\delta }\left( \left\Vert u_{2}\left(
t\right) \right\Vert _{-\theta _{2}}^{2}+\left\Vert \phi _{1}\right\Vert
_{2s+1}^{2}+\left\Vert \phi _{2}\right\Vert _{2s+1}^{2}\right) , \\
\mathcal{N}_{2}\left( t\right) & :=C_{\delta }\left( 1+\left\Vert
u_{2}\left( t\right) \right\Vert _{\beta }^{2}\right) \left\Vert A_{1}\phi
_{1}\left( t\right) -A_{1}\phi _{2}\left( t\right) \right\Vert
_{L^{2}}^{2}+C_{\delta }\left\Vert u\left( t\right) \right\Vert _{-\theta
_{2}}^{2}\left\Vert u_{1}\left( t\right) \right\Vert _{\beta }^{2}.
\end{align*}%
On account of (\ref{eq36})-(\ref{eq41}), we can choose a sufficiently small $%
\delta \sim \min \left( c_{A_{0}},c_{A_{1}}\right) >0$ to deduce%
\begin{equation}
\frac{d}{dt}\mathcal{P}\left( t\right) \leq \mathcal{N}_{1}\left( t\right) 
\mathcal{P}\left( t\right) +\mathcal{N}_{2}\left( t\right) ,
\end{equation}%
for all $t\geq 0$. Hence, integrating (\ref{ineq8bis}) with respect to time
on the interval $\left( 0,t\right) $, noting that $\mathcal{P}\left(
0\right) =0,$ and exploiting (\ref{uniq_stab}), (\ref{higher-est}), we can
find a positive continuous function $\lambda :\mathbb{R}_{+}\rightarrow 
\mathbb{R}_{+},$ $\lambda \left( 0\right) >0$, such that%
\begin{equation}
\left\Vert \left( \omega \left( t\right) ,\psi \left( t\right) \right)
\right\Vert _{\mathbb{V}_{\beta ,s}}^{2}=\mathcal{P}\left( t\right) \leq
\lambda \left( t\right) \left\Vert \varphi _{01}-\varphi _{02}\right\Vert
_{V^{-\theta _{2}}\times W^{1}}^{2},  \label{eq42}
\end{equation}%
for all $t\geq 0$.

\emph{Step 2} (The final argument). Thanks to Theorem \ref{t:glob-nondis},
Propositions \ref{t:attr-exist}, \ref{t:to-H1} and the transitivity property
of the exponential attraction \cite[Theorem 5.1]{FGMZ04}, the set $\mathcal{C%
}_{\beta ,s}$ is positively invariant for $S_{\theta _{2}}\left( t\right) $
and attracts any bounded set of $\mathcal{Y}_{\theta _{2}}$ exponentially
fast. Moreover, by Theorem \ref{t:stab} and Lemma \ref{time_reg} the map $%
\left( t,\varphi _{0}\right) \mapsto S_{\theta _{2}}\left( t\right) \varphi
_{0}$ is H\"{o}lder continuous on $\left[ 0,T\right] \times \mathcal{C}%
_{\beta ,s}$, provided that $\mathcal{C}_{\beta ,s}$ is endowed with the
metric topology of $V^{-\theta _{2}}\times W^{2s-1}$. Also, by virtue of
estimates (\ref{eq35}), (\ref{eq42}), $S_{\theta _{2}}\left( t\right) $
enjoys the smoothing property. Thus, using also the bound (\ref{higher-est})
and exploiting a well-known abstract result (see Appendix, Theorem \ref%
{abstract}), the existence of an exponential attractor with the properties
stated in Theorem \ref{expo-thm2}\ follows. The proof is finished.
\end{proof}

We finish this subsection with a simple corollary.

\begin{corollary}
Let the assumptions of Theorem \ref{expo-thm2} be satisfied. The global
attractor $\mathcal{A}_{0}$ has finite fractal dimension.
\end{corollary}

\begin{remark}
Note that Theorem \ref{expo-thm2} applies to the 3D NSV-AC system provided
that $s\in (0.5,0.75)$.
\end{remark}

\subsection{Convergence to steady states}

\label{s:convss}

In this subsection, we show that any global-in-time bounded solution to the
model (\ref{e:op}) converges to a single equilibrium as time tends to
infinity. The proof of the main result is based on a suitable version of the
Lojasiewicz--Simon theorem and the regularity results provided in Sections %
\ref{ss:attr-diss}, \ref{ss:attr-nondiss}. The question of such convergence
is usually a delicate matter since it is well known that the topology of the
set of stationary solutions of (\ref{e:op}) can be non-trivial even when $%
u\equiv 0$. In particular, there may be a continuum of stationary solutions
for (\ref{e:op}) even in the simplest cases, for instance when $\Omega $ is
a disk, see e.g., \cite{GG1, Ha} (cf. also \cite{Ab, Ab2, ZWH}). Since weak
solutions $\left( u\left( t\right) ,\phi \left( t\right) \right) $\ for (\ref%
{e:op}) become strong for times $t\geq t_{1}$ (for some $t_{1}=t_{1}\left(
\theta \right) >0$), we can confine ourselves to considering only strong
solutions.

We summarize the regularity results of the previous subsections in the
following proposition.

\begin{proposition}
\label{p:reg}Let the assumptions of Theorem \ref{t:attr-existb} when $\theta
>0$, and Theorem \ref{expo-thm2} when $\theta =0,$ be satisfied. Moreover,
assume that $g\in L_{tb}^{2}\left( \mathbb{R}_{+};V^{\beta -\theta }\right) $%
. For every $R>0$, there exists $C_{\ast }=C_{\ast }\left( R\right) >0$,
independent of time, such that,%
\begin{equation}
\sup_{\varphi _{0}\in \mathcal{B}_{\mathbb{V}_{\beta }}\left( R\right)
}\left\Vert S_{\theta _{2}}\left( t\right) \varphi _{0}\right\Vert _{\mathbb{%
V}_{\beta }}^{2}+\int_{t}^{t+1}\left( \left\Vert u\left( s\right)
\right\Vert _{\theta +\beta }^{2}+\left\Vert A_{1}^{3/2}\phi \left( s\right)
\right\Vert _{L^{2}}^{2}\right) ds\leq C_{\ast },\text{ for all }t\geq 0,
\label{reg4.1}
\end{equation}%
when $\theta >0$, and%
\begin{equation}
\sup_{\varphi _{0}\in \mathcal{B}_{\mathbb{V}_{\beta ,s}}\left( R\right)
}\left\Vert S_{\theta _{2}}\left( t\right) \varphi _{0}\right\Vert _{\mathbb{%
V}_{\beta ,s}}^{2}+\int_{t}^{t+1}\left( \left\Vert u\left( s\right)
\right\Vert _{\beta }^{2}+\left\Vert A_{1}^{\left( 2s+1\right) /2}\phi
\left( s\right) \right\Vert _{L^{2}}^{2}\right) ds\leq C_{\ast },\text{ for
all }t\geq 0,  \label{reg4.1bis}
\end{equation}%
when $\theta =0$. Here $\mathcal{B}_{\mathcal{X}}\left( R\right) $ denotes
the ball in $\mathcal{X}$ of radius $R,$ centered at $0$.
\end{proposition}

Note that $\mathbb{V}_{\beta ,1}=\mathbb{V}_{\beta }.$ For the sake of
simplicity and notation, below we will make the following convention: $%
\mathbb{V}_{\beta ,s}=\mathbb{V}_{\beta }$ when $\theta >0$ (and, thus,
always assume $s=1$), while for $\theta =0$, we recall that $s\in (\frac{1}{2%
},1]$ is possibly sufficiently small. Next, we characterize the structure of
the $\omega $-limit set for problem (\ref{e:op}) corresponding to any
initial datum $\varphi _{0}=\left( u_{0},\phi _{0}\right) \in \mathbb{V}%
_{\beta ,s}$. Recall that $\omega $-limit set is defined as follows:%
\begin{equation*}
\omega \left( \varphi _{0}\right) =\left\{ \left( u_{\ast },\phi _{\ast
}\right) \in \mathbb{V}_{\beta ,s}:\exists t_{n}\rightarrow \infty \text{
such that }\lim_{n\rightarrow \infty }\left\Vert \left( u\left( t_{n}\right)
,\phi \left( t_{n}\right) \right) -\left( u_{\ast },\phi _{\ast }\right)
\right\Vert _{\mathbb{V}_{\beta ,s}}=0\right\} .
\end{equation*}%
Clearly, $\omega \left( \varphi _{0}\right) $ is nonempty by virtue of
Proposition \ref{p:reg}.

\begin{lemma}
\label{ss:omega}Le the assumptions of Proposition \ref{p:reg} be satisfied,
and suppose that $g$ also obeys the following condition:%
\begin{equation}
\int_{t}^{\infty }\left\Vert g\left( s\right) \right\Vert _{-\theta -\theta
_{2}}^{2}ds\lesssim \left( 1+t\right) ^{-\left( 1+\delta \right) },\text{
for all }t\geq 0,  \label{g-assumpt}
\end{equation}%
for some constant $\delta >0.$ Then, the $\omega $-limit set $\omega \left(
\varphi _{0}\right) $ is a subset of%
\begin{equation*}
\mathcal{L}=\{\left( 0,\phi _{\ast }\right) :\phi _{\ast }\in D(A_{1}^{s}),%
\text{ }\phi _{\ast }\in \left[ -1,1\right] \text{ and (\ref{*}) holds}\},
\end{equation*}%
for some $s\in (\frac{1}{2},1]$, where (\ref{*}) is the following problem:%
\begin{equation}
A_{1}\phi _{\ast }+f\left( \phi _{\ast }\right) =0.  \label{*}
\end{equation}%
Moreover, we have%
\begin{equation}
\lim_{t\rightarrow \infty }\left\Vert u\left( t\right) \right\Vert _{-\theta
_{2}}=0,\text{ }\lim_{t\rightarrow \infty }\left\Vert A_{1}\phi \left(
t\right) +f\left( \phi \left( t\right) \right) \right\Vert _{L^{2}}=0.
\label{reg4.2}
\end{equation}
\end{lemma}

\begin{proof}
First, we have by Proposition \ref{p:reg}, the corresponding energy
inequality (\ref{e:exist-1}), and assumption (\ref{g-assumpt}) that%
\begin{equation}
\int_{0}^{\infty }y\left( s\right) ds<\infty ,\text{ }y:=\left\Vert
u\right\Vert _{-\theta _{2}}^{2}+\left\Vert \mu \right\Vert _{L^{2}}^{2},
\label{reg4.3}
\end{equation}%
where $\mu =A_{1}\phi +f\left( \phi \right) $. It follows from (\ref{reg4.1}%
)-(\ref{reg4.1bis}), (\ref{g-assumpt}) and a higher-order differential
inequality for the function $y\left( t\right) $:%
\begin{equation}
\frac{d}{dt}y\left( t\right) \leq C\left( 1+\left\Vert g\left( t\right)
\right\Vert _{-\theta -\theta _{2}}^{2}\right) ,  \label{reg4.3bis}
\end{equation}%
that%
\begin{equation*}
\lim_{t\rightarrow \infty }\left\Vert u\left( t\right) \right\Vert _{-\theta
_{2}}+\left\Vert \mu \left( t\right) \right\Vert _{L^{2}}=0.
\end{equation*}%
Hence, for any $\left( u_{\ast },\phi _{\ast }\right) \in \omega \left(
\varphi _{0}\right) $ we have $u_{\ast }\equiv 0.$ The assertion (\ref%
{reg4.2}) is also immediate. Moreover, by (\ref{reg4.1}), (\ref{reg4.1bis})
and (\ref{reg4.3}) it is easy to see that $\phi _{\ast }\in D(A_{1}^{s}),$ $%
\phi _{\ast }\in \left[ -1,1\right] $, and that the following inequality
holds:%
\begin{equation}
\lim_{n\rightarrow \infty }\left\Vert A_{1}^{s}\left( \phi \left(
t_{n}\right) -\phi _{\ast }\right) \right\Vert _{L^{2}}=0.  \label{reg4.4}
\end{equation}%
Let $\psi \in D\left( A_{1}^{1-s}\right) $. Then, $\phi _{\ast }$ satisfies (%
\ref{*}) on account of (\ref{reg4.2}), (\ref{reg4.4}), and the basic
inequality%
\begin{align*}
& \left\vert \left\langle A_{1}\phi _{\ast }+f\left( \phi _{\ast }\right)
,\psi \right\rangle \right\vert \\
& \leq \left( \left\Vert A_{1}^{s}\left( \phi \left( t_{n}\right) -\phi
_{\ast }\right) \right\Vert _{L^{2}}+\left\Vert f\left( \phi \left(
t_{n}\right) -f\left( \phi _{\ast }\right) \right) \right\Vert
_{L^{2}}+\left\Vert A_{1}^{s}\left( \phi \left( t_{n}\right) \right)
+f\left( \phi \left( t_{n}\right) \right) \right\Vert _{L^{2}}\right)
\left\Vert A_{1}^{1-s}\psi \right\Vert _{L^{2}},
\end{align*}%
Passing now to the limit as $t_{n}\rightarrow \infty $, the proof of Lemma %
\ref{ss:omega} is concluded. Finally, we only briefly sketch the details for
getting (\ref{reg4.3bis}) in the case $\theta >0$ (the case $\theta =0$ is
similar). We observe that for smoth solutions the function $y\left( t\right) 
$ satisfies the identity%
\begin{align}
& \frac{dy}{dt}+2\left\langle A_{0}u,Nu\right\rangle +\left\Vert
A_{1}^{1/2}\mu \right\Vert _{L^{2}}^{2}  \label{reg4.3es1} \\
& =2\left\langle g,Nu\right\rangle +2\left\langle R_{0}\left( \mu ,\phi
\right) ,Nu\right\rangle -\left\langle A_{1}^{1/2}B_{1}\left( u,\phi \right)
,A_{1}^{1/2}\mu \right\rangle  \notag \\
& -\left\langle f^{^{\prime }}\left( \phi \right) B_{1}\left( u,\phi \right)
,\mu \right\rangle -\left\langle f^{^{\prime }}\left( \phi \right) \mu ,\mu
\right\rangle .  \notag
\end{align}%
We can bound the nonlinear terms on the right-hand side of (\ref{reg4.3es1}%
), by using the following facts:

\begin{itemize}
\item $R_{0}\left( \mu ,\phi \right) \in L^{2}\left( 0,\infty ;V^{-\theta
-\theta _{2}}\right) ,$ as a mapping from $L^{2}\times W^{1}\rightarrow
V^{-\theta -\theta _{2}}$ due to (\ref{reg4.3}) and (\ref{reg4.1}), since $%
\theta +\theta _{2}\geq 1$ (see Lemma \ref{l:hole}). It follows that%
\begin{equation*}
\left\vert \left\langle R_{0}\left( \mu ,\phi \right) ,Nu\right\rangle
\right\vert \leq \delta \left\Vert u\right\Vert _{\theta -\theta
_{2}}^{2}+C_{\delta }\left\Vert \mu \right\Vert _{L^{2}}^{2}\left\Vert \phi
\right\Vert _{2}^{2}.
\end{equation*}

\item $B_{1}\left( u,\phi \right) \in L^{2}\left( 0,\infty ;L^{2}\left(
\Omega \right) \right) ,$ as a mapping from $V^{\theta +\theta _{2}}\times
W^{1}\rightarrow L^{2}\left( \Omega \right) ,$ due to (\ref{reg4.3}) and (%
\ref{reg4.1}). Since $\phi $ is also bounded, we have%
\begin{equation*}
\left\vert \left\langle f^{^{\prime }}\left( \phi \right) B_{1}\left( u,\phi
\right) ,\mu \right\rangle \right\vert \leq \delta \left\Vert u\right\Vert
_{\theta -\theta _{2}}^{2}+C_{\delta }\left\Vert f^{^{\prime }}\left( \phi
\right) \right\Vert _{L^{\infty }}^{2}\left\Vert \phi \right\Vert
_{2}^{2}\left\Vert \mu \right\Vert _{L^{2}}^{2}.
\end{equation*}

\item For the third term on the right-hand side of (\ref{reg4.3es1}), we
employ the same strategy of proof used in Theorem \ref{t:reg}, (\ref{est14}%
)-(\ref{est15}) in the case $n=3$, and (\ref{est14bis})-(\ref{est15q}) for $%
n=2.$ For instance, \ when $n=3$ we derive%
\begin{align*}
\left\vert \left\langle A_{1}^{1/2}B_{1}\left( u,\phi \right)
,A_{1}^{1/2}\mu \right\rangle \right\vert & \lesssim \left( \left\Vert
Nu\right\Vert _{L^{6}}\left\Vert \nabla \phi \right\Vert
_{W^{1,3}}+\left\Vert Nu\right\Vert _{1}\left\Vert \nabla \phi \right\Vert
_{L^{\infty }}\right) \left\Vert A_{1}^{1/2}\mu \right\Vert _{L^{2}} \\
& =I_{1}+I_{2},
\end{align*}%
where%
\begin{equation*}
\begin{array}{l}
I_{1}\leq 2\delta \left\Vert A_{1}^{1/2}\mu \right\Vert
_{L^{2}}^{2}+C_{\delta }\left( \left\Vert u\right\Vert _{-\theta
_{2}}^{8}\left\Vert \phi \right\Vert _{1}^{2}+\left\Vert u\right\Vert
_{-\theta _{2}}^{2}\left\Vert \phi \right\Vert _{1}^{3/2}\left\Vert
f^{^{\prime }}\left( \phi \right) \right\Vert _{L^{\infty }}^{2}\right) , \\ 
I_{2}\leq 2\delta \left\Vert A_{1}^{1/2}\mu \right\Vert
_{L^{2}}^{2}+C_{\delta }\left( \left\Vert u\right\Vert _{-\theta
_{2}}^{4}\left\Vert \phi \right\Vert _{2}^{2}+\left\Vert u\right\Vert
_{-\theta _{2}}^{2}\left\Vert \phi \right\Vert _{2}\left\Vert f^{^{\prime
}}\left( \phi \right) \right\Vert _{L^{\infty }}\left\Vert \phi \right\Vert
_{1}\right) .%
\end{array}%
\end{equation*}

\item The bound on the final term is basic%
\begin{equation*}
\left\vert \left\langle f^{^{\prime }}\left( \phi \right) \mu ,\mu
\right\rangle \right\vert \leq \left\Vert f^{^{\prime }}\left( \phi \right)
\right\Vert _{L^{\infty }}\left\Vert \mu \right\Vert _{L^{2}}^{2}.
\end{equation*}
\end{itemize}

We can now choose $\delta \sim \min \left( c_{A_{0}},1\right) >0$
sufficiently small in all these estimates, and use the coercitivity of the
operator $A_{0}$ together with (\ref{e:exist-1}) to handle the term $%
\left\langle g,Nu\right\rangle $. Then recalling that for a strong solution, 
$y\left( t\right) \leq C$ uniformly, for all $t\geq 0$ and reporting all the
preceeding estimates in (\ref{reg4.3es1}) we can easily obtain the claim (%
\ref{reg4.2}) on account of the application of \cite[Lemma 6.2.1]{Zh}.
\end{proof}

Consequently, for the energy functional%
\begin{equation*}
\mathcal{E}\left( u,\phi \right) :=\frac{1}{2}\left\langle u,Nu\right\rangle
+\widehat{\mathcal{E}}\left( \phi \right) ,
\end{equation*}%
where%
\begin{equation*}
\widehat{\mathcal{E}}\left( \phi \right) :=\frac{1}{2}||A_{1}^{1/2}\phi
||_{L^{2}}^{2}+\int_{\Omega }F\left( \phi \right) dx,
\end{equation*}%
the following statement holds. Note that, by the above Lemma \ref{ss:omega}, 
$\phi _{\ast }$ is a critical point of $\widehat{\mathcal{E}}$ over $%
D(A_{1}^{1/2})\cap L^{\infty }\left( \Omega \right) .$

\begin{proposition}
\label{p:const-energ}There exists a constant $e_{\infty }\in \mathbb{R}$
such that $\widehat{\mathcal{E}}\left( \phi _{\ast }\right) =e_{\infty },$
for all $\left( 0,\phi _{\ast }\right) \in \mathcal{L}$, and we have%
\begin{equation}
\lim_{t\rightarrow \infty }\mathcal{E}\left( u\left( t\right) ,\phi \left(
t\right) \right) =e_{\infty }.  \label{e_const}
\end{equation}%
Moreover, the functional $\Phi \left( t\right) $ is decreasing along all
strong trajectories $\left( u\left( t\right) ,\phi \left( t\right) \right) $
and, for all $t\geq 0,$%
\begin{equation*}
\frac{d}{dt}\Phi \left( t\right) \leq -\left( \frac{c_{A_{0}}}{2}\left\Vert
u\left( t\right) \right\Vert _{\theta -\theta _{2}}^{2}+\left\Vert A_{1}\phi
\left( t\right) +f\left( \phi \left( t\right) \right) \right\Vert
_{L^{2}}^{2}\right) ,
\end{equation*}%
where%
\begin{equation}
\Phi \left( t\right) :=\mathcal{E}\left( u\left( t\right) ,\phi \left(
t\right) \right) +\frac{\left\Vert N\right\Vert _{-\theta _{2};\theta
_{2}}^{2}}{2c_{A_{0}}}\int_{t}^{\infty }\left\Vert g\left( s\right)
\right\Vert _{-\theta -\theta _{2}}^{2}ds.  \label{big_phi}
\end{equation}
\end{proposition}

Even though we are dealing with an asymptotically decaying force due to (\ref%
{g-assumpt}), we cannot conclude that each strong solution of (\ref{e:op})
converges to a \emph{single} equilibrium, for $\mathcal{L}$ can be a
continuum (see, e.g., \cite{Ha}). However, we can establish this fact when
the nonlinear function $f$ is real analytic. The version of the \L %
ojasiewicz-Simon inequality we need is given by the following lemma (see 
\cite{CJ, JB}).

\begin{lemma}
\label{l3} For the above setting, let $f$ be real analytic. There exist
constants $\zeta \in (0,1/2)$ and $C_{L}>0,$ $\eta >0$ depending on $\left(
0,\phi _{\ast }\right) $ such that, for any $\phi \in D\left(
A_{1}^{s}\right) ,$ if $\left\Vert \left( \phi -\phi _{\ast }\right)
\right\Vert _{1}\leq \eta ,$ denoting by $\widehat{\mathcal{E}}^{\prime }$
the Fr\'{e}chet derivative of $\widehat{\mathcal{E}}$, we have%
\begin{equation}
C_{L}||\widehat{\mathcal{E}}^{^{\prime }}\left( \phi \left( t\right) \right)
||_{-1}\geq |\widehat{\mathcal{E}}\left( \phi \left( t\right) \right) -%
\widehat{\mathcal{E}}\left( \phi _{\ast }\right) |^{1-\zeta }.
\label{3.2bis}
\end{equation}
\end{lemma}

The result below is concerned with the convergence of trajectories of
problem (\ref{e:op}) to single equilibria, which shows, in a strong form,
their (global)\ asymptotic stability.

\begin{theorem}
\label{conv-equil} Let the assumptions of Lemma \ref{ss:omega} hold. In
addition, assume that $f$ is real analytic. For any given initial datum $%
\varphi _{0}=\left( u_{0},\phi _{0}\right) \in \mathbb{V}_{\beta ,s}$, the
corresponding solution $\varphi \left( t\right) =\left( u\left( t\right)
,\phi \left( t\right) \right) =S_{\theta _{2}}\varphi _{0}$ to problem (\ref%
{e:op})\ converges to a \emph{single} equilibrium $\left( 0,\phi _{\ast
}\right) $\ in the strong topology of $\mathbb{V}_{\beta ,s},$ as time goes
to infinity. Moreover, there exists $\xi \in \left( 0,1\right) ,$ depending
on $\phi _{\ast }$\ and the other physical parameters of the problem, such
that 
\begin{equation}
\left\Vert u\left( t\right) \right\Vert _{-\theta _{2}}+||A_{1}^{1/2}\left(
\phi (t)-\phi _{\ast }\right) ||_{L^{2}}\lesssim (1+t)^{-\xi },
\label{conv_rate}
\end{equation}%
for all $t\geq 0.$
\end{theorem}

\begin{proof}
We prove the case $\theta >0$ (the case $\theta =0$ is analogous and follows
with minor modifications). We adapt the ideas of \cite{CJ, JB, GG1} (cf.
also \cite{Ab, GG2, ZWH}) to prove the claim. On account of the first
statement of (\ref{reg4.2}), it suffices to prove the claim only for the
phase-field component $\phi $. First, by Lemma \ref{ss:omega} the omega
limit set $\omega \left( \varphi _{0}\right) $ is a non-empty and compact
subset of $\mathcal{Y}_{\theta _{2}}$.$\ $Secondly, we can choose a
sufficiently large time $t_{1}>0$ such that for all $t\geq t_{1},$ we have $%
||A_{1}^{1/2}\left( \phi \left( t\right) -\phi _{\ast }\right)
||_{L^{2}}<\eta ,$ such that the conclusion of Lemma \ref{l3} holds with $%
\widehat{\mathcal{E}}^{^{\prime }}\left( \phi \right) =A_{1}\phi +f\left(
\phi \right) $ provided that we choose even a smaller constant $\zeta \in (0,%
\frac{1}{2})\cap $ $(0,\delta \left( 1+\delta \right) ^{-1})$. Without loss
of generality, in what follows we let $t_{1}=1$ and assume $\delta <1$.
Next, define%
\begin{equation*}
\Sigma :=\{t\geq 1:||A_{1}^{1/2}\left( \phi \left( t\right) -\phi _{\ast
}\right) ||_{L^{2}}\leq \eta /3\}
\end{equation*}%
and observe that $\Sigma $ is unbounded by Lemma \ref{ss:omega}. For every $%
t\in \Sigma ,$ we define%
\begin{equation*}
\tau \left( t\right) =\sup \{t^{^{\prime }}\geq t:\sup_{s\in \lbrack
t,t^{^{\prime }}]}||A_{1}^{1/2}\left( \phi \left( t\right) -\phi _{\ast
}\right) ||_{L^{2}}<\eta \}.
\end{equation*}%
By continuity, $\tau \left( t\right) >t$ for every $t\in \Sigma $. Let now $%
t_{0}\in \Sigma $ and divide the interval $J:=[t_{0},\tau \left(
t_{0}\right) )$ into two subsets%
\begin{equation*}
\Sigma _{1}:=\left\{ t\in J:\Upsilon \left( t\right) \geq \left(
\int_{t}^{\tau \left( t_{0}\right) }\left\Vert g\left( s\right) \right\Vert
_{-\theta -\theta _{2}}^{2}ds\right) ^{1-\zeta }\right\} ,\text{ }\Sigma
_{2}:=J\backslash \Sigma _{1},
\end{equation*}%
where%
\begin{equation*}
\Upsilon \left( t\right) :=\left\Vert u\left( t\right) \right\Vert _{\theta
-\theta _{2}}+\left\Vert A_{1}\phi \left( t\right) +f\left( \phi \left(
t\right) \right) \right\Vert _{L^{2}}.
\end{equation*}

Define further%
\begin{equation*}
\widehat{\Phi }\left( t\right) :=\mathcal{E}\left( u\left( t\right) ,\phi
\left( t\right) \right) -\widehat{\mathcal{E}}\left( \phi _{\ast }\right) +%
\frac{\left\Vert N\right\Vert _{-\theta _{2};\theta _{2}}^{2}}{2c_{A_{0}}}%
\int_{t}^{\tau \left( t_{0}\right) }\left\Vert g\left( s\right) \right\Vert
_{-\theta -\theta _{2}}^{2}ds
\end{equation*}%
and notice that $\widehat{\Phi }\left( t\right) $ differs from $\Phi \left(
t\right) $ in (\ref{big_phi}) only by a constant. Hence, for every $t\in J$
we have%
\begin{equation*}
\frac{d}{dt}\widehat{\Phi }\left( t\right) \lesssim -\Upsilon ^{2}\left(
t\right) \leq 0
\end{equation*}%
so that $\widehat{\Phi }$ is a decreasing function. Moreover, for every $%
t\in J$ we have%
\begin{align}
\frac{d}{dt}\left( |\widehat{\Phi }\left( t\right) |^{\zeta }sgn(\widehat{%
\Phi }\left( t\right) \right) & =\zeta |\widehat{\Phi }\left( t\right)
|^{\zeta -1}\frac{d}{dt}\widehat{\Phi }\left( t\right)  \label{rel4.4} \\
& \lesssim -|\widehat{\Phi }\left( t\right) |^{\zeta -1}\Upsilon ^{2}\left(
t\right) ,  \notag
\end{align}%
which implies that the functional $sgn(\widehat{\Phi }\left( t\right) )|%
\widehat{\Phi }\left( t\right) |^{\zeta }$ is decreasing on $J$. By (\ref%
{3.2bis}) and Proposition \ref{p:const-energ} (indeed, $\left\Vert
u\right\Vert _{-\theta _{2}}^{2\left( 1-\zeta \right) }\lesssim \left\Vert
u\right\Vert _{\theta -\theta _{2}}$ since $\zeta <\frac{1}{2}$), for every $%
t\in \Sigma _{1}$ we have%
\begin{align}
|\widehat{\Phi }\left( t\right) |^{1-\zeta }& \leq \left\vert \mathcal{E}%
\left( u\left( t\right) ,\phi \left( t\right) \right) -\widehat{\mathcal{E}}%
\left( \phi _{\ast }\right) \right\vert ^{1-\zeta }+\left( \frac{\left\Vert
N\right\Vert _{-\theta _{2};\theta _{2}}^{2}}{2c_{A_{0}}}\int_{t}^{\tau
\left( t_{0}\right) }\left\Vert g\left( s\right) \right\Vert _{-\theta
-\theta _{2}}^{2}ds\right) ^{1-\zeta }  \label{rel4.5} \\
& \lesssim \Upsilon \left( t\right) ,  \notag
\end{align}%
which together with equation (\ref{rel4.4}) yields%
\begin{equation}
-\frac{d}{dt}\left( |\widehat{\Phi }\left( t\right) |^{\zeta }sgn(\widehat{%
\Phi }\left( t\right) \right) \gtrsim \Upsilon \left( t\right) .
\label{rel4.6}
\end{equation}%
Moreover, exploiting (\ref{rel4.6}) we have%
\begin{align}
\int_{\Sigma _{1}}\Upsilon \left( s\right) ds& \lesssim -\int_{\Sigma _{1}}%
\frac{d}{ds}\left( |\widehat{\Phi }\left( s\right) |^{\zeta }sgn(\widehat{%
\Phi }\left( s\right) \right) ds  \label{rel4.7} \\
& \lesssim \left( |\widehat{\Phi }\left( t_{0}\right) |^{\zeta }+|\widehat{%
\Phi }\left( \tau \left( t_{0}\right) \right) |^{\zeta }\right) ,  \notag
\end{align}%
where we interpret the term involving $\tau \left( t_{0}\right) $ on the
right hand side of (\ref{rel4.7}) as $0$ if $\tau \left( t_{0}\right)
=\infty $ (recall (\ref{e_const})). On the other hand, if $t\in \Sigma _{2},$
using assumption (\ref{g-assumpt}) we obtain%
\begin{equation}
\Upsilon \left( t\right) \leq \left( \int_{t}^{\tau \left( t_{0}\right)
}\left\Vert g\left( s\right) \right\Vert _{-\theta -\theta
_{2}}^{2}ds\right) ^{1-\zeta }\lesssim \left( 1+t\right) ^{-\left( 1-\zeta
\right) \left( 1+\delta \right) },  \label{rel4.8}
\end{equation}%
so once again the function $\Upsilon $ is dominated by an integrable
function on $\Sigma _{2}$ since $\zeta \left( 1+\delta \right) <\delta $.
Combining the inequalities (\ref{rel4.7}), (\ref{rel4.8}), we deduce that $%
\Upsilon $ is absolutely integrable on $J$ and%
\begin{equation}
\lim_{t_{0}\rightarrow \infty ,t_{0}\in \Sigma }\int_{t_{0}}^{\tau \left(
t_{0}\right) }\Upsilon \left( s\right) ds=0.  \label{rel4.9}
\end{equation}%
On the other hand, recalling estimates (\ref{reg4.1}), (\ref{reg4.3}), from
the second equation of (\ref{e:op}) it follows that%
\begin{align}
\left\Vert \partial _{t}\phi \left( t\right) \right\Vert _{L^{2}}& \leq
\left\Vert B_{1}\left( u\left( t\right) ,\phi \left( t\right) \right)
\right\Vert _{L^{2}}+\left\Vert A_{1}\phi \left( t\right) +f\left( \phi
\left( t\right) \right) \right\Vert _{L^{2}}  \label{rel4.10} \\
& \leq \left\Vert u\left( t\right) \right\Vert _{\theta -\theta
_{2}}\left\Vert A\phi \left( t\right) \right\Vert _{L^{2}}+\left\Vert
A_{1}\phi \left( t\right) +f\left( \phi \left( t\right) \right) \right\Vert
_{L^{2}}  \notag \\
& \lesssim \Upsilon \left( s\right) ,  \notag
\end{align}%
since $N:V^{-2\theta _{2}}\rightarrow V^{0}$ is bounded, and $V^{-\theta
_{2}}\subseteq V^{-2\theta _{2}},$ $\theta _{2},\theta \geq 0.$
Consequently, we also have%
\begin{equation}
\lim_{t_{0}\rightarrow \infty ,t_{0}\in \Sigma }\int_{t_{0}}^{\tau \left(
t_{0}\right) }\left\Vert \partial _{t}\phi \left( s\right) \right\Vert
_{L^{2}}ds=0.  \label{rel4.11}
\end{equation}%
For $t\in \Sigma ,$ note that the inequality%
\begin{equation*}
\left\Vert \phi \left( t\right) -\phi _{\ast }\right\Vert _{L^{2}}\leq
\int_{t}^{t_{0}}\left\Vert \partial _{t}\phi \left( s\right) \right\Vert
_{L^{2}}ds+\left\Vert \phi \left( t_{0}\right) -\phi _{\ast }\right\Vert
_{L^{2}}
\end{equation*}%
implies that $\tau \left( t_{0}\right) =\infty ,$ for some $t_{0}\in \Sigma $%
. Indeed, let us assume for a second that the latter statement is not true.
Then, by definition of $\tau \left( t_{0}\right) $, $||A_{1}^{1/2}\left(
\phi \left( t\right) -\phi _{\ast }\right) ||_{L^{2}}=\eta $ for every $%
t_{0}\in \Sigma $. Let now $\left\{ t_{n}\right\} \subset \Sigma $ be an
unbounded sequence such that%
\begin{equation*}
\lim_{t_{n}\rightarrow \infty }||A_{1}^{1/2}\left( \phi \left( t_{n}\right)
-\phi _{\ast }\right) ||_{L^{2}}=0.
\end{equation*}%
By compactness and Lemma \ref{ss:omega}, we can now pass to a subsequence of 
$\left\{ t_{n}\right\} $ if necessary to conclude that one can find $%
\widetilde{\phi }\in \omega \left( \varphi _{0}\right) $ such that $||%
\widetilde{\phi }-\phi ||_{1}=\eta >0$ and $\lim_{t_{n}\rightarrow \infty }||%
\widetilde{\phi }-\phi \left( \tau \left( t_{n}\right) \right) ||_{1}=0.$
Then, the above inequality gives%
\begin{equation*}
0<||\widetilde{\phi }-\phi ||_{L^{2}}\leq \lim_{t_{n}\rightarrow \infty
}\{\int_{t_{n}}^{\tau \left( t_{n}\right) }\left\Vert \partial _{t}\phi
\left( s\right) \right\Vert _{L^{2}}ds+\left\Vert \phi \left( t_{n}\right)
-\phi \right\Vert _{L^{2}}\}=0,
\end{equation*}%
which is a contradiction. Hence, we conclude that we must have $\tau \left(
t_{0}\right) =\infty $, for some $t_{0}$ sufficiently large. Thus, the above
arguments (see, in particular, (\ref{rel4.11})) imply that $\left\Vert
\partial _{t}\phi \right\Vert _{L^{2}}$ is absolutely integrable on $%
[t_{0},\infty ),$ which implies that the limit of $\phi \left( t\right) $
exists, as time goes to infinity. By compactness and a basic interpolation
inequality, we have $\phi \left( t\right) \rightarrow \phi _{\ast }$ in the
strong topology of $D\left( A_{1}\right) .$ Hence, $\omega \left( \varphi
_{0}\right) =\left\{ \left( 0,\phi _{\ast }\right) \right\} ,$ as claimed.
The estimate of the rate of convergence in (\ref{conv_rate}) is a
straight-forward consequence of (\ref{rel4.4})-(\ref{rel4.5}) and the
definition of $\Phi $ (see, e.g., \cite[Theorem 5.7]{GG1}). We leave the
details to the interested reader. The proof of Theorem \ref{conv-equil} is
complete.
\end{proof}

\begin{remark}
\label{sing-rem}All the results of the previous sections and subsections can
be extended for singular potentials of of the form (\ref{sing-pot}). Let us
suppose that $F\in C^{2}(-1,1)$, set $f=F^{\prime }$ and assume%
\begin{equation}
\lim_{r\rightarrow \pm 1}f\left( r\right) =\pm \infty \;\text{ and }%
\lim_{r\rightarrow \pm 1}f^{\prime }\left( r\right) =+\infty .  \label{6.2}
\end{equation}%
It is easy to see that the derivative of $F$ defined by (\ref{sing-pot})
satisfies (\ref{6.2}). Although the potential $f$ is singular, we can still
use the results derived in the previous sections, since the solutions to our
problem are smooth and strictly separated thanks to \cite[Theorem 6.1]{GG1},
if additionally we assume $\left( u_{0},\phi _{0}\right) \in \mathcal{Y}%
_{\theta _{2}}^{\delta }$, for some $\delta \in \left( 0,1\right) .$ Here,
we have defined%
\begin{equation*}
\mathcal{Y}_{\theta _{2}}^{\delta }:=\left\{ \left( u_{0},\phi _{0}\right)
\in V^{-\theta _{2}}\times (D(A_{1}^{1/2})\cap L^{\infty }\left( \Omega
\right) ):\left\Vert \phi _{0}\right\Vert _{L^{\infty }}\leq \delta
<1\right\} .
\end{equation*}%
Consequently, for $\phi \in \mathcal{Y}_{\theta _{2}}^{\delta }$, $f(\phi )$
and any of its higher-order derivatives are bounded provided that we replace 
$\mathcal{Y}_{\theta _{2}}$ by $\mathcal{Y}_{\theta _{2}}^{\delta }$
everywhere in the paper. Thus, the arguments used in the previous
(sub)sections are still valid in the present case.
\end{remark}

\section{Remarks on a regularized family for the NSE and MHD models}

\label{rem-mhd}

As in Section \ref{s:prelim}, consider the following system%
\begin{equation}
\left\{ 
\begin{array}{l}
\partial _{t}u+A_{0}u+B_{0}(u,u)=g, \\ 
u\left( 0\right) =u_{0},%
\end{array}%
\right.  \label{MHD}
\end{equation}%
on the time interval $[0,T]$. Bearing in mind the model (\ref{MHD}), we are
mainly interested in bilinear maps $B_{0}$ of the form (\ref{e:b-def}). We
recall that the formulation (\ref{MHD}) here includes not only various
regularized models for the Navier-Stokes (NSE) equation but also certain
(regularized or un-regularized) magnetohydrodynamics (MHD) models (see \cite%
{HLT}).

In this section, we show how to close a gap in the proof of \cite[Section 5,
Corollary 5.4]{HLT} whose assumptions can only be verified in the case $%
\theta >0$. Note that when $\theta =0$, the assumptions of \cite[Theorem
5.1, (b)]{HLT} do not longer provide the existence of a compact absorbing
set as claimed on \cite[pg. 550]{HLT} (see also Remark \ref{attr-rem}).
Besides, it is well-known that in the non-dissipative case when $\theta =0,$
the regularized model for the Navier-Stokes equation looses its parabolic
character, see Section \ref{ss:attr-nondiss} for related discussions.

Our analogue result of \cite[Corollary 5.4]{HLT}\ in the case $\theta =0$ is
as follows.

\begin{theorem}
\label{smooth-mhd}Let the following conditions hold and recall that $\theta
_{1},\theta _{2}\in \mathbb{R}.$

(i) $b_{0}(w,v,Nv)=0,$ for any $v,w\in V^{-\theta _{2}}$;

(ii) $b_{0}:V^{\bar{\sigma}_{1}}\times V^{\bar{\sigma}_{2}}\times V^{\bar{%
\gamma}}\rightarrow \mathbb{R}$ is bounded for some $\bar{\sigma}%
_{i}<-\theta _{2}$, $i=1,2$, and $\bar{\gamma}\in \mathbb{R}$;

(iii) $b_{0}:V^{\sigma _{1}}\times V^{-\theta _{2}}\times V^{\sigma
_{2}}\rightarrow \mathbb{R}$ is bounded for some $\sigma _{1}\leq -\theta
_{2}$ and $\sigma _{2}\leq \theta _{2}$ with $\sigma _{1}+\sigma _{2}\leq 0$;

(iv) $\left\langle A_{0}v,Nv\right\rangle \geq c_{A_{0}}\Vert v\Vert
_{-\theta _{2}}^{2}$ for any $v\in V^{-\theta _{2}}$, with a constant $%
c_{A_{0}}>0$;

In addition, for some $\beta >-\theta _{2},$ let the following conditions
hold:

(v) $\left\langle A_{0}v,\left( I-\Delta \right) ^{\beta }v\right\rangle
\geq c_{A_{0}}\Vert v\Vert _{\beta }^{2},$ for any $v\in V^{\beta },$ for
some $c_{A_{0}}>0$;

(vi) $b_{0}:V^{-\theta _{2}}\times V^{\beta }\times V^{-\beta }\rightarrow 
\mathbb{R}$ is bounded;

(vii) $g\in V^{\beta }$ is time independent.

Then, there exists a compact attractor $\mathcal{A}_{reg}\subset V^{-\theta
_{2}}$ for the system (\ref{MHD}) which attracts the bounded sets of $%
V^{-\theta _{2}}$. Moreover, $\mathcal{A}_{reg}$ is connected and it is the
maximal bounded attractor in $V^{-\theta _{2}}.$
\end{theorem}

\begin{remark}
(\textbf{a}) For instance, the 3D Navier-Stokes-Voight (NSV) equation
satisfies all the requirements (i)-(vii) of Theorem~\ref{smooth-mhd}. Note
that the existence of the global attractor and several properties for the 3D
NSV were derived in \cite{KT1, KT2}. The assumptions (i)-(iv) of the above
theorem insure that there is a continuous (nonlinear)\ semigroup for problem
(\ref{MHD}),%
\begin{equation}
\begin{array}{ll}
S_{NS}(t):V^{-\theta _{2}}\rightarrow V^{-\theta _{2}},\text{ }t\geq 0, & 
\\ 
u_{0}\mapsto u\left( t\right) . & 
\end{array}%
\end{equation}

(\textbf{b}) The final assumptions (v)-(vii) are used to derive that $S_{NS}$
is asymptotically smooth, as in the proof of Theorem \ref{smooth}. Indeed,
actual modifications of the proof of Theorem \ref{smooth} are not necessary
(i.e., we can let $\phi \equiv 0$ in (\ref{e:op}) in order to completely
decouple the system); we use a suitable decomposition of the velocity $%
u=u_{d}+u_{c},$ such that \thinspace $u_{d}\left( t\right) =S_{NS}^{d}\left(
t\right) u_{0}$ solves%
\begin{equation}
\left\{ 
\begin{array}{l}
\partial _{t}u_{d}+A_{0}u_{d}+B_{0}(u,u_{d})=0, \\ 
u_{d}\left( 0\right) =u_{0},%
\end{array}%
\right.  \label{NS-decay}
\end{equation}%
and $u_{c}\left( t\right) =S_{NS}^{c}\left( t\right) u_{0}$ solves%
\begin{equation}
\left\{ 
\begin{array}{l}
\partial _{t}u_{c}+A_{0}u_{c}+B_{0}(u,u_{c})=g, \\ 
u_{c}\left( 0\right) =0.%
\end{array}%
\right.  \label{NS-comp}
\end{equation}
\end{remark}

We can now obtain the following result which is also new in the literature.
For instance, it allows us to derive the existence of an exponential
attractor for the 3D NSV model, which was not reported anywhere else. Recall
again that $\theta =0$.

\begin{theorem}
\label{expo-att-mhd}Let the assumptions of Theorem \ref{smooth-mhd},
(i)-(v), (vii) be satisfied, and instead of (vi), assume that $%
b_{0}:V^{-\theta _{2}}\times V^{-\theta _{2}}\times V^{-\beta }\rightarrow 
\mathbb{R}$ is bounded. Then $\left( S_{NS},V^{-\theta _{2}}\right) $
possesses an exponential attractor $\mathcal{M}_{reg}\subset V^{\beta }$
which is bounded in $V^{\beta }$. Thus, by definition, we have

(a) $\mathcal{M}_{reg}$ is compact and semi-invariant with respect $%
S_{NS}\left( t\right) ,$ that is,%
\begin{equation*}
S_{NS}\left( t\right) \left( \mathcal{M}_{reg}\right) \subseteq \mathcal{M}%
_{reg},\quad \forall \,t\geq 0.
\end{equation*}

(b) The fractal dimension $\dim _{F}\left( \mathcal{M}_{reg},V^{-\theta
_{2}}\right) $ of $\mathcal{M}_{reg}$ is finite and an upper bound can be
computed explicitly.

(c) $\mathcal{M}_{reg}$ attracts exponentially fast any bounded subset $B$
of $V^{-\theta _{2}}$, that is, there exist a positive nondecreasing
function $Q$ and a constant $\rho >0$ such that 
\begin{equation*}
dist_{V^{-\theta _{2}}}\left( S_{NS}\left( t\right) B,\mathcal{M}%
_{reg}\right) \leq Q(\Vert B\Vert _{V^{-\theta _{2}}})e^{-\rho t},\quad
\forall \,t\geq 0.
\end{equation*}%
Both $Q$ and $\rho $ can be explicitly calculated.
\end{theorem}

\begin{proof}
The proof requires only minor modifications from that of Theorem \ref%
{expo-thm2}. We briefly sketch the main two (crucial)\ steps, the existence
of a compact exponentially attracting set in $V^{-\theta _{2}},$ and the
smoothing property for $\left( S_{NS},V^{-\theta _{2}}\right) $.

\noindent \emph{Step (i)}. We claim that there exists a bounded subset of $%
V^{\beta }$ which is positively invariant for $S_{NS}(t)$ and attracts any
bounded set $B$\ of $V^{-\theta _{2}}$ exponentially fast. We work with the
semigroup decomposition (\ref{NS-decay})-(\ref{NS-comp}). Fist it is clear
that $\left\Vert u_{d}\left( t\right) \right\Vert _{-\theta _{2}}^{2}\leq
e^{-c_{A_{0}}t}\left\Vert u_{0}\right\Vert _{-\theta _{2}}^{2}$, for all $%
t\geq 0$. In order to show the bound for the $u_{c}$-component, we pair the
first equation of (\ref{NS-comp}) with $\Lambda ^{2\beta }u_{c}$, and then
use the fact that $b_{0}:V^{-\theta _{2}}\times V^{-\theta _{2}}\times
V^{-\beta }\rightarrow \mathbb{R}$ is bounded. After standard
transformations, we finally have the following inequality:%
\begin{equation}
\frac{d}{dt}\left\Vert u_{c}\right\Vert _{\beta }^{2}+c_{A_{0}}\left\Vert
u_{c}\right\Vert _{\beta }^{2}\lesssim \left\Vert u\right\Vert _{-\theta
_{2}}^{2}\left( \left\Vert u_{d}\right\Vert _{-\theta _{2}}^{2}+\left\Vert
u\right\Vert _{-\theta _{2}}^{2}\right) +\left\Vert g\right\Vert _{\beta
}^{2}\leq C,  \label{5.2}
\end{equation}%
for some positive constant $C$ which depends only on $\mathcal{B}_{\theta
_{2}},$ $c_{A_{0}}$ and $g$. Hence, noting once again that $u_{c}\left(
0\right) =0$, application of Lemma \ref{Gineq} (Appendix) yields%
\begin{equation}
\left\Vert u_{c}\left( t\right) \right\Vert _{\beta }^{2}\leq R_{\beta
},\quad \text{\ for all }t\geq 0,  \label{5.3}
\end{equation}%
for some constant $R_{\beta }>0$ which is independent of time and initial
data. Let now%
\begin{equation*}
\mathcal{D}_{\beta }:=\left\{ u\in V^{\beta }:\Vert u\Vert _{\beta }\leq
R_{\beta }\right\} .
\end{equation*}%
On the other hand, recalling that $\left\Vert S_{NS}^{d}\left( t\right)
u_{0}\right\Vert _{-\theta _{2}}\lesssim e^{-\rho t}$, for all $t\geq 0$, (%
\ref{5.3}) implies that $\mathcal{D}_{\beta }$ attracts $\mathcal{B}_{\theta
_{2}}$exponentially fast, that is,%
\begin{equation}
dist_{V^{-\theta _{2}}}\left( S_{NS}\left( t\right) \mathcal{B}_{\theta
_{2}},\mathcal{D}_{\beta }\right) \leq \left\Vert S_{NS}\left( t\right)
u_{0}-u_{c}\left( t\right) \right\Vert \lesssim e^{-\rho t},\text{ \ for all 
}t\geq 0.  \label{5.4}
\end{equation}%
Since by \cite[Theorem 5.1]{HLT} we already know that for every nonempty
bounded subset $B$ of $V^{-\theta _{2}}$, $dist_{V^{-\theta
_{2}}}(S_{NS}(t)B,\mathcal{B}_{\theta _{2}})\lesssim e^{-\zeta t}$, for all $%
t\geq 0$, as usual we can appeal to (\ref{5.4}) and the transitivity
property of the exponential attraction \cite[Theorem 5.1]{FGMZ04} to infer%
\begin{equation}
dist_{V^{-\theta _{2}}}(S_{NS}(t)B,\mathcal{D}_{\beta })\lesssim e^{-\kappa
t},\text{ for all }t\geq 0,  \label{5.5}
\end{equation}%
for some $\kappa >0$ depending only on $\rho ,\zeta $. Note that (\ref{5.5})
entails that $\mathcal{D}_{\beta }$ is a compact (exponentially) attracting
set in $V^{-\theta _{2}}$ for $S_{NS}(t)$. By enlarging $R_{\beta }>0$ if
necessary, the claim follows easily.

\noindent \emph{Step (ii)}. One argues as follows: for $u_{0i}\in \mathcal{D}%
_{\beta },$ let $u_{i}=S_{NS}\left( t\right) u_{0i},$ $i=1,2$, be the
corresponding solutions of (\ref{MHD}), and decompose $u\left( t\right)
=u_{1}\left( t\right) -u_{2}\left( t\right) $ such that $u\left( t\right)
=v\left( t\right) +\omega \left( t\right) ,$ where $v\left( t\right) $ solves%
\begin{equation}
\left\{ 
\begin{array}{l}
\partial _{t}v+A_{0}v+B_{0}(u_{2},v)=0, \\ 
v\left( 0\right) =u_{01}-u_{02},%
\end{array}%
\right.  \label{5.6}
\end{equation}%
and $\omega \left( t\right) $ solves%
\begin{equation}
\left\{ 
\begin{array}{l}
\partial _{t}\omega +A_{0}\omega +B_{0}(u,u_{1})+B_{0}\left( u_{2},\omega
\right) =0, \\ 
\omega \left( 0\right) =0,%
\end{array}%
\right.  \label{5.7}
\end{equation}%
As in the proof of Theorem \ref{expo-thm2} one can show that%
\begin{equation*}
\left\Vert v\left( t\right) \right\Vert _{-\theta _{2}}^{2}\leq
e^{-c_{A_{0}}t}\left\Vert u_{01}-u_{02}\right\Vert _{-\theta _{2}}^{2},\text{
for all }t\geq 0
\end{equation*}%
and%
\begin{equation*}
\left\Vert \omega \left( t\right) \right\Vert _{\mathbb{\beta }}^{2}\leq
\lambda \left( t\right) \left\Vert u_{01}-u_{02}\right\Vert _{-\theta
_{2}}^{2},\text{ for all }t\geq 0,
\end{equation*}%
for some positive continuous function $\lambda :\mathbb{R}_{+}\rightarrow 
\mathbb{R}_{+},$ $\lambda \left( 0\right) >0$. These final estimates entail
that the mapping $S_{NS}\left( t\right) $ enjoys the smoothing property in
the sense of Theorem \ref{abstract} (Appendix), assumption (H4). We leave
the other (minor) details for the interested reader to check.
\end{proof}

As a consequence of this result, we also have the following.

\begin{corollary}
With the assumptions of Theorem \ref{expo-att-mhd}, the global attractor $%
\mathcal{A}_{reg}$ is bounded in $V^{\beta }$ and has finite fractal
dimension.
\end{corollary}

Finally, we state the analogue of Theorem \ref{expo-att-mhd} for problem (%
\ref{MHD}) in the case $\theta >0$. Note that in contrast to the case
studied in Section \ref{ss:attr-diss} (see the proof of Theorem \ref%
{expo-thm}) we can construct the exponential attractor directly on the set $%
V^{\beta },$ for some $\beta >-\theta _{2}.$

\begin{theorem}
\label{expo-att-mhd2}Let the assumptions of Theorems \ref{t:exist} and \ref%
{t:stab} be satisfied for some $\theta >0.$ In addition, for some $\beta \in
(-\theta _{2},\theta -\theta _{2}]$, let the following conditions hold.

(i) $b_{0}:V^{\beta }\times V^{\beta }\times V^{\theta -\beta }\rightarrow 
\mathbb{R}$ is bounded;

(ii) $g\in V^{\beta -\theta }$ is time independent.

Then $\left( S_{NS},V^{-\theta _{2}}\right) $ possesses an exponential
attractor $\widetilde{\mathcal{M}}\subset V^{\beta }$ which is bounded in $%
V^{\beta }$. Thus, by definition, we have

(a) $\widetilde{\mathcal{M}}$ is compact and semi-invariant with respect $%
S_{NS}\left( t\right) ,$ that is,%
\begin{equation*}
S_{NS}\left( t\right) (\widetilde{\mathcal{M}})\subseteq \widetilde{\mathcal{%
M}},\quad \forall \,t\geq 0.
\end{equation*}

(b) The fractal dimension $\dim _{F}(\widetilde{\mathcal{M}},V^{-\theta
_{2}})$ of $\widetilde{\mathcal{M}}$ is finite and an upper bound can be
computed explicitly.

(c) $\widetilde{\mathcal{M}}$ attracts exponentially fast any bounded subset 
$B$ of $V^{-\theta _{2}}$, that is, there exist a positive nondecreasing
function $Q$ and a constant $\rho >0$ such that 
\begin{equation*}
dist_{V^{-\theta _{2}}}(S_{NS}\left( t\right) B,\widetilde{\mathcal{M}})\leq
Q(\Vert B\Vert _{V^{-\theta _{2}}})e^{-\rho t},\quad \forall \,t\geq 0.
\end{equation*}%
Both $Q$ and $\rho $ can be explicitly calculated.
\end{theorem}

\begin{proof}
We sketch the proof by only showing that the semigroup $S_{NS}\left(
t\right) $ enjoys the smoothing property. By \cite[Theorem 5.1]{HLT}, we
recall that there exists a time $t_{4}>0$ such that%
\begin{equation}
\sup_{t\geq t_{4}}\left( \left\Vert u\left( t\right) \right\Vert _{\beta
}^{2}+\int_{t}^{t+1}\left\Vert u\left( s\right) \right\Vert _{\beta +\theta
}^{2}ds\right) \leq C,  \label{6.1}
\end{equation}%
for some positive constant $C$ independent of time and the initial data.
Setting $v=u_{1}-u_{2}$, and recalling that each $u_{i}$ is a solution of (%
\ref{MHD}), we observe that $v$ solves the following problem%
\begin{equation}
\partial _{t}v+A_{0}v+B_{0}\left( v,u_{1}\right) +B_{0}\left( u_{2},v\right)
=0.  \label{eq6.2}
\end{equation}%
First, owing to \cite[Theorem 3.5]{HLT}, for all $t\geq 0$ we have the
following estimate%
\begin{equation}
\left\Vert v\left( t\right) \right\Vert _{-\theta
_{2}}^{2}+\int_{0}^{t}\left\Vert v\left( s\right) \right\Vert _{\theta
-\theta _{2}}^{2}ds\leq \lambda \left( t\right) \left\Vert v\left( 0\right)
\right\Vert _{-\theta _{2}}^{2},  \label{eq6.2bis}
\end{equation}%
for some positive continuous function $\lambda \in C\left( \mathbb{R}_{+},%
\mathbb{R}_{+}\right) ,$ with $\lambda \left( 0\right) >0$. Pairing (\ref%
{eq6.2}) with $\Lambda ^{2\beta }v$, we infer%
\begin{equation}
\frac{1}{2}\frac{d}{dt}\Vert v\Vert _{\beta }^{2}+c_{A_{0}}\Vert v\Vert
_{\theta +\beta }^{2}\leq b_{0}(v,u_{1},\Lambda ^{2\beta
}v)+b_{0}(u_{2},v,\Lambda ^{2\beta }v).  \label{eq6.3}
\end{equation}%
The first two terms are bounded using the assumption (i) above. We have%
\begin{align}
b_{0}(v,u_{1},\Lambda ^{2\beta }v)& \leq C\delta ^{-1}\Vert u_{1}\Vert
_{\beta }^{2}\Vert v\Vert _{\beta }^{2}+\delta \left\Vert v\right\Vert
_{\beta +\theta }^{2},\text{ }  \label{eq6.4} \\
b_{0}(u_{2},v,\Lambda ^{2\beta }v)& \leq C\delta ^{-1}\Vert u_{2}\Vert
_{\beta }^{2}\Vert v\Vert _{\beta }^{2}+\delta \left\Vert v\right\Vert
_{\beta +\theta }^{2},  \notag
\end{align}%
for any $\delta >0$. Inserting the above estimates from (\ref{eq6.4}) into (%
\ref{eq6.3}), and choosing $\delta =c_{A_{0}}/2>0$ sufficiently small, we
deduce the following inequality%
\begin{equation}
\frac{d}{dt}\Vert v\left( t\right) \Vert _{\beta }^{2}\leq C\delta
^{-1}\left( \Vert u_{1}\Vert _{\beta }^{2}+\Vert u_{2}\Vert _{\beta
}^{2}\right) \Vert v\left( t\right) \Vert _{\beta }^{2}\leq C\Vert v\left(
t\right) \Vert _{\beta }^{2},\text{ for all }t\geq t_{4},  \label{eq6.5}
\end{equation}%
owing once more to the uniform estimate (\ref{6.1}). Multiplying now both
sides of this inequality by $\overline{t}:=t-t_{4}$ and integrating the
resulting relation over $\left( t_{4},t\right) ,$ thanks to (\ref{eq6.2bis})
we get%
\begin{equation*}
\Vert v\left( t\right) \Vert _{\beta }^{2}\leq C\left( \overline{t}+1\right)
\left( \overline{t}\right) ^{-1}e^{Ct}\left\Vert v\left( 0\right)
\right\Vert _{-\theta _{2}}^{2},
\end{equation*}%
for all $\overline{t}>0,$ which entails the required smoothing property for $%
S_{NS}\left( t\right) $. From this point on, the rest of the proof goes
essentially as in the proof of Theorem \ref{expo-thm}. We leave the details
to the interested reader.
\end{proof}

\begin{remark}
We can also argue as in the proof of Lemma \ref{l:smooth-compset} to deduce
further regularity properties for $\widetilde{\mathcal{M}}.$ Finally, note
that the assumptions of the last two theorems apply to 3D Navier-Stokes-$%
\alpha $ (NS-$\alpha $) equations, the 3D Leray-$\alpha $ models, the
modified 3D Leray-$\alpha $ models, the simplified 3D Bardina models,~the 3D
Navier-Stokes-Voight (NSV) equations, and many other models not explicitly
stated anywhere in the literature.
\end{remark}

\section{Remarks on a simplified Ericksen--Leslie model for liquid crystals}

\label{Leslie}

If we take $\phi $ as a vector, say $d\in \mathbb{R}^{n}$, then our
regularized system (\ref{e:pde}) can be used to model the motion of liquid
crystal flows in an $n$-dimensional compact Riemannian manifold $\Omega $,
with $n=2,3$:%
\begin{equation}
\left\{ 
\begin{array}{l}
\partial _{t}u+A_{0}u+(Mu\cdot \nabla )(Nu)+\chi \nabla (Mu)^{T}\cdot
(Nu)+\nabla p=-\varepsilon \text{div}(\nabla d\otimes \nabla d)+g, \\ 
\text{div}\left( u\right) =0, \\ 
\partial _{t}d+Nu\cdot \nabla d+\gamma \left( A_{1}d+f\left( d\right)
\right) =0, \\ 
u\left( 0\right) =u_{0}, \\ 
d\left( 0\right) =d_{0},%
\end{array}%
\right.  \label{EL-ls}
\end{equation}%
where, for instance, $A_{0}$, $A_{1}$, $M$, and $N$ are the same bounded
linear operators as (\ref{param}). In this context, the positive constants $%
\varepsilon ,$ $\gamma $ stand for the competition between kinetic energy
and potential energy, and the macroscopic elastic relaxation time (Deborah
number) for the molecular orientation field, respectively. Generally
speaking, the system (\ref{EL-ls}) may be viewed as a macroscopic continuum
description of the time evolutions of liquid crystal materials influenced by
both the flow field $u\left( x,t\right) $, and the microscopic orientational
configuration $d\left( x,t\right) $. The first and second equations combine
the laws describing the incompressible (regularized) flow of fluid with an
extra nonlinear coupling term, which is the induced elastic stress from the
elastic energy through the transport, represented by the third equation
which is a second-order parabolic equation with $f\left( d\right) =\nabla
_{d}F\left( d\right) $. Here $F\left( d\right) =(\left\vert d\right\vert
^{2}-1)^{2}$ is a potential function used to approximate the constraint $%
\left\vert d\right\vert =1$, see \cite{Er, Le}.

Problem (\ref{EL-ls}) with the following choice of operators%
\begin{equation}
A_{0}=P\left( -\nu \Delta \right) ,\text{ }M=N=I\text{ and }\chi =0
\label{paraspe}
\end{equation}%
has been investigated on a \textit{two-dimensional} compact Riemannian
manifold in \cite{Sh}, where the existence of a global attractor was also
proved. Global existence and regularity results in bounded domains $\Omega
\subset \mathbb{R}^{2}$ were also derived in \cite{LL} for the first time
(see also \cite{EGI, EGI2, HW}\ for related results in three dimensions).
The longtime behavior of the system (\ref{EL-ls}), under the hypothesis (\ref%
{paraspe}) and various boundary conditions, was also investigated recently
in \cite{Bos, GW, Wu} in the case of bounded domains $\Omega \subset \mathbb{%
R}^{n}$, $n=2,3$. In all these cases, a maximum principle as stated in
Proposition \ref{maxp} holds for the $d$-component of any weak solution $%
\left( u,d\right) $ of (\ref{EL-ls}). On account of this fact, all the
results on well-posedness, regularity and singular perturbations, the
existence of global and exponential attractors, and convergence to single
equilibria, as stated in Sections \ref{s:well}, \ref{s:pert} and \ref%
{s:longbehav} remain valid without any essential modifications for the
family of regularized problems (\ref{EL-ls}). To avoid redundancy, we
refrain from explicitly stating these results and their obvious proofs. In
particular, we recover the result on existence of global attractors for the
Lagrangian averaged liquid crystal equations, which consists of the
Navier-Stokes-$\alpha $-model coupled with the equation for the orientation
parameter $d$ from (\ref{EL-ls}). This case was reported in \cite[Theorem 7.1%
]{Sh}.

\section{Appendix}

\label{ss:app}

In this section, we include some supporting material on Gr\"{o}nwall-type
inequalities, Sobolev inequalities, some definitions and abstract results.
The first lemma is a slight generalization of the usual Gr\"{o}nwall-type
inequality \cite{T}; its proof is quite elementary and thus omitted.

\begin{lemma}
\label{Gineq}Let $\mathcal{E}:\mathbb{R}_{+}\rightarrow \mathbb{R}_{+}$ be
an absolutely continuous function satisfying 
\begin{equation*}
\frac{d}{dt}\mathcal{E}(t)+2\eta \mathcal{E}(t)\leq h(t)\mathcal{E}%
(t)+l\left( t\right) +k,
\end{equation*}%
where $\eta >0$, $k\geq 0$ and $\int_{s}^{t}h\left( \tau \right) d\tau \leq
\eta (t-s)+m$, for all $t\geq s\geq 0$ and some $m\in \mathbb{R}$, and $%
\int_{t}^{t+1}l\left( \tau \right) d\tau \leq \gamma <\infty $. Then, for
all $t\geq 0$, 
\begin{equation*}
\mathcal{E}(t)\leq \mathcal{E}(0)e^{m}e^{-\eta t}+\frac{2\gamma e^{m+\eta }}{%
e^{\eta }-1}+\frac{ke^{m}}{\eta }.
\end{equation*}
\end{lemma}

With $s,p\in \mathbb{R}_{+}$, let \thinspace $W^{s,p}$ be the standard
Sobolev space on an $n$-dimensional compact Riemannian manifold with $n\geq
2 $. The following result states the classical Gagliardo-Nirenberg-Sobolev
inequality (cf. \cite{BCL, He, Ni} and \cite{CM1, CM2}).

\begin{lemma}
\label{GNS}Let $0\leq k<m$ with $k,m\in \mathbb{N}$ and numbers $p,q,q\in %
\left[ 1,\infty \right] $ satisfy%
\begin{equation*}
k-\frac{n}{p}=\tau \left( m-\frac{n}{q}\right) -\left( 1-\tau \right) \frac{n%
}{r}.
\end{equation*}%
Then there exists a positive constant $C$ independent of $u$ such that%
\begin{equation*}
\left\Vert u\right\Vert _{W^{k,p}}\leq C\left\Vert u\right\Vert
_{W^{m,q}}^{\tau }\left\Vert u\right\Vert _{L^{r}}^{1-\tau },
\end{equation*}%
with $\tau \in \left[ \frac{k}{m},1\right] $ provided that $m-k-\frac{n}{r}%
\notin \mathbb{N}_{0},$ and $\tau =\frac{k}{m}$ provided that $m-k-\frac{n}{r%
}\in \mathbb{N}_{0}.$
\end{lemma}

We state here a standard result on pointwise multiplication of functions in
Sobolev spaces (see~\cite{Ma04a}; cf. also \cite{HLT}).

\begin{lemma}
\label{l:hole} Let $s$, $s_{1}$, and $s_{2}$ be real numbers satisfying 
\begin{equation*}
s_{1}+s_{2}\geq 0,\qquad min(s_{1},s_{2})\geq s,\qquad \text{and}\qquad
s_{1}+s_{2}-s>\frac{n}{2},
\end{equation*}%
where the strictness of the last two inequalities can be interchanged if $%
s\in \mathbb{N}_{0}$. Then, the pointwise multiplication of functions
extends uniquely to a continuous bilinear map 
\begin{equation*}
H^{s_{1}}\otimes H^{s_{2}}\rightarrow H^{s}.
\end{equation*}
\end{lemma}

We have the following definition of exponential attractor (also known as
inertial set).

\begin{definition}
\label{robust_expo}Let $\left( S\left( t\right) ,\mathcal{K}\right) $ be a
dynamical system on a given Banach space $\mathcal{K}$. A set $\mathcal{M}%
\subset \mathcal{K}$ is said to be an exponential attractor (also known as
inertial set) for the semigroup $S\left( t\right) $ provided that the
following statements hold:

$\left( i\right) $ The sets $\mathcal{M}$ are positively invariant with
respect to the semigroup $S\left( t\right) ,$ that is, $S\left( t\right) 
\mathcal{M}\subseteq \mathcal{M}$, for all $\,t\geq 0.$

$\left( ii\right) $ The fractal dimension of the sets $\mathcal{M}$ is
finite, that is, $\dim _{F}\left( \mathcal{M},\mathcal{K}\right) \leq
C<\infty ,$ where $C>0$ can be computed explicitly.

$\left( iii\right) $ Each $\mathcal{M}$ attracts exponentially any bounded
subset of $\mathcal{K},$ that is, there exist a positive constant $\rho $
and a monotone nonnegative function $Q,$ such that, for every bounded subset 
$B$ of $\mathcal{K},$ we have%
\begin{equation*}
dist_{\mathcal{K}}\left( S\left( t\right) B,\mathcal{M}\right) \leq
Q(\left\Vert B\right\Vert _{\mathcal{K}})e^{-\rho t},
\end{equation*}%
where $\displaystyle dist_{\mathcal{K}}\left( X,Y\right) :=\sup_{x\in
X}\inf_{y\in Y}\left\Vert x-y\right\Vert _{\mathcal{K}}$ is the Hausdorff
semidistance.
\end{definition}

We report the following basic abstract results (see \cite[Theorem 4.4]{GMPZ}%
, \cite{GG1, GG2}; cf. also \cite{MZ}) which are needed in order to prove
Theorem \ref{expo-thm} when $\theta >0$ and Theorem \ref{expo-thm2} in the
case $\theta =0$.

\begin{theorem}
\label{t4.5}Let ${\mathcal{X}}_{1}$ and ${\mathcal{X}}_{2}$ be two Banach
spaces such that ${\mathcal{X}_{2}}$ is compactly embedded in ${\mathcal{X}}%
_{1}.$ Let $X_{0}$ be a bounded subset of $\mathcal{X}_{2}$ and consider a
nonlinear map $\Sigma :X_{0}\rightarrow X_{0}$ satisfying the smoothing
property 
\begin{equation}
\left\Vert \Sigma \left( x_{1}\right) -\Sigma \left( x_{2}\right)
\right\Vert _{{\mathcal{X}}_{2}}\leq d\left\Vert x_{1}-x_{2}\right\Vert _{{%
\mathcal{X}_{1}}},  \label{3.56}
\end{equation}%
for all $x_{1},x_{2}\in X_{0},$ where $d>0$ depends on $X_{0}.$ Then the
discrete dynamical system $(X_{0},\Sigma ^{n})$ possesses a discrete
exponential attractor $\mathcal{E}_{M}^{\ast }\subset {\mathcal{X}_{2}},$
that is, a compact set in ${\mathcal{X}_{1}}$ with finite fractal dimension
such that 
\begin{equation}
\Sigma \left( \mathcal{E}_{M}^{\ast }\right) \subset \mathcal{E}_{M}^{\ast },
\label{3.57}
\end{equation}%
\begin{equation}
dist_{{\mathcal{X}_{1}}}\left( \Sigma ^{n}\left( X_{0}\right) ,\mathcal{E}%
_{M}^{\ast }\right) \leq d_{X}e^{-\rho _{\ast }n},\quad n\in {\mathbb{N}},
\label{3.58}
\end{equation}%
where $d_{X}$ and $\rho _{\ast }$ are positive constants independent of $n,$
with the former depending on $X_{0}.$
\end{theorem}

\begin{theorem}
\label{abstract}Let $\mathcal{K}$, $\mathcal{K}_{c}$ be two Banach spaces
such that $\mathcal{K}_{c}$ is compactly embedded in $\mathcal{K}$, and let $%
\left( S\left( t\right) ,\mathcal{K}\right) $ be a dynamical system. Assume
the following hypotheses hold:

\noindent (H1) There exists a bounded subset $\mathbb{B}\subset \mathcal{K}$
which is positively invariant for $S(t)$ and attracts any bounded set of $%
\mathcal{K}$ exponentially fast.

\noindent (H2) There exists a positive constant $C$ independent of time such
that%
\begin{equation*}
\left\Vert S\left( t\right) \varphi _{1}-S\left( t\right) \varphi
_{1}\right\Vert _{\mathcal{K}}\leq \rho \left( t\right) \left\Vert \varphi
_{1}-\varphi _{2}\right\Vert _{\mathcal{K}},
\end{equation*}%
for every $t\geq 0$, and every $\varphi _{1},\varphi _{2}\in \mathbb{B}$,
where $\rho :\mathbb{R}_{+}\rightarrow \mathbb{R}_{+}$ is some continuous
function with $\rho \left( 0\right) >0.$

\noindent (H3) There exist a positive constant $C,$ $\kappa \in \left(
0,1\right) $ and a time $t^{\ast }>0$ such that%
\begin{equation*}
\left\Vert S(t)\varphi _{0}-S(\tilde{t})\varphi _{0}\right\Vert _{\mathcal{K}%
}\leq C|t-\tilde{t}|^{\kappa },
\end{equation*}%
for all $t,\tilde{t}\in \left[ t^{\ast },2t^{\ast }\right] $ and any $%
\varphi _{0}\in \mathbb{B}.$

\noindent (H4) For every $\varphi _{01},\varphi _{02}\in \mathbb{B}$, $S(t)$
can be decomposed as follows:%
\begin{equation*}
S(t)\varphi _{01}-S(t)\varphi _{02}=D\left( t\right) \left( \varphi
_{01},\varphi _{02}\right) +N\left( t\right) \left( \varphi _{01},\varphi
_{02}\right)
\end{equation*}%
where, for all $t\geq 0$, we have%
\begin{equation*}
\left\{ 
\begin{array}{l}
\left\Vert D\left( t\right) \left( \varphi _{01},\varphi _{02}\right)
\right\Vert _{\mathcal{K}}\leq Ke^{-\varkappa t}\left\Vert \varphi
_{01}-\varphi _{02}\right\Vert _{\mathcal{K}}, \\ 
\left\Vert N\left( t\right) \left( \varphi _{01},\varphi _{02}\right)
\right\Vert _{\mathcal{K}_{c}}\leq \rho \left( t\right) \left\Vert \varphi
_{01}-\varphi _{02}\right\Vert _{\mathcal{K}},%
\end{array}%
\right.
\end{equation*}%
for some positive constants $\varkappa ,K$ independent of time, and some
positive continuous function $\rho :\mathbb{R}_{+}\rightarrow \mathbb{R}_{+}$%
, $\rho \left( 0\right) >0.$

Then, if (H1)-(H4) are satisfied, there exists an exponential attractors $%
\mathcal{M}$ for $\left( S\left( t\right) ,\mathcal{K}\right) $ in the sense
of Definition \ref{robust_expo}, (i)-(iii).
\end{theorem}

\end{document}